\documentclass[pdflatex]{article}
\usepackage{amsmath}
\usepackage{amssymb}
\usepackage{mathabx}
\usepackage{amsfonts}
\usepackage{amsthm}
\usepackage{bm}
\usepackage{lmodern}
\usepackage{fancybox}
\usepackage{latexsym}
\usepackage{bbm}

\usepackage{cases}
\DeclareFontFamily{U}{mathc}{}
\DeclareFontShape{U}{mathc}{m}{it}%
{<->s*[1.03] mathc10}{}

\DeclareMathAlphabet{\mathscr}{U}{mathc}{m}{it}

\usepackage{graphics}
\usepackage{enumerate}
\usepackage[numbers]{natbib}
\usepackage{mathscinet}
\usepackage{setspace}
\usepackage{graphicx}
\usepackage{ulem}
\usepackage{mathrsfs}

\usepackage[inline]{enumitem}

\usepackage{hyperref}
\hypersetup{colorlinks=true,linkcolor=blue,
citecolor=blue}

\newsavebox{\toy}
\savebox{\toy}{\framebox[0.65em]{\rule{0cm}{1ex}}}
\newcommand{\QED}{\usebox{\toy}\end{demo}}

\newcommand{\toybox}{\framebox[0.65em]{\rule{0cm}{0.65ex}}}

\makeatletter
\renewcommand{\th@plain}{%
  \let\thm@indent\indent
  \thm@headfont{\bfseries}%
  \thm@notefont{\normalfont}%
  \itshape
}
\makeatother

\theoremstyle{plain} 
\newtheorem{theorem}{Theorem}[section]
\newtheorem{lemma}[theorem]{Lemma}
\newtheorem{proposition}[theorem]{Proposition}
\newtheorem{corollary}[theorem]{Corollary}
\newtheorem{remark}[theorem]{Remark}

\newtheorem{definition}[theorem]{Definition}

\makeatletter
             \@addtoreset{equation}{section}
\makeatother

\usepackage{geometry}
\newcommand{\dis}{\displaystyle}

\newcommand{\F}{{\mathcal F}}
\usepackage{comment}

\newcommand{\N}{{\mathbb{N}}}
\newcommand{\Z}{{\mathbb{Z}}}

\newcommand{\R}{{\mathbb{R}}}

\newcommand{{\rd}}{\R^d}

\newcommand{\IP}{{\mathbb P}}
\newcommand{\IE}{\mathbb E}

\newcommand{\lan}{\langle}
\newcommand{\ran}{\rangle}

\renewcommand{\b}{\beta}

\newcommand{\e}{\varepsilon}

\newcommand{\dd}{\text{\rm d}}             

\newcommand{\sZ}{\mathsf{Z}}
\newcommand{\scZ}{\mathscr{Z}}
\newcommand{\ind}{\mathbbm{1}}

\newcommand{\lr}[1]{\left\lfloor{#1}\right\rfloor}
\newcommand{\bsm}{\begin{smallmatrix}}
\newcommand{\esm}{\end{smallmatrix}}

\makeatletter
\newcommand{\subsubsubsection}{\@startsection{paragraph}{4}{\z@}%
  {1.0\Cvs \@plus.5\Cdp \@minus.2\Cdp}%
  {.1\Cvs \@plus.3\Cdp}%
  {\reset@font\sffamily\normalsize}
}
\makeatother
\begin{document}
\author{Makoto Nakashima}
\title{Martingale measure associated with the critical $2d$ stochastic heat flow}
\date{}
\maketitle

\begin{abstract}
In \cite{CSZ23},  finite-dimensional convergence was established for rescaled random fields arising from directed polymers in  random
environment. The resulting scaling limit is known as the critical $2d$ stochastic heat flow (SHF).

In this paper, we develop a stochastic calculus framework for the critical $2d$ SHF. 
We first show that the critical $2d$ SHF admits a continuous measure-valued version. 
We then formulate a martingale problem for the $2d$ SHF, analogous to the one for super-Brownian motion. 
This allows us to construct an associated orthogonal martingale measure in the sense of Walsh \cite[Chapter 2]{Wal86}.

As an application, we establish a lower bound on the roughness of the
critical $2d$ SHF through its quadratic variation.
\end{abstract}

\vspace{1em}
\textbf{2020 Mathematics Subject Classification.}
Primary 60H17; Secondary 60H05, 60G44, 60G57.

\vspace{1em}\textbf{Keywords:} Critical $2d$ stochastic heat flow, drected polymers in random environment, martingale problem, martingale measure.

\section{Introduction and main results}

In this paper, we identify a martingale structure associated with the critical $2d$ stochastic heat flow (SHF), introduced by Caravenna, Sun, and Zygouras \cite{CSZ23}.
The critical $2d$ SHF is a random measure-valued flow indexed by two time parameters. 
In the present work, we mainly fix the initial time and regard it as a one-parameter measure-valued stochastic process. 
More precisely, for an initial density $\phi$, we write $\mathscr{Z}_t^{\vartheta,\phi}$ for the corresponding measure-valued process.

The critical $2d$ SHF arises as the scaling limit of two-dimensional directed polymers in random environment.
It may be viewed as a candidate solution, in a singular critical regime, to the two-dimensional multiplicative stochastic heat equation\begin{align}
  \partial_t u(t,x)=\frac{1}{2}\Delta u(t,x)+\beta u(t,x)\dot{\mathcal{W}}(t,x),\label{eq:SHE}\tag{SHE}
\end{align}
where $\beta\in\mathbb{R}$ and $\dot{\mathcal{W}}$ is the space-time white noise on $[0,\infty)\times \R^2$.
In two spatial dimensions, \eqref{eq:SHE} lies at a critical singular
regime. In particular, the standard theories developed for singular
SPDEs, such as the theory of regularity structures, do not apply
directly.

Caravenna, Sun, and Zygouras constructed the critical $2d$ SHF as the scaling limit of rescaled partition functions of directed polymers and proved convergence in the sense of finite-dimensional distributions. 
More recently, Tsai established process-level convergence of
rescaled solutions to mollified stochastic heat equations in the
vague topology \cite{Tsa24}.
These results provide the limiting object, but they do not by themselves yield a stochastic calculus for it.

A stochastic calculus framework is central to the analysis of
measure-valued processes and stochastic PDEs and provides a basis for
further developments. Our first result, Theorem \ref{thm:ContVer},
shows that the critical two-dimensional SHF admits a continuous
measure-valued version. We then investigate the martingale structure
of this version.
At a purely formal level suggested by the multiplicative stochastic
heat equation, for a suitable initial density $\phi$ and test function $\psi$, one expects that \begin{align*}
&\int_{\R^2}\mathscr{Z}^{\vartheta,\phi}(t,x)\psi(x)\dd x-\int_{\R^2}\phi(x)\psi(x)\dd x\\
&=\int_0^t \int_{\R^2}\frac{1}{2}\Delta \mathscr{Z}^{\vartheta,\phi}(s,x)\psi(x)\dd x\dd s+\int_0^t\int_{\R^2} \beta \mathscr{Z}^{\vartheta,\phi}(s,x)\psi(x)\dot{\mathcal{W}}(s,x)\dd s\dd x,
\end{align*}
so that the last term should define a martingale. 

Theorem \ref{thm:SHFMP} makes this formal decomposition rigorous by
formulating a martingale problem for the continuous measure-valued
version and identifying the quadratic variation of its martingale
part. In particular,
\[
\int_{\R^2}\psi(x)\mathscr Z_t^\vartheta(\dd x)
\]
is a semimartingale for every suitable test function $\psi$.

The resulting family of martingales is compatible across spatial
test functions. Theorem \ref{thm:martinalemeasure} shows that it defines
an orthogonal martingale measure in the sense of Walsh
\cite{Wal86}, thereby providing a stochastic calculus framework for
the critical two-dimensional SHF.

Finally, we use the quadratic-variation formula to study the local
roughness of the critical two-dimensional SHF.
Theorem \ref{thm:peaks} establishes the occurrence of large local
peaks, from which the lower bound on roughness in
Theorem \ref{thm:regular} follows.

\subsection{$\mathscr{Z}^\vartheta$ as a measure-valued process}

After fixing the initial time at $0$, for each initial condition $\phi\in C_c^+(\R^2)$ or $\phi\in C_b^+(\R^2)$, we consider the measure-valued process
\begin{align*}
\scZ^{\vartheta}(\phi,\dd y):=\left\{\scZ^{\vartheta}_t(\phi,\dd y)\right\}_{t\geq 0},
\end{align*}
where
\begin{align*}
\scZ^{\vartheta,\phi}_t(\dd y)=\scZ^{\vartheta}_t(\phi,\dd y):=\int_{\R^2} \phi(x)\scZ^\vartheta_{t}(\dd x,\dd y).
\end{align*}

We begin by recalling several basic facts about measure-valued processes; see \cite{Per02} for a standard reference on Dawson-Watanabe superdiffusions.

We equip the space $M_{LF}(\mathbb{R}^2)$ of locally finite measures on $\mathbb{R}^2$ with the topology of vague convergence, and the space $M_F(\mathbb{R}^2)$ of finite measures on $\mathbb{R}^2$ with the topology of weak convergence.

For a Polish space $E$, let $C(E):=C(\R_{+},E)$ denote the space of $E$-valued paths equipped with the topology of uniform convergence on compact sets.

Since  $M_{LF}(\R^2)$ and $M_F(\R^2)$ are  Polish spaces \cite[3.1.1 and 3.1.4]{Daw93},  it follows that $C(M_{LF}(\R^2))$ and $C(M_F(\R^2))$ are also Polish spaces.

For an $M_F(\R^2)$- or $M_{LF}(\R^2)$-valued process
$X=\{X_t\}_{t\geq0}$, let
\begin{align}
\mathcal F_t^{X,0}
:=
\sigma\left(X_s:0\leq s\leq t\right),
\qquad
\mathcal F_t^X
:=
\bigcap_{u>t}
\left(
\mathcal F_u^{X,0}\vee\mathcal N\label{eq:filtration}
\right),
\end{align}
where $\mathcal N$ denotes the collection of all
$\mathbb P$-null sets in $\mathcal F$.
Then $\{\mathcal F_t^X\}_{t\geq0}$ is the usual augmentation of the
natural filtration generated by $X$.

The following theorem establishes the existence of a continuous measure-valued version  of  $\scZ^{\vartheta,\phi}_t$.
\begin{theorem}\label{thm:ContVer}
\begin{enumerate}[label=(\arabic*)]
\item For each $\phi\in C_c^{+}(\R^2)$, there exists a continuous $M_F(\R^2)$-valued process $\scZ^{\vartheta,\phi}=\{\scZ_t^{\vartheta,\phi}\}_{t\geq 0}$ whose  finite-dimensional distributions coincide with those of the  critical $2d$ SHF.
\item For each $\phi\in C_b^{+}(\R^2)$, there exists a continuous $M_{LF}(\R^2)$-valued process $\scZ^{\vartheta,\phi}=\{\scZ_t^{\vartheta,\phi}\}_{t\geq 0}$ whose  finite-dimensional distributions coincide with those of the  critical $2d$ SHF.
\end{enumerate}

\end{theorem}

\begin{remark}
In \cite{Tsa24}, Tsai proved that the critical $2d$ stochastic heat flow ${\scZ}^{\vartheta}$ admits a continuous $M_{LF}(\R^2)$-valued version for each $\vartheta\in \R$. 
More precisely, he shows that the SHF $\{\mathcal{Z}^\vartheta_{s,t}\}_{-\infty<s\leq t<\infty}$ is continuous as a process taking values in $M_{LF}(\R^2\times\R^2)$.  
\end{remark}

We next formulate a martingale problem for the measure-valued process $\scZ^{\vartheta,\phi}$, in a form analogous  to the martingale problem for super-Brownian motion.

\begin{theorem}\label{thm:SHFMP}
Let $\phi\in C_c^+(\R^2)$, $\psi\in C_b^2(\R^2)$, and $\vartheta\in \R$. Define the continuous process $\mathscr{Z}_t^{\vartheta,\phi}(\psi):=\int \psi(x)\mathscr{Z}_t^{\vartheta,\phi}(\dd x)$. 

Define \begin{align}
\mathscr{M}_t^{\vartheta,\phi}(\psi):=\scZ^{\vartheta,\phi}_t(\psi)-\int_{\R^2} \phi(x)\psi(x)\dd x-\int_0^t \scZ^{\vartheta,\phi}_s\left(\frac{1}{2}\Delta\psi\right)\dd s\label{eq:Martingalepart}
\end{align}
for $t\geq 0$. 
Then, $\mathscr{M}_t^{\vartheta,\phi}(\psi)$ is a continuous $\left\{\F_t^{\scZ^{\vartheta}}\right\}_{t\geq 0}$-martingale with
$\mathscr{M}_0^{\vartheta,\phi}(\psi)=0$ and quadratic variation
\begin{align}
&\left\langle \mathscr{M}^{\vartheta,\phi}(\psi) \right\rangle_t=\lim_{\e\to 0}\frac{4\pi}{\log \frac{1}{\e}}\int_0^t \int_{\R^2} \left(\scZ^{\vartheta,\phi}_u(p_\e(\cdot-z))\right)^2\psi(z)^2\dd z\dd u,\label{eq:QuadConv}
\end{align}
where 
\begin{align}
\scZ^{\vartheta,\phi}_t(p_\e(\cdot-x))=\int_{\R^2}\int_{\R^2} \phi(y)p_\e(z-x)\scZ^{\vartheta}_t(\dd y,\dd z).\label{eq:QVSBHF}
\end{align}
where the convergence in \eqref{eq:QuadConv} holds uniformly on
compact time intervals in probability.

The same conclusions hold for $\mathscr{Z}^{\vartheta,\phi}(\psi)$ with $\phi\in C_b^+(\R^2)$ and $\psi\in C_c^2(\R^2)$. 
\end{theorem}

\begin{remark}
We remark that the martingale problem in Theorem \ref{thm:SHFMP} does not uniquely
characterize the law of the process.
Indeed, if the superscripts $\vartheta$ is suppressed in \eqref{eq:Martingalepart}-\eqref{eq:QVSBHF}, then the formulation of  the martingale problem no longer contains the parameter $\vartheta$. 
However,  the distributions of the corresponding processes depend on $\vartheta\in\R$. In particular, one can verify  \begin{align*}
\mathscr{Z}_t^{\vartheta,\phi}(1)=\int_{\R^2}\phi(x)\dd x+\mathscr{M}_t^{\vartheta,\phi}(1)\overset{d}{\ne }\int_{\R^2}\phi(x)\dd x+\mathscr{M}_t^{\vartheta',\phi}(1)=\mathscr{Z}_t^{\vartheta',\phi}(1)
\end{align*}
for $\vartheta\ne \vartheta'$ by looking at their variances (see \eqref{eq:varianceZ}). In particular, we know that the deterministic process \begin{align*}
\mathscr{Z}_t^{-\infty,\phi}(\psi):=\int_{\R^2}\dd x\int_{\R^2}\dd y\phi(x)p_t(x,y)\psi(y)
\end{align*}
solves the same martingale problem, with identically zero martingale
part and quadratic variation.
Thus, the martingale problem \eqref{eq:Martingalepart}-\eqref{eq:QVSBHF} admits a family of solutions  $\{\mathscr{Z}^{\vartheta,\phi}(\cdot)\}_{\vartheta\in\mathbb{R}\cup\{-\infty\}}$. 

In \cite{Tsa24}, the uniqueness of $\mathscr{Z}^\vartheta$ relies on an explicit representation of its fourth moment. This suggests that the martingale problem may become well-posed if additional information, such as the  fourth moment structure of the flow $\mathscr{Z}^\vartheta(\dd x,\dd y)$,  is imposed.
However, this remains an open problem.
\end{remark}

\begin{remark}
  In \cite{Che26}, Chen formulated a martingale problem for the critical $2d$ SHF that differs from the one in Theorem \ref{thm:SHFMP}. He also provides an explicit expression for the quadratic variation process, which depends on the parameter $\vartheta\in \R$.
\end{remark}

\begin{remark}
Let $u$ be a formal solution to the stochastic heat equation \begin{align*}
\partial_t u=\frac{1}{2}\Delta u+\beta u \dot{W}, \quad u(0,x)=\phi(x)
\end{align*}
and suppose that it admits a continuous function-valued
realization.
Then, the quadratic variation  process $\langle \int_{\R^2}u(\cdot,x)\psi(x)\dd x\rangle_t$ is formally given by 
\begin{align}
\int_0^t \int_{\R^2}\beta^2 u^2(s,x)\psi(x)^2\dd x\dd s.\label{eq:QVprocForm}
\end{align}
This reflects the effect of space-time white noise. Theorem \ref{thm:SHFMP} provides a rigorous formulation of \eqref{eq:QVprocForm}.

It is known  \cite{CSZ24,CSZ25} that $\scZ^{\vartheta,\phi}_t(\dd y)$ is singular with respect to the Lebesgue measure, 
and hence $\scZ_t^{\vartheta,\phi}(p_\e(\cdot-x))$ diverges at some points. Consequently, a renormalization factor $\frac{1}{\log \frac{1}{\e}}$ is required in \eqref{eq:QuadConv}.
\end{remark}

\begin{remark}
For standard super-Brownian motions $\{X_t\}_{t\geq 0}$,  the quadratic variation $\langle X(\psi)\rangle_t$ is given by \begin{align*}
\langle X(\psi)\rangle_t=\gamma\int_0^t X_s(\psi^2)\dd s, \quad t\geq 0
\end{align*}
for some $\gamma>0$. It is explicitly determined by the process $\{X_t\}_{t\geq 0}$.

By contrast, in the case of super-Brownian motion with a single-point catalyst in  $d=1$, 
the quadratic variation can be expressed in terms of a density field obtained as a limit of $\int_0^t X_s(p_\e(\cdot-y))\dd s$ for $d=1$ \cite{DF94}. 
\end{remark}

By the definition of the cross-variation of martingales, we obtain the following.
\begin{corollary}\label{cor:CrossVar}
Let $\phi\in C_c^+(\R^2)$, $\psi_1,\psi_2\in C_b^2(\R^2)$, and $\vartheta\in \R$. Let  $\mathscr{M}^{\vartheta,\phi}_t(\psi_i)$ ($i=1,2$) be the martingales defined in \eqref{eq:Martingalepart} corresponding to $\psi_1,\psi_2$.   
Then, we have 
\begin{align}
&\left\langle \mathscr{M}^{\vartheta,\phi}(\psi_1),\mathscr{M}^{\vartheta,\phi}(\psi_2) \right\rangle_t=\lim_{\e\to 0}\frac{4\pi}{\log \frac{1}{\e}}\int_0^t \int_{\R^2} \left(\scZ^{\vartheta,\phi}_u(p_\e(\cdot-z))\right)^2\psi_1(z)\psi_2(z)\dd z\dd u,\label{eq:CrossConv}
\end{align}
where the convergence in \eqref{eq:CrossConv} is  locally uniform in probability.

The same conclusion holds when
$\phi\in C_b^+(\R^2)$ and
$\psi_1,\psi_2\in C_c^2(\R^2)$.

\end{corollary}

\begin{remark}\label{rem:Ito}
Theorem \ref{thm:SHFMP} provides the semimartingale representation of $\mathscr{Z}^{\vartheta,\phi}(\psi)$. We can therefore  apply It\^{o}'s formula to  $f\left(\mathscr{Z}_t^{\vartheta,\phi}(\psi)\right)$ for  $f\in C_b^2(\R)$. 

This yields \begin{align*}
f\left(\mathscr{Z}_t^{\vartheta,\phi}(\psi)\right)&=f\left(\mathscr{Z}_0^{\vartheta,\phi}(\psi)\right)+\int_0^t f'\left(\mathscr{Z}_s^{\vartheta,\phi}(\psi)\right)\mathscr{Z}_s^{\vartheta,\phi}\left(\frac{1}{2}\Delta\psi\right)\dd s\\
&+\int_0^t f'\left(\mathscr{Z}_s^{\vartheta,\phi}(\psi)\right)\dd\mathscr{M}_s^{\vartheta,\phi}\left(\psi\right)+\frac{1}{2}\int_0^t f''\left(\mathscr{Z}_s^{\vartheta,\phi}(\psi)\right)\dd \left\langle \mathscr{M}^{\vartheta,\phi}\left(\psi\right)\right\rangle_s.
\end{align*}
However, it is unclear whether \begin{align*}
&\int_0^t f''\left(\mathscr{Z}_s^{\vartheta,\phi}(\psi)\right)\dd \left\langle \mathscr{M}^{\vartheta,\phi}\left(\psi\right)\right\rangle_s\\
&\hspace{3em}=\lim_{\e\to 0}\frac{4\pi}{\log \frac{1}{\e}}\int_0^t\int_{\R^2}f''\left(\mathscr{Z}_s^{\vartheta,\phi}(\psi)\right)\left(\mathscr{Z}_s^{\vartheta,\phi}(p_\e(\cdot-z))\right)^2\psi(z)^2\dd z\dd s
\end{align*}
 holds in probability. 
Furthermore, it remains open whether  $\left\langle \mathscr{M}^{\vartheta,\phi}\left(\psi\right)\right\rangle_s$ is absolutely continuous in time with respect to the Lebesgue measure.

\end{remark}

\subsection{A martingale measure and stochastic integral associated with $\mathscr{Z}^\vartheta$}

We now focus on the martingale part $\mathscr{M}^{\vartheta,\phi}(\psi)$ defined in \eqref{eq:Martingalepart} for $\vartheta\in \R$ and $\phi\in C_c^+(\R^2)$.

Let $\mathcal{M}_c^2$ denote the set of continuous   $\mathcal{F}^{\mathscr{Z}^{\vartheta,\phi}}_t$-martingales $M$ satisfying $M_0=0$ and $\IE[M^2_t]<\infty$ for each $t>0$.   

Then,  $\mathscr{M}^{\vartheta,\phi}$ defines a linear map from  $C_b^2(\R^2)$ into $\mathcal{M}_c^2$, associating each $\psi$ with the process $\mathscr{M}^{\vartheta,\phi}_t(\psi)$.

The following theorem establishes a linear extension of this map. 

\begin{theorem}\label{thm:martinalemeasure}
Let $\phi\in C_c^+(\R^2)$,  and $\vartheta\in\R$. Then, 
the linear map
\[
\mathscr M^{\vartheta,\phi}:
C_b^2(\R^2)\longrightarrow\mathcal M_c^2
\]
admits a unique linear extension
\[
\mathscr M^{\vartheta,\phi}:
\mathcal B_b(\R^2)\longrightarrow\mathcal M_c^2
\]
such that, for every $\psi\in\mathcal B_b(\R^2)$ and $t>0$,
\begin{align}
\left\langle \mathscr{M}^{\vartheta,\phi}(\psi)\right\rangle_t=\lim_{\e\to 0}\frac{4\pi}{\log\frac{1}{\e}}\int_0^t \int_{\R^2}\left(\mathscr{Z}_s^{\vartheta,\phi}\left(p_\e(\cdot-z)\right)\right)^2 \psi(z)^2\dd z\dd s,\label{eq:QuadGenlim}
\end{align}
 where  $\mathcal{B}_b(\R^2)$ is the set of bounded Borel measurable functions and  the convergence in \eqref{eq:QuadGenlim} holds  uniformly on compact time intervals in probability.
 
 The same result holds for $\phi\in C_b^+(\R^2)$ after replacing $C_b^2(\R^2)$ and $\mathcal{B}_b(\R^2)$   by $C_c^2(\R^2)$ and $\mathcal{B}_c(\R^2)$, respectively, 
 where $\mathcal{B}_c(\R^2)$ denotes the set of bounded Borel measurable functions with compact support. 

\end{theorem}

For indicator functions $\psi(x):=\mathbf{1}_A(x)$ with $A\in \mathcal{B}(\mathbb{R}^2)$, we write $\mathscr{M}^{\vartheta,\phi}_t(A):=\mathscr{M}^{\vartheta,\phi}_t(\mathbf{1}_A)$.

\begin{remark}
Theorem \ref{thm:martinalemeasure} together with Corollary \ref{cor:CrossVar} implies that $\mathscr{M}^{\vartheta,\phi}(\cdot)$ defines an orthogonal martingale measure in the sense of Walsh \cite[Chapter 2]{Wal86}. 
Indeed, for each $A\in \mathcal{B}(\R^2)$, $\mathscr{M}^{\vartheta,\phi}_t(A)$ is a continuous martingale, and  if $A,B\in \mathcal{B}(\R^2)$ with $A\cap B=\emptyset$, then \begin{align*}
\left\langle \mathscr{M}^{\vartheta,\phi}(A),\mathscr{M}^{\vartheta,\phi}(B) \right\rangle_t=0\quad \text{a.s.~for all $t>0$.}
\end{align*}
By the general theory in \cite[Corollary 2.8]{Wal86}, there exists a random $\sigma$-finite measure $\nu^{\vartheta,\phi}(\dd x,\dd s)$ on $\R^2\times [0,\infty)$ such that \begin{align*}
\nu^{\vartheta,\phi}((A\cap B)\times [0,t])=\left\langle \mathscr{M}^{\vartheta,\phi}(A),\mathscr{M}^{\vartheta,\phi}(B) \right\rangle_t
\end{align*} 
for $A,B\in \mathcal{B}(\R^2)$ with $E\left[\left\langle \mathscr{M}^{\vartheta,\phi}(C) \right\rangle_t\right]<\infty$ ($C=A,B$).

Moreover, stochastic integrals with respect to this martingale measure can then
be defined following \cite[Chapter 2]{Wal86}. 
\end{remark}

\subsection{A lower bound on the roughness of
$\mathscr{Z}^{\vartheta,\phi}_t(\dd x)$}\label{sec:reg}
As an application of the quadratic-variation formula and
Theorem \ref{thm:martinalemeasure}, we derive a lower bound on the
roughness of $\mathscr{Z}^{\vartheta,\phi}_t(\dd x)$.

It is known that $\mathscr{Z}^{\vartheta,\phi}_t(\dd x)$ is singular with respect to the Lebesgue measure \cite[Theorem 1.1]{CSZ25}. 
More precisely, the averaged density $\frac{1}{\pi\delta^2}\mathscr{Z}^{\vartheta,\phi}_t\left(B(x,\delta)\right)$ converges to $0$ for Lebesgue-a.e.~$x\in \R^2$ as $\delta\to 0$.

As a consequence, on a set of positive probability, there exists $(t,x)\in (0,\infty)\times \R^2$ such that \begin{align*}
\frac{1}{\pi \delta^2}\mathscr{Z}^{\vartheta,\phi}_t\left(B(x,\delta)\right)\to \infty \text{ as $\delta\to 0$.}
\end{align*}
Furthermore, \cite[Theorem 1.5]{CSZ25} shows that for any $t>0$ and $\vartheta\in \R$, $\mathscr{Z}_t^{\vartheta,\phi}(\dd x)$ belongs to  $\mathcal{C}^{0-}:=\bigcap_{\e>0}\mathcal{C}^{-\e}$, 
where $\mathcal{C}^{-\e}$ denotes the negative Besov-H\"older space of order $-\e$ in the sense of \cite[Definition 2.1]{FM17}. 

We now derive estimates describing  ``the peaks" of $\mathscr{Z}^{\vartheta,\phi}_t\left(\dd x\right)$.

Fix $f\in C_c^+(\R^2)$ with $\int_{\R^2}f(x)\dd x=1$. We define the random set \begin{align*}
\mathcal{T}_{s,t}(A,f,\lambda,\e):=\left\{(r,x):r\in [s,t],\, x\in A,\, \mathscr{Z}_{r}^{\vartheta,\phi}(f_\e(\cdot-x))\geq \lambda\log \frac{1}{\e}\right\},
\end{align*} 
 \text{for $\lambda>0$, $\e\in (0,1)$}, where $A\in \mathcal{B}(\R^2)$ is a Borel set of finite Lebesgue measure and $f_\e(x)=\frac{1}{\e}f\left(\frac{x}{\sqrt{\e}}\right)$.

\begin{theorem}\label{thm:peaks}
Fix $\vartheta\in \R$ and $\phi\in C_c^+(\R^2)$ with $\phi\not\equiv 0$. Let $A\subset\R^2$ be a nonempty open set of finite Lebesgue measure and $0\leq s<t<\infty$. 
\begin{enumerate}[label=(\arabic*)]
\item\label{item:LARGE1} We have \begin{align*}
\varlimsup_{\lambda\to \infty}\varlimsup_{\e\to 0}\iint_{\mathcal{T}_{s,t}(A,f,\lambda,\e)}\mathscr{Z}^{\vartheta,\phi}_u(f_\e(\cdot-x))\dd x\dd u=0\quad \text{a.s.}
\end{align*}
\item\label{item:LARGE2} 
There exists a deterministic decreasing sequence  $\{\e_n\}_{n\geq 1}$ with $\e_n\to 0$ such that \begin{align*}
P\left(\bigcup_{\lambda>0}\varlimsup_{n\to \infty}\mathcal{T}_{s,t}(A,f,\lambda,\e_n)\ne  \emptyset\right)>0.
\end{align*}
 
\end{enumerate}
\end{theorem}

\begin{remark}
Theorem \ref{thm:peaks} \ref{item:LARGE1} shows that the contribution
of the high-level set
$\mathcal T_{s,t}(A,f,\lambda,\e)$ to the integrated mollified mass
vanishes as $\lambda\to\infty$.

On the other hand,
Theorem \ref{thm:peaks} \ref{item:LARGE2} shows that peaks of order
$\log(1/\e)$ nevertheless occur with positive probability along a
deterministic sequence $\e_n\to0$.

Several recent works have further investigated the intermittency and
small-scale mass distribution of the critical two-dimensional SHF.
Liu and Zygouras \cite{LZ26} determined the leading logarithmic order
of the moments of the mass assigned to shrinking balls.
Ganguly and Nam obtained sharp high-moment and upper-tail estimates
\cite{GN25}, and subsequently established a precise form of
logarithmic intermittency and mass concentration at small spatial
scales \cite{GN26}. Related comparison inequalities for positive and
negative moments of locally averaged SHF were recently obtained by
Liu and Xie \cite{LX26}.

\end{remark}

Thus, it is natural to expect that $\mathscr{Z}^{\vartheta,\phi}_t(\dd x)$ does not belong to the ``logarithmic H\"older space" $\mathscr{C}^{0,b}$ for any $-1<b<0$.

We say $\xi\in \mathscr{S}'(\R^2)$ belongs to $\mathscr{C}^{s,b}$ for $s<0$ and $b\in\R$, or $s=0$ and $b<0$  if 
\begin{itemize}
\item  $\xi$ belongs to the dual of $\mathcal{C}^r$ with $r=\begin{cases}-\lfloor s\rfloor &\text{if }b\geq 0\\ -\lfloor s\rfloor  +1 \quad &\text{if }b<0\end{cases}$, and
\item for any compact set $K\subset\mathbb{R}^2$ \begin{align*}
\|\xi\|_{s,b;K}:=\sup_{\e \in (0,1]}\sup_{\bsm f\in \mathcal{C}^r,\\ \mathrm{supp}f\subset B(0,1),\|f\|_r\leq 1\esm}\sup_{x\in K}\left|\frac{(1+\log \frac{1}{\e})^b}{\e^s }\left\langle \xi,\frac{1}{\e^2}f\left(\frac{\cdot-x}{\e}\right)\right\rangle \right|<\infty.
\end{align*}
\end{itemize}
Here, we denote by $\lr{s}$ the greatest integer $n\leq s$ for $s\in \R$.
\begin{remark}
For the definition of $\mathcal{C}^r$ and $\|\cdot\|_r$, see \cite[Definition 3.7]{Hai14}. It is shown   in \cite{Hai14} that  $\mathcal{C}^\alpha$ is closely related to the Besov space $B^\alpha_{\infty,\infty}$.
The space $\mathscr{C}^{s,b}$ can be viewed as a modification of $\mathcal{C}^\alpha$. 
To our knowledge, its relationship to generalized Besov spaces (see e.g.~\cite{Mou01,FL04,Alm05}) has not been investigated. \end{remark}

We have the following result.
\begin{theorem}\label{thm:regular}
Fix $\vartheta\in \R$, and $\phi\in C_c^+(\R^2)$ with $\phi\not\equiv 0$.
Then,   $\mathscr Z^{\vartheta,\phi}$ is not a continuous $ \mathscr{C}^{0,b}$-valued process with respect to the uniform-on-compact topology for any $-1<b<0$.
\end{theorem}

\begin{proof}[Proof of Theorem \ref{thm:regular}]
Choose a nonempty bounded open set $A\subset\R^2$ and a compact set
$K$ containing $\overline A$.

Fix $f\in C_c^\infty(\mathbb R^2)$ 
 with $f\in \mathscr{C}^\gamma$, $\mathrm{supp}(f)\subset B(0,1)$ and $\|f\|_\gamma\leq 1$.

By \ref{item:LARGE2} in Theorem \ref{thm:peaks}, for $s=0$ and $b>-1$,  \begin{align*}
&\sup_{0\leq u\leq T}\sup_{\e \in (0,1]}\sup_{x\in K}\left|{\left(1+\log \frac{1}{\e}\right)^b}\mathscr{Z}_u^{\vartheta,\phi}\left(\frac{1}{\e^2}f\left(\frac{\cdot-x}{\e}\right)\right) \right|=\infty 
\end{align*}
with positive probability for any $T>0$.  
\end{proof}

\begin{remark}

Theorem \ref{thm:regular} provides a lower bound on the roughness of $\mathscr{Z}_t^{\vartheta,\phi}(\dd x)$. The corresponding upper bound  still remains open.
 It is natural to expect that $\mathscr{Z}_t^{\vartheta,\phi}(\dd x)$ belongs to even more singular spaces.
\end{remark}

\begin{remark}
It is of interest  to identify a random set $\Lambda$ on which $\mathscr{Z}_t^{\vartheta,\phi}(\dd x)$ concentrates, i.e.~$\mathscr{Z}_t^{\vartheta,\phi}(A)=\mathscr{Z}_t^{\vartheta,\phi}(\R^2)$ almost surely.  

In the Gaussian multiplicative chaos theory, the GMC measure concentrates on $\gamma$-thick points of the underlying log-correlated Gaussian field, 
and the Hausdorff dimension of its support is given by $2-\frac{\gamma^2}{2}$ \cite{Kah85} (See also \cite{RV14}). 

\end{remark}

\subsection{Organization of the paper}

The paper is organized as follows.

Section \ref{sec:prelim} introduces the discrete approximation and reviews the moment estimates of $\mathscr{Z}^\vartheta$, based on \cite{CSZ19b,CSZ24}.

Section  \ref{sec:proofcontQuad}  gives the proofs of Theorems \ref{thm:ContVer} and \ref{thm:SHFMP}, relying on technical estimates proved in the appendices.

Sections \ref{sec:Quadconv} and \ref{sec:peaks} provide the proofs of Theorem \ref{thm:martinalemeasure} and Theorem \ref{thm:peaks}, respectively.

Appendices \ref{sec:ChaosExpansion}--\ref{sec:momentsbound} contain the remaining technical arguments needed to complete the proofs of Theorem \ref{thm:ContVer} and Theorem \ref{thm:SHFMP}.

\section{Preliminaries}\label{sec:prelim}

\subsection{Setting and known results}
In this paper, we consider the discrete model introduced in \cite{CSZ23}.

Let $\{\omega_{n,x}\}_{n\in \Z,x\in \Z^2}$ be i.i.d.~random variables on a probability space $(\Omega, \mathcal{F},\mathbb{P})$, satisfying \begin{align}
&\mathbb{E}[\omega_{n,x}]=0,\quad \mathbb{E}[\omega_{n,x}^2]=1,\quad \lambda(\beta) :=\log \mathbb{E}\left[e^{\beta \omega_{n,x}}\right]<\infty\quad  \text{for small }\beta>0.\label{eq:omegacond}
\end{align}

Now, we introduce the random fields according to \cite{CSZ23,CSZ24}.

Let $\{S_n\}_{n\geq 0}$ be an irreducible, symmetric, and aperiodic random walk on $\Z^2$ whose increment $\xi:=S_1-S_0$ has mean $0$ and identity covariance matrix. Let $P$ and $E$ denote probability and expectation for $S$.

We assume that $\xi$ has finite support, i.e.~there exists a finite set $\Sigma\subset \Z^2$ such that $P(\xi\in \Sigma)=1$.

Let us define the point-to-point partition functions for directed polymers in random environment as follows:\begin{align}
&Z_{N}^\beta(x,y):=E\left[\left.\exp\left(\sum_{i=1}^{N}\left(\beta \omega_{i,S_i}-\lambda(\beta)\right)\right)\ind _{S_N=y}\right|S_0=x\right]		\notag\\
&\overline{Z}_{N}^{\beta}(x,y):=E\left[\left.\exp\left(\sum_{i=1}^{N-1}\left(\beta \omega_{i,S_i}-\lambda(\beta)\right)\right)\ind_{S_N=y}\right|S_0=x\right]		\notag
\end{align}
for $N\geq 0$ and $x,y\in \Z^2$, $\beta\geq 0$, where we use the convention $\sum_{n=1}^{0}\{\dots\}:=0$.

Note that \begin{align*}
\mathbb{E}\left[Z_{N}^\beta(x,y)\right]=\mathbb{E}\left[\overline{Z}_{N}^\beta(x,y)\right]=P(S_N=y|S_0=x)=q_{N}(y-x),
\end{align*}
where  we denote by $q_n(x)$ the transition probability  of the underlying random walk starting at $0$, i.e.\begin{align*}
q_n(x):=P(S_n=x|S_0=0)
\end{align*}
for $x\in\Z^2$, $n\geq 0$.

For $x\in\R^2$, let $\lfloor x\rfloor$ denote a nearest point in
$\Z^2$, with ties resolved lexicographically. When no confusion can
arise, we suppress integer-part or nearest-lattice-point notation.

Fix $N\geq 1$. We define the rescaled measure-valued random measures by \begin{align*}
\sZ^{\b}_{N;t}:=\frac{1}{N}\sum_{x,y\in\Z^2}Z^{\b}_{\lr{Nt}}\left(x,y\right)\delta_{\frac{x}{\sqrt{N}}}\otimes \delta_{\frac{y}{\sqrt{N}}},
\end{align*}
for $t\geq 0$. Note that for $\phi\in C_c(\R^2)$, $\psi\in C_b(\R^2)$, and $t\geq 0$ \begin{align}
&\sZ_{N;t}^{\beta}(\phi,\psi)\notag\\
&:=\frac{1}{N}\sum_{x\in\Z^2}\sum_{y\in \Z^2}\phi\left(\frac{x}{\sqrt{N}}\right)Z_{\lr{Nt}}^\beta(x,y)\psi\left(\frac{y}{\sqrt{N}}\right)\notag\\
&=\frac{1}{N}\sum_{x\in \Z^2}E\left[\left.\psi\left(\frac{S_{\lr{Nt}}}{\sqrt{N}}\right)\exp\left( \sum_{i=1}^{\lr{Nt}}\left(\beta\omega_{i,S_{i}}-\lambda(\beta)\right)\right)\right|S_0=x\right]\phi\left(\frac{x}{\sqrt{N}}\right).\notag
\end{align}
We introduce the rescaled random flows $\overline{\mathsf{Z}}_{N;t}^\beta$ and $\overline{\mathsf{Z}}_{N;t}^\beta(\phi,\psi)$ by replacing $Z^\beta$ with its modified partition function  $\overline{Z}^\beta$.

We note that for each $\phi\in C_c(\R^2)$ and $\psi\in C_b(\R^2)$ and $t\geq 0$, \begin{align}
\overline{Z}^\beta_{N;t}(\phi,\psi)\overset{d}{=}\overline{Z}^\beta_{N;t}(\psi,\phi).\label{eq:reversal}
\end{align} 
We emphasize that this identity in distribution holds at each fixed time, but not at the process level.

\begin{remark}
The sequence of random measure-valued processes $\sZ_{N}^\beta=\left\{\sZ_{N;t}^{\beta}\right\}_{t\geq 0}$ is essentially the same as the one considered in \cite{CSZ23,CSZ24}. 
The differences are as follows:\begin{enumerate}
\item In \cite{CSZ23,CSZ24}, the authors used $\overline{\sZ}_N^{\beta}$ instead of $\sZ_N^\beta$.
\item In \cite{CSZ23,CSZ24}, random measures are absolutely continuous with respect to the Lebesgue measure by replacing  Dirac measures with uniform measures on squares. 
\end{enumerate} 
However, these differences are negligible for convergence. 
\end{remark}

To obtain a non-trivial limit, we rescale the disorder strength by setting  $\beta=\beta_N$ appropriately.

Let $S'$ be an independent copy of $S$ and $R_N$ be the expectation of the number of collisions of  $S$ and $S'$ up to time $N$: \begin{align*}
R_N:=\sum_{n=1}^N P(S_n=S_n')=\sum_{n=1}^N\sum_{x\in \Z^2}q_n(x)^2=\sum_{n=1}^{N}q_{2n}(0).
\end{align*}
Then, the local limit theorem below gives the asymptotic behavior \begin{align*}
R_N= \frac{\log N}{4\pi}(1+o(1)).\notag
\end{align*} 

\begin{theorem}[{Local limit theorem \cite[P7.9]{Spi76}, \cite[Theorem 2.3.5 and Theorem 2.3.11]{LL10}}]
Let $q_n$ be the transition probability defined above. Then,  we have 
 \begin{align}
q_n(x)=p_n(x)+O\left(\frac{1}{n^2}\right)\qquad\text{uniformly in }x\in\Z^2\notag
\end{align}
for $n\in \mathbb{N}$, where \begin{align*}
p_t(y-x)=p_t(x,y)=\frac{1}{2\pi t}\exp\left(-\frac{|y-x|^2}{2t}\right)\qquad t>0, x,y\in\R^2
\end{align*}  
is the Gaussian heat kernel on $\R^2$ with mean $0$ and covariance matrix $t I$.
\end{theorem}

We choose  $\beta=\beta_N$ such that \begin{align}
\sigma_N^2:=e^{\lambda(2\beta_N)-2\lambda(\beta_N)}-1=\frac{1}{R_N}\left(1+\frac{\vartheta}{\log N}(1+o(1))\right) \quad\text{as $N\to \infty$}\label{eq:betan}
\end{align}
for some $\vartheta\in \R$.

\begin{theorem}\label{thm:SHF} \cite[Theorem 1.1, Theorem 9.1]{CSZ23} \cite[Theorem 6.1]{CSZ24}
The family of random measures $\mathsf{Z}^{\b_N}_N=\left\{\mathsf{Z}^{\b_N}_{N;t}\right\}_{t\geq 0}$ in $M_{LF}(\R^2)\times M_{F}(\R^2)$ converges in finite-dimensional distributions to a unique limit \begin{align*}
\mathscr{Z}^\vartheta=\{\scZ^\vartheta_{t}(\dd x,\dd y)\}_{t\geq 0}.
\end{align*}
Moreover, the distribution of $\mathscr{Z}^\vartheta$ is independent of the choice of $\{\omega_{n,x}\}_{n,x}$.
\end{theorem}

\begin{remark}
We will see that $\mathbb{E}\left[\left(\mathsf{Z}^{\b_N}_{N;t}(\phi,\psi)-\overline{\mathsf{Z}}^{\b_N}_{N;t}(\phi,\psi)\right)^2\right]\to 0$ in Remark \ref{rem:ZZbar}, so $\left\{\overline{\mathsf{Z}}^{\b_N}_{N;t}\right\}$ also converges to $\mathscr{Z}^\vartheta$. 
\end{remark}

\begin{remark}
 To be precise, in \cite[Theorem 1.1, Theorem 9.1]{CSZ23} \cite[Theorem 6.1]{CSZ24}, 
 Caravenna, Sun, and Zygouras considered the flow version $\left\{\sZ^\beta_{N;s,t}\right\}_{-\infty<s\leq t<\infty}$ 
and showed that $\left\{\sZ^\beta_{N;s,t}\right\}_{-\infty<s\leq t<\infty}$ in $M_{LF}(\R^2)\times M_F(\R^2)$ converges in finite-dimensional distributions to a unique limit $\left\{\mathscr{Z}^\vartheta_{s,t}\right\}_{-\infty<s\leq t<\infty}$, i.e.~for $-\infty< s_i\leq t_i<\infty$,  $\phi_i\in C_c(\R^2)$, and $\psi_i\in C_b(\R^2)$ ($i=1,\dots,k$)
\begin{align*}
\left\{ \mathsf{Z}^{\b_N}_{N;s_i,t_i}(\phi_i,\psi_i)\right\}_{i=1,\dots,k} \Rightarrow\left\{\scZ^\vartheta_{s_i,t_i}(\phi_i,\psi_i)\right\}_{i=1,\dots,k},
\end{align*}
where \begin{align*}
\scZ^\vartheta_{s,t}(\phi,\psi)&=\int\int \phi(x)\psi(y)\scZ^\vartheta_{s,t}(\dd x,\dd y).
\end{align*}
\end{remark}

\subsection{Variance and its limit}
In this subsection, we recall an explicit formula for the variance of $\mathsf{Z}_{N;t}^{\beta_N,\phi}(\psi)$ and its limit as $N\to \infty$, which yields
the variance of $\mathscr{Z}_t^{\vartheta,\phi}(\psi)$. 

We note that \begin{align*}
\mathbb{E}\left[\mathsf{Z}_{N;t}^{\beta_N,\phi}(\psi)\right]=\frac{1}{N}\sum_{x\in \Z^2}\sum_{y\in \Z^2}\phi\left(\frac{x}{\sqrt{N}}\right)q_{Nt}(y-x)\psi\left(\frac{y}{\sqrt{N}}\right)
\end{align*}
and that as $N\to \infty$,
\begin{align}
 \mathbb{E}\left[\mathsf{Z}_{N;t}^{\beta_N,\phi}(\psi)\right]\to \int_{\R^2}\phi(x)\dd x \int_{\R^2}p_{t}(y-x)\psi(y)\dd y.\notag
\end{align}

By the standard second-moment method for directed polymers in random environment, we obtain \begin{align*}
&\mathbb{E}\left[\mathsf{Z}_{N;t}^{\beta_N,\phi}(\psi)^2\right]\\
&=\frac{1}{N^2}\sum_{x,x'\in \Z^2}\phi\left(\frac{x}{\sqrt{N}}\right)\phi\left(\frac{x'}{\sqrt{N}}\right)E^{\otimes}_{x,x'}\left[e^{(\lambda(2\beta_N)-2\lambda(\b_N))\sum_{i=1}^{Nt} \mathbf{1}\{S_i=S_i'\}}\psi\left(\frac{S_{Nt}}{\sqrt{N}}\right)\psi\left(\frac{S'_{Nt}}{\sqrt{N}}\right)\right],
\end{align*}
where $(S,P_x)$ and $(S',P_{x'})$ are independent random walks starting at $x$ and $x'$, respectively, and with increments distributed as $\xi$.
Observe that the exponential weight in the second moment expands as \begin{align}
e^{(\lambda(2\beta_N)-2\lambda(\b_N))\sum_{i=1}^{Nt} \mathbf{1}\{S_i=S_i'\}}&=\prod_{i=1}^{Nt} \left(1+\sigma_N^2 \mathbf{1}\{S_i=S_i'\}\right)\notag\\
&=1+\sum_{k=1}^{Nt}\sigma_N^{2k} \sum_{1\leq i_1<\dots<i_k\leq Nt}\prod_{j=1}^k   \mathbf{1}\{S_{i_j}=S_{i_j}'\}.\notag
\end{align}
Thus, it follows that \begin{align*}
&\mathrm{Var}(\mathsf{Z}^{\b_N,\phi}_{N;t}(\psi))\\
&=\frac{1}{N^2}\sum_{x_0,x_0'\in \Z^2}\phi\left(\frac{x_0}{\sqrt{N}}\right)\phi\left(\frac{x_0'}{\sqrt{N}}\right)\\
&\qquad \sum_{k=1}^{Nt}\sigma_N^{2k} \sum_{0=i_0<i_1<\dots<i_k\leq Nt}\sum_{x_1,\dots,x_k\in \Z^2}\sum_{y,y'\in\Z^2}q_{i_{0},i_1}(x_{0},x_{1})q_{i_{0},i_1}(x_{0}',x_{1})\\
&\qquad\qquad  \prod_{j=2}^k  \left(q_{i_{j-1},i_j}(x_{j-1},x_{j})^2\right)q_{i_k,Nt}(x_k,y)q_{i_k,Nt}(x_k,y')\psi\left(\frac{y}{\sqrt{N}}\right)\psi\left(\frac{y'}{\sqrt{N}}\right),
\end{align*}
with the convention $\prod_{k=2}^1\cdots=1$, where $q_{i,j}(x,y)=q_{j-i}(y-x)$ for $0\leq i\leq j$ and $x,y\in\Z^2$.

To analyze this quantity further, we use the weighted local renewal functions introduced in \cite{CSZ19a} with slight modifications to the original definitions: For each $N\geq 1$, 
\begin{align}
U_N(0,x)&=\sigma_N^2\delta_{x,0}=\sigma_N^2\mathbf{1}_{\{x=0\}},\label{eq:UN}\\
U_N(n,x)&=\sigma_N^4 q_{0,n}(0,x)^2\notag\\
&+\sum_{k\geq 1} (\sigma_N^2)^{k+2}\hspace{-1em}\sum_{\bsm 0<n_1<\dots<n_k<n\\  x_1,\dots,x_k\in \Z^2\esm} q_{0,n_1}(0,x_1)^2\left(\prod_{j=2}^{k}q_{n_{j-1},n_{j}}(x_{j-1},x_j)^2\right)q_{n_k,n}(x_k,x)^2\qquad \notag
\intertext{for $n\geq 1$, and}
U_N(n)&=\sum_{x\in\Z^2} U_N(n,x).\notag
\end{align}
Then, \begin{align}
\mathrm{Var}(\mathsf{Z}^{\b_N,\phi}_{N;t}(\psi))&=\frac{1}{N^2}\sum_{x_0,x_0'\in \Z^2}\sum_{\bsm 1\leq i\leq j\leq {Nt}\\ x,y,y',z\in \Z^2\esm}\phi\left(\frac{x_0}{\sqrt{N}}\right)\phi\left(\frac{x_0'}{\sqrt{N}}\right)q_{i}(x_0,x)q_{i}(x_0',x)\label{eq:varUN}\\
&\cdot U_{N}(j-i,z-x)q_{j,{Nt}}(z,y)q_{j,{Nt}}(z,y')\psi\left(\frac{y}{\sqrt{N}}\right)\psi\left(\frac{y'}{\sqrt{N}}\right).\notag
\end{align}

Thus,  $U_N(n,x)$ plays a key role in controlling the modulus of continuity of $\{\mathsf{Z}_{N;t}^{\b_N,\phi}(\psi)\}_{t\geq 0}$.  

We review some known results on $U_N$ and $\mathsf{Z}_{N;t}^{\b_N,\phi}(\psi)$ obtained in \cite{CSZ19a,CSZ19b}. 
We introduce the following weighted Green's function of the Dickman subordinator:\begin{align*}
&G_{\vartheta}(t):=\int_0^\infty \frac{e^{(\vartheta-\gamma)s}st^{s-1}}{\Gamma(s+1)}\dd s 
\intertext{for $t\in (0,\infty)$, together with its space-time version}
&G_{\vartheta}(t,x):=G_{\vartheta}(t)p_{\frac{t}{2}}(x),
\end{align*}
for $t\in (0,\infty)$ and $x\in \R^2$, where $\gamma$ denotes the Euler--Mascheroni constant.

\begin{theorem}[{\cite[Theorem 1.4 and 2.3]{CSZ19a}, \cite[Proposition 8.4]{CSZ24}}]\label{thm:UN}
Suppose $\vartheta\in \R$. Then,  \begin{align}
U_N(n)=\frac{\sigma_N^2\log N}{N}G_{\vartheta}\left(\frac{n}{N}\right)\left(1+o(1)\right),\quad \textrm{uniformly for $\delta N\leq n\leq TN$,}\label{eq:UNapproxG}
\end{align}
for  $0<\delta<T<\infty$. Also,
\begin{align}
&U_N(n)\leq C\frac{\sigma_N^2\log N}{N}G_{\vartheta}\left(\frac{n}{N}\right),\quad 1\leq n\leq TN\label{eq:UNboundG}
\intertext{and}
&U_N(n,x)=\frac{\sigma_N^2\log N}{N^2}G_\vartheta\left(\frac{n}{N},\frac{x}{\sqrt{N}}\right)\left(1+o(1)\right),\label{eq:UNxapproxG}\\
&\hspace{10em}\textrm{uniformly for $\delta N\leq n\leq TN$ and $|x|\leq   \frac{\sqrt{N}}{\delta}$.}\notag
\end{align}
\end{theorem}

In this paper, we often omit the subscript $N$ and denote by \begin{align*}
U_{m,n}(x,y)=U_{N}(n-m,y-x),\qquad  U_{m,n}=\sum_{y\in \Z^2}U_{m,n}(0,y)
\end{align*}
for $0\leq m\leq n<\infty$ and $x,y\in \Z^2$.

\begin{proposition}[{\cite[Proposition 1.6]{CSZ19a}, \cite[Proposition 2.3]{CSZ19b}, \cite[Proposition 8.2]{CSZ24}}]\label{pro:Gbound}
For fixed $\vartheta\in \R$, \begin{align*}
G_\vartheta(t)=\frac{1}{t\left(\log \frac{1}{t}\right)^2}+\frac{2\vartheta}{t\left(\log \frac{1}{t}\right)^3}+O\left(\frac{1}{t\left(\log \frac{1}{t}\right)^4}\right),\quad \text{as $t\to 0$}. 
\end{align*}
Also, for $T>0$, there exists $c_{\vartheta,T}\in (0,\infty)$ such that \begin{align}
G_\vartheta(t)\leq \widehat{G}_{\vartheta,T}(t):=\frac{c_{\vartheta,T}}{t\left(\log \frac{e^2T}{t}\right)^2},\quad t\in (0,T].\label{eq:GGhat}
\end{align}
In particular, $\widehat{G}_{\vartheta,T}$ is decreasing in $t\in (0,T]$.
\end{proposition}

Combining  \eqref{eq:UNboundG} and \eqref{eq:GGhat}, we obtain that \begin{align}
U_N(n)\leq \frac{1}{N}\frac{C_{\vartheta,T}}{\frac{n}{N}\left(\log \frac{e^2TN}{n}\right)^2}\quad \text{for }1\leq n\leq TN.\label{eq:UNGbound}
\end{align}

The following is a modification of Theorem 1.2 in \cite{CSZ19b} and Theorem 6.1 in \cite{CSZ24}.
\begin{theorem}\label{thm:variancelimit}
Let $\vartheta\in \R$. Suppose $\b_N$ satisfies \eqref{eq:betan}. Then, we have for every $t\geq0$,  each $\phi\in C_c(\R^2)$ and $\psi\in C_b(\R^2)$
\begin{align}
&\mathrm{Var}(\mathscr{Z}^{\vartheta,\phi}_t(\psi))=\lim_{N\to \infty}\mathrm{Var}(\mathsf{Z}_{N;t}^{\b_N,\phi}(\psi))\label{eq:varianceZ}\\
&=4\pi \int_{0<u<v<t}\dd u\dd v\int_{\R^2\times \R^2}\dd x\dd y\Phi_u(x)^2 G_{\vartheta}(v-u,y-x)\Psi_{t-v}(y)^2,\notag
\end{align}
where  we set $\Phi_s(x)=\int_{\R^2}\phi(y)p_{s}(x-y)\dd y$ and $\Psi_s(x)=\int_{\R^2}\psi(y)p_{s}(x-y)\dd y$ for $s>0$ and $x\in \R^2$.
\end{theorem}

The higher moments of $\mathscr{Z}^{\vartheta,\phi}_t(\psi)$ are given explicitly in \cite[Theorem 1.1]{GQT21} and \cite[Theorem 9.6]{CSZ24}. 
To state the formula, we introduce some notation:
Let $h\geq 2$ be an integer. 
For $I=\{i,j\}$ ($1\leq i<j\leq h$), we define \begin{align*}
(\R^2)^h_I:=\{\bm{x}=(x_1,\dots,x_h)\in (\R^2)^h:x_i=x_j\}
\end{align*}
which is identified with $(\R^2)^{h-1}$. We denote by $\idotsint_{(\R^2)^h_I}f(\bm{x}) \dd \bm{x}_I$ the integral of an integrable function on $(\R^2)^h_I$ with respect to  the Lebesgue measure.

We define \begin{align*}
G^{I}_{\vartheta,t}(\bm{x},\bm{y}):=\left(\prod_{\ell\in \{1,\dots,h\}\backslash I}p_t(y_\ell-x_\ell)\right)G_\vartheta(t,y_i-x_i)
\end{align*}
for $\bm{x},\bm{y}\in (\R^{2})^h_I, t>0$, and $I=\{i,j\}$. Also, we define \begin{align*}
&\mathscr{Q}^{I,J}_t(\bm{y},\bm{x})=\prod_{i=1}^hp_t(x_i-y_i)\qquad \bm{y}\in(\R^2)^h_I, \bm{x}\in (\R^2)^h_J,\\
&\mathscr{Q}^{*,J}_t(\bm{y},\bm{x})=\prod_{i=1}^hp_t(x_i-y_i)\qquad \bm{y}\in(\R^2)^h, \bm{x}\in (\R^2)^h_J,\\
&\mathscr{Q}^{I,*}_t(\bm{y},\bm{x})=\prod_{i=1}^hp_t(x_i-y_i)\qquad \bm{y}\in(\R^2)^h_I, \bm{x}\in (\R^2)^h
\end{align*}
for  $t>0$, and $I,J$ with $|I|=|J|=2$.

\begin{theorem}\label{thm:highermoments}
Fix $\vartheta\in\R$. Let $h\geq 2$ be an integer. Then, for $\phi_1,\dots,\phi_h\in C_c(\R^2)$, $\psi_1,\dots,\psi_h\in C_b(\R^2)$, and $t\geq 0$ \begin{align}
&E\left[\prod_{i=1}^h \mathscr{Z}^{\vartheta,\phi_i}_{t}(\psi_i)\right]=\idotsint\limits_{(\R^2)^{h}\times (\R^2)^{h}}\phi^{\otimes h}(\bm{z})\mathscr{K}^{(h)}_{t}(\bm{z},\bm{w})\psi^{\otimes h}(\bm{w})\dd \bm{z}\dd \bm{w},\label{eq:highermomentsrepre}
\end{align}
where  we write $\bm{x}=(x_1,\dots,x_h)\in (\R^2)^h$, 
$\phi^{\otimes h}(\bm{x}):=\prod_{i=1}^h \phi_i(x_i)$, 
$\psi^{\otimes h}(\bm{x}):=\prod_{i=1}^h \psi_i(x_i)$, and \begin{align*}
&\mathscr{K}_t^{(h)}(\bm{z},\bm{w}):=\prod_{i=1}^h p_t(w_i-z_i)\\
&+\sum_{m=1}^\infty (4\pi)^m \sum_{\bsm I_1,\dots,I_m\subset \{1,\dots,h\}\\ |I_\ell|=2,I_\ell\ne I_{\ell+1}\esm}\quad\idotsint\limits_{0<a_1<b_1<\dots<a_m<b_m<t}\dd\bm{a}\dd\bm{b}\idotsint\limits_{\bsm \bm{x}^{(\ell)},\bm{y}^{(\ell)}\in (\R^2)^{h}_{I_\ell}\\ \text{for }\ell=1,\dots,m\esm }\prod_{\ell=1}^m \dd \bm{x}_{I_\ell}\dd \bm{y}_{I_\ell}\\
&\mathscr{Q}^{*,I_1}_{a_1}(\bm{z},\bm{x}^{(1)})G^{I_1}_{\vartheta,b_1-a_1}(\bm{x}^{(1)},\bm{y}^{(1)})
\\
&\hspace{1em}\left(\prod_{\ell=2}^m \mathscr{Q}_{a_\ell-b_{\ell-1}}^{I_{\ell-1},I_\ell}(\bm{y}^{(\ell-1)},\bm{x}^{(\ell)})G^{I_\ell}_{\vartheta,b_\ell-a_\ell}(\bm{x}^{(\ell)},\bm{y}^{(\ell)})\right)\mathscr{Q}^{I_m,*}_{t-b_m}(\bm{y}^{(m)},\bm{w}).
\end{align*}

\end{theorem} 

These explicit moment computations will serve as the basis for the convergence analysis in subsequent sections.

\section{Proofs of Theorem \ref{thm:ContVer} and Theorem \ref{thm:SHFMP}}\label{sec:proofcontQuad}

In this section, we prove Theorem \ref{thm:ContVer} and \ref{thm:SHFMP}, which establish  continuity and a martingale problem for the critical $2d$ stochastic heat flow.

Our approach closely follows the framework for measure-valued martingale problems developed in \cite[Chapter 4]{Per02}, which provides a general method to characterize superprocesses via their generator and associated martingale problem. In particular, we adapt Perkins' method to the setting of the critical $2d$ stochastic heat flow, with suitable modifications.

\subsection{$\mathsf{Z}_{N;\cdot}^{\phi}$ as a measure-valued process}\label{subsec:measureval}

Fix $\phi\in C_c^+(\R^2)$ and $\psi\in C_b(\R^2)$. 
We note that the argument in this subsection can be adapted to $\mathsf{Z}_{N;\cdot}^\phi$ as $M_{LF}(\R^2)$-valued process, 
that is to the case $\phi\in C_b^+(\R^2)$ and $\psi\in C_c(\R^2)$.

By the universality of the convergence result recalled in Theorem \ref{thm:SHF}, 
we may assume that  $\{\omega_{n,x}\}_{n\in \Z,x\in \Z^2}$ are i.i.d.~random variables satisfying
\begin{align}
\IP(\omega_{n,x}=1)=\IP(\omega_{n,x}=-1)=\frac{1}{2}.\label{eq:Bernoulli}
\end{align}
Indeed, the limiting process $\{\scZ_{t}^{\vartheta}(\phi,\psi)\}$ does not depend on  the specific  distribution of $\{\omega_{n,x}\}$, as long as the condition \eqref{eq:omegacond} is satisfied. 
In this case, one can verify that 
$\left\{\xi_{n,x}^{\beta_N}:=e^{\beta_N\omega_{n,x}-\lambda(\beta_N)}-1\right\}_{(n,x)\in \N\times \Z^2}$ are i.i.d.~random variables 
with $\IP\left(\xi_{n,x}^{\beta_N}=\sigma_{N} \right)=\IP\left(\xi_{n,x}^{\beta_N}=-\sigma_{N}\right)=\frac{1}{2}$. 
In the Bernoulli case, the quantity $\sigma_N$ defined through \eqref{eq:betan} satisfies $\sigma_N=\tanh(\beta_N)$.
In particular, we have
\begin{align}
&\IE\left[\left(\xi^{\b_N}_{i,x}\right)^{2m-1}\right]=0,\quad \IE\left[\left(\xi^{\b_N}_{i,x}\right)^{2m}\right]=\sigma_{N}^{2m}\label{eq:omegamoment}
\end{align}
for every integer $m\geq1$. 

\eqref{eq:omegamoment} will help us estimate the chaos expansions of moments (see Section \ref{sec:ChaosExpansion}), but it is not crucial.

We fix $\beta_N$ as in \eqref{eq:betan}.  For notation simplicity, we write 
\begin{align*}
\sZ^{N,\phi}_{t}(\dd y)&:=\sZ^{\b_N}_{N;t}(\phi,\dd y),\quad \sZ^{N,\phi}_{t}(\psi):=\sZ_{N;t}^{\b_N}(\phi,\psi).
\end{align*}

We now regard $\{\sZ_{t}^{N,\phi}(\psi)\}_{t\geq 0}$ as  a discrete semimartingale as in the construction of super-Brownian motion from a critical branching Brownian motion \cite[II.~4]{Per02}.

Let $\{\mathcal{F}_n\}$ be the filtration generated by $\{\omega_{i,x}:x\in \Z^2,0\leq i\leq n\}$, and set $\overline{\mathcal{F}}_{t}^{N}:=\mathcal{F}_{Nt}$ for $t=\frac{k}{N}$ ($k\in \mathbb{N}_0$). We denote $
f_N(x)=f\left(\frac{x}{\sqrt{N}}\right) $
for $f\in C(\R^2)$.

Then, we decompose the increment into a drift term and a martingale term as   \begin{align*}
&\mathsf{Z}_{\frac{k+1}{N}}^{N,\phi}(\psi)-\mathsf{Z}_{\frac{k}{N}}^{N,\phi}(\psi)\\
&=\frac{1}{N}\sum_{x\in \Z^2} \phi_N(x)\\
&\times \left(\sum_{\tilde{y}\in \Z^2}E_x\left[e_{k+1,S_{k+1}}\prod_{i=1}^{k}e_{i,S_i}\mathbf{1}_{S_{k+1}=\tilde{y}}\right]\psi_N(\tilde{y}) -\sum_{y\in \Z^2}E_x\left[\prod_{i=1}^{k}e_{i,S_i}\mathbf{1}_{S_{k}=y}\right]\psi_N(y)\right)\\
&=\frac{1}{N}\sum_{x\in \Z^2} \phi_N(x)\sum_{\tilde{y}\in \Z^2} E_x\left[\prod_{i=1}^k e_{i,S_i}\mathbf{1}_{S_{k+1}=\tilde{y}}\right] \psi_N(\tilde{y})\left(e_{k+1,\tilde{y}}-1\right)\\
&+\frac{1}{N}\sum_{x\in \Z^2} \phi_N(x)\sum_{y\in \Z^2} Z_{k}^{\beta_N} (x,y)\left(\sum_{\tilde{y}\in \Z^2}q_1(y,\tilde{y})\psi_N(\tilde{y})-\psi_N(y)\right).
\end{align*}
where 
we set $e_{n,x}=e^{\beta_N\omega_{n,x}-\lambda(\b_N)}=\xi_{n,x}^{\b_N}+1$ for $n\geq 0$ and $x\in \Z^2$.

Let $\Delta_N$ denote the discrete Laplacian, defined by \begin{align*}
\Delta_N \psi_N\left(x\right)=N\left(\sum_{y\in \Z^2}q_1(x,y)\psi_N\left(y\right)-\psi_N\left(x\right)\right)
\end{align*}
and we write  \begin{align*}
\delta M_{k}^{N,\phi}(\psi)=\frac{1}{N}\sum_{\tilde{y}\in \Z^2} \overline{Z}_{N;k+1}^{\phi}(\tilde{y}) \psi_N(\tilde{y})\left(e_{k+1,\tilde{y}}-1\right),
\end{align*}
where 
we set \begin{align*}
  \overline{Z}^{\phi}_{N;n}(y)=E\left[\left.\phi_N\left(S_{n}\right)\exp\left( \sum_{i=1}^{n-1}\left(\beta_N\omega_{i,S_{n-i}}-\lambda(\beta_N)\right)\right)\right|S_0=y\right].
\end{align*}
We note that $E\left[\delta M_{k}^{N,\phi}(\psi)|\mathcal{F}_{k}\right]=0$ since $\overline Z_{N;k+1}^{\phi}(y)$ is $\mathcal{ F}_k$-measurable and
$\mathbb E[e_{k+1,y}-1]=0$.

Thus, we have \begin{align*}
&\mathsf{Z}_{\frac{k+1}{N}}^{N,\phi}(\psi)-\mathsf{Z}_{\frac{k}{N}}^{N,\phi}(\psi)=\frac{1}{N}\mathsf{Z}_{\frac{k}{N}}^{N,\phi} \left(\Delta_N \psi_N\right)+\delta M_{k}^{N,\phi} (\psi)
\end{align*}
and hence \begin{align}
\mathsf{Z}^{N,\phi}_{\frac{n}{N}}(\psi)
&=\mathsf{Z}_{0}^{N,\phi}(\psi)+\frac{1}{N}\sum_{k=0}^{n-1} \mathsf{Z}_{\frac{k}{N}}^{N,\phi} \left(\Delta_N \psi_N\right)+\sum_{k=0}^{n-1}\delta M_{k}^{N,\phi}(\psi)\label{eq:semimart}\\
&=\mathsf{Z}_{0}^{N,\phi}(\psi)+\int_{0}^\frac{n}{N} \mathsf{Z}_{t}^{N,\phi} \left(\Delta_N \psi_N\right)\dd t+\sum_{k=0}^{n-1}\delta M_{k}^{N,\phi}(\psi).\notag
\end{align}
We remark that \begin{align*}
M_{n}^{N,\phi}(\psi):=\sum_{k=0}^{n-1}\delta M_{k}^{N,\phi}(\psi)
\end{align*}
is an $\mathcal{F}_n$-martingale with $M_{0}^{N,\phi}(\psi)=0$ and its quadratic variation \begin{align*}
&\left\langle M^{N,\phi}(\psi)\right\rangle_{n}
=\frac{\sigma_N^2}{N^2}\sum_{k=1}^{n}\sum_{y\in \Z^2}\overline{Z}_{N;k}^{\phi}(y)^2\psi_N\left(y\right)^2.
\end{align*}

In particular, \begin{align}
\left\langle M^{N,\phi}(\psi)\right\rangle_{Nt}&=\frac{\sigma_N^2}{N^2}\sum_{k=1}^{Nt}\sum_{y\in \Z^2}{\overline{Z}}_{N;k}^{\phi}(y)^2\psi_N(y)^2
=\int_0^{t}\dd s {\frac{\sigma_N^2}{N}\sum_{y\in \Z^2}{\overline{Z}}_{N;Ns}^{\phi}(y)^2\psi_N(y)^2}.\label{eq:QVprocess}
\end{align}

\begin{remark}\label{rem:ZZbar}
The same argument above shows that \begin{align*}
\IE\left[\left.\mathsf{Z}_{N,t}^{\beta_N}(\phi,\psi)\right|\mathcal{F}_{\lr{Nt}-1}\right]=\overline{\mathsf{Z}}_{N,t}^{\beta_N}(\phi,\psi)
\end{align*}
and hence \begin{align*}
&\IE\left[\left(\mathsf{Z}_{N,t}^{\beta_N}(\phi,\psi)-\overline{\mathsf{Z}}_{N,t}^{\beta_N}(\phi,\psi)\right)^2\right]=\IE\left[\frac{\sigma_N^2}{N^2}\sum_{y\in \Z^2}\overline{Z}_{N,Nt}^\phi(y)^2\psi_N(y)^2\right]\to 0
\end{align*}
for each $t\geq 0$, where the convergence follows from Theorem \ref{thm:UN} and Proposition \ref{pro:Gbound}. 
\end{remark}

\subsection{Continuity of $\scZ^{\vartheta,\phi}$}\label{subsec:Cont}

To prove Theorem \ref{thm:ContVer}, we use several key lemmas from \cite{Per02}, \cite{Kal17} and \cite{Mit83}, which provide a general framework for establishing the continuity of measure-valued processes via $C$-tightness.

Let $D(E):=D(\R_+,E)$ denote the space of c\`adl\`ag $E$-valued paths equipped with the Skorokhod $J_1$-topology for a Polish space $E$.  We note that $D(E)$ is also a Polish space.

\begin{definition}
Let $E$ be a Polish space. The collection of processes $\{X^\alpha:\alpha\in  I\}$ with paths in $D(E)$ is said to be $C$-relatively compact in $D(E)$ if and only if it is relatively compact in $D(E)$ and all weak limit points are a.s.~continuous. 
\end{definition}

\begin{definition}
\begin{enumerate}[label=(\arabic*)] 
\item Let $E$ be a Polish space. We say that $D\subset C_b(E)$ is separating if and only if for any $\mu,\nu\in M_{F}(E)$, $\mu(\phi)=\nu(\phi)$ for all $\phi\in D$ implies $\mu=\nu$. 

\item Let $E$ be a Polish space. We say that $D\subset C_c(E)$ is separating if and only if for any  $\mu,\nu\in M_{LF}(E)$, $\mu(\phi)=\nu(\phi)$ for all $\phi\in D$ implies $\mu=\nu$. 
\end{enumerate}
\end{definition}

Theorem \ref{thm:ContVer} follows once we verify the following conditions.

\begin{theorem}\label{thm:Crel}{\cite{Mit83,Per02}}
\begin{enumerate}[label=(\arabic*)]
\item\label{item:Crel1} Let $E$ be   a Polish space and $D\subset C_b(E)$ be a separating class in $C_b(E)$ containing $1$. 
A sequence of c\`adl\`ag $M_{F}(E)$-valued processes $\{X^N\}$ is $C$-relatively compact in $D(M_{F}(E))$ if and only if 
the following conditions hold:
\begin{enumerate}[label=$(\text{$1$-}\mathrm{\roman*})$]
\item \label{item:CrelMF1} For all $\e>0$, $T>0$, there exists a compact set $K=K_{\e,T}$ in $E$ such that \begin{align*}
\sup_N P\left(\sup_{t\leq T}X_t^N\left(K_{\e,T}^c\right)>\e\right)<\e.
\end{align*} 
\item \label{item:CrelMF2} 
for all $\phi\in D$, $\left\{X^N(\phi)\right\}$ is  $C$-relatively compact  in $D(\R)$.
\end{enumerate}

If in addition, $D$ is closed under addition, then the above equivalence holds when ordinary relative compactness in $D$ replaces $C$-relative compactness in both the hypothesis and conclusion.
\item \label{item:Crel2} Let $E$ be   a Polish space and $D\subset C_c(E)$ be a separating class in $C_c(E)$. 
A sequence of c\`adl\`ag $M_{LF}(E)$-valued processes $\{X^N\}$ is $C$-relatively compact in $D(M_{LF}(E))$ if and only if 
for all $\phi\in D$, $\left\{X^N(\phi)\right\}$ is  $C$-relatively compact  in $D(\R)$.

If in addition, $D$ is closed under addition, then the above equivalence holds when ordinary relative compactness in $D$ replaces $C$-relative compactness in both the hypothesis and conclusion. 
\end{enumerate}
\end{theorem}

Returning to our setting, we take $D=C_b^2(\R^2)$ for $M_F(\R^2)$ and $D=C_c^2(\R^2)$ for $M_{LF}(\R^2)$ as the separating class of test functions for $M_{F}(\R^2)$ and $M_{LF}(\R^2)$, respectively.

\subsubsection{Condition \ref{item:CrelMF1} in Theorem \ref{thm:Crel} }
First, we verify the condition \ref{item:CrelMF1} in Theorem \ref{thm:Crel}.
\begin{lemma}\label{lem:tightequiv}
For all $\e>0$, $T>0$, there exists a compact set $K=K_{\e,T}\subset\R^2$ such that \begin{align*}
\sup_N P\left(\sup_{t\leq T}\sZ_t^{N,\phi}\left(K_{\e,T}^c\right)>\e\right)<\e.
\end{align*}
\end{lemma}
\begin{proof}
For the purpose of this proof, let
$\widehat{\sZ}_t^{N,\phi}$ denote the corresponding measure on the
path space of the random walk, defined by
\begin{align*}
\widehat{\sZ}_t^{N,\phi}(\dd S)
:=\frac{1}{N}
\sum_{x\in\Z^2}\phi_N(x)
E_x\left[\prod_{i=1}^{\lfloor Nt\rfloor}e_{i,S_i}
:\dd S\right].
\end{align*}
Then, $\sZ_t^{N,\phi}$ is the pushforward of
$\widehat{\sZ}_t^{N,\phi}$ under the endpoint map
$S\mapsto S_{\lfloor Nt\rfloor}/\sqrt N$.

Fix $\e>0$ and $T>0$. Then, by the invariance principle, we can take a compact set $K\subset \R^2$ such that \begin{align*}
\sup_{N\geq 1}\frac{1}{N}\sum_{x\in \Z^2}\phi_N\left(x\right)P_x\left(\frac{S_{n}}{\sqrt{N}}\in K^c \text{ for some $n\leq {NT}$}\right)<\e^2.
\end{align*}

Then, it is clear that 
\begin{align*}
\sZ^{N,\phi}_{t}(K^c)\leq \widehat{\sZ}^{N,\phi}_{t}\left(\left\{\frac{S_n}{\sqrt{N}}\in K^c \text{ for some $n\leq {NT}$}\right\}\right)
\end{align*}
and the right-hand side is an $\F_n$-martingale. Therefore, we have \begin{align*}
\IP\left(\sup_{t\leq T}\sZ^{N,\phi}_{t}(K^c)>\e\right)&\leq \IP\left(\sup_{t\leq T}\widehat\sZ^{N,\phi}_{t}\left(\left\{\frac{S_n}{\sqrt{N}}\in K^c \text{ for some $n\leq {NT}$}\right\}\right)>\e\right)\notag\\
&\leq \frac{1}{\e}\IE\left[\widehat\sZ^{N,\phi}_{T}\left(\left\{\frac{S_n}{\sqrt{N}}\in K^c \text{ for some $n\leq {NT}$}\right\}\right)\right]\notag\\
&=\frac{1}{\e}\frac{1}{N}\sum_{x\in\Z^2}\phi_N(x)P_x\left(\frac{S_{n}}{\sqrt{N}}\in K^c \text{ for some $n\leq NT$}\right)<\e,
\end{align*} 
where we have used Doob's maximal inequality in the second inequality.
\end{proof}

\subsubsection{Conditions \ref{item:CrelMF2} and \ref{item:Crel2} in Theorem \ref{thm:Crel}}
Next, we verify the conditions \ref{item:CrelMF2} and \ref{item:Crel2} in  Theorem \ref{thm:Crel} for $\{\mathsf{Z}_{t}^{N,\phi}\}$.

To see  the $C$-relative compactness in $D([0,T],\R)$ of  $\{\mathsf{Z}^{N,\phi}_{t}(\psi)\}_{t\geq 0}$, it is enough to verify the following two conditions \begin{enumerate}[label=(C-\arabic*)]
\item\label{item:C1} the $C$-relative compactness of $\left\{\int_{0}^{\frac{{Nt}}{N}} \mathsf{Z}_{s}^{N,\phi} \left(\Delta_N \psi_N\right)\dd s\right\}$ in $D([0,T],\R)$, and 
\item\label{item:C2} the $C$-relatively compactness of  $\left\{M^{N,\phi}_{{Nt}}(\psi)\right\}$ in $D([0,T],\R)$.
\end{enumerate}

As a corollary of \cite[Theorem 1.11]{CSZ25}, we can obtain a boundedness of the higher moments  of $\mathsf{Z}_u^{N,\phi}$ as follows.   
\begin{lemma}\label{lem:supZbdd} 
Fix $h\geq 2$. Then, for any $t\geq 0$, any $\phi\in C_c(\mathbb{R}^2)$ and $\psi\in C_b(\R^2)$, we have \begin{align*}
\sup_{N\geq 1}\mathbb{E}\left[\sup_{u\leq t}\left|\mathsf{Z}_u^{N,\phi}(\psi)\right|^h\right]<\infty.
\end{align*}
The same result holds for $\phi\in C_b(\mathbb{R}^2)$ and $\psi\in C_c(\R^2)$.
\end{lemma}

\begin{proof}
  By positivity,
\[
\left|\mathsf Z_u^{N,\phi}(\psi)\right|
\leq
\|\psi\|_\infty
\mathsf Z_u^{N,|\phi|}(1),
\]
Therefore, it is enough to see the case $\phi\in C_c^+(\R^2)$ and $\psi\equiv 1$ due to \eqref{eq:reversal}. Moreover, it is straightforward to see that for each $N\geq 1$, $\left\{\sZ^{N,\phi}_{\frac{k}{N}}(1)\right\}_{k\geq 0}$ is an $\F_n$-martingale. Therefore, Doob's maximal inequality yields that \begin{align*}
\IE\left[\sup_{u\leq t}\mathsf{Z}_{u}^{N,\phi}(1)^h\right]\leq \left(\frac{h}{h-1}\right)^h\IE\left[\mathsf{Z}_{t}^{N,\phi}(1)^h\right]\quad t\geq 0.
\end{align*}
\cite[Theorem 1.11]{CSZ25} shows the right-hand side is bounded.

\end{proof}

\begin{proof}[Proof of \ref{item:C1}] 
Fix $T\geq 0$, $h\geq 4$ and let $\phi\in C_c^+(\R^2)$ and $\psi\in C_b^2(\R^2)$. Then, we have \begin{align*}
\mathbb{E}\left[\left|\int_{s}^{t}\mathsf{Z}_{u}^{N,\phi}\left(\Delta_N\psi_N\right)\dd u\right|^h\right] 
&\leq \|\Delta_N\psi_N \|_\infty^h \mathbb{E}\left[\left(\int_{s}^{t}\mathsf{Z}_{u}^{N,\phi}\left(1\right)\dd u\right)^h\right]\\
&\leq \|\Delta_N\psi_N \|_\infty^h 
|t-s|^{h-1}
\int_s^t
\IE\left[
\mathsf Z_u^{N,\phi}(1)^h
\right]\dd u\\
&\leq C|t-s|^h.
\end{align*}
for $0\leq s\leq t$, where we have used Lemma \ref{lem:supZbdd} in the last inequality. 
The Kolmogorov criterion implies the $C$-relative compactness of the drift terms.

For the case where  $\phi\in C_b^+(\R^2)$ and $\psi\in C_c^2(\R^2)$, 
we can use the same argument by replacing $\phi\equiv 1$ 
and dividing $\Delta_N \psi_N=\left(\Delta_N \psi_N\right)\vee 0+\left(\Delta_N \psi_N\right)\wedge 0$.
\end{proof}

It remains only to verify condition \ref{item:C2} in order to complete the proof of Theorem \ref{thm:ContVer}.
 To prove \ref{item:C2}, we apply the general theory for $C$-tightness of martingales given in \cite[Lemma II.4.5.]{Per02} or the conclusion of \cite[Theorem VI.4.13, Theorem VI.3.26]{JS87} 

\begin{lemma}\label{lem:CrelM}
Let $\left\{M^N_{t},\mathcal{F}^N_{t}\right\}_{t=\frac{k}{N},k\in \mathbb{N}_0}$ be martingales with $M^N_{0}=0$. Let \begin{align*}
\left\lan M^N\right\ran_{t}
:=\sum_{\begin{smallmatrix}s\in \frac{1}{N}\mathbb{N}\\ 0\leq s< t\end{smallmatrix}} \mathbb{E}\left[ \left(\left.M^N_{s+\frac{1}{N}}-M^N_{s}\right)^2\right|\mathcal{F}_{s}^N\right],
\end{align*}
and extend $M^N_\cdot$ and $\lan M^N\ran_{\cdot}$ to $[0,\infty)$ as right-continuous step functions.

 Then, the following statements hold:
\begin{enumerate}[label=$(\arabic*)$]
\item\label{item:CrelM1} 
Suppose the following two conditions hold:
\begin{enumerate}[label=$(\mathrm{C}\text{-$2$-}\mathrm{\roman*})$]
\item\label{item:Crel1-1} $\left\{\left\{\left\lan M^{N}_{\cdot}\right\ran_{t}\right\}_{t\geq 0}:N\geq 1\right\}$ is $C$-relatively compact in $D(\R)$.
\item\label{item:Crel1-2} \begin{align}
\sup_{0\leq t\leq T}\left|M^{N}_{t+\frac{1}{N}}-M^{N}_{t}\right|\to 0 \quad \text{in probability for all $T>0$}. \label{eq:jumcontrol}
\end{align}
\end{enumerate}
Then, $M^{N}_{\cdot}$ is $C$-relatively compact in $D(\R)$.
\item\label{item:CrelM2} If, in addition, \begin{align}
\left\{\left(M^{N}_{s}\right)^2+\left\lan M^{N}\right\ran_{s}:N\geq 1\right\}\quad \text{is uniformly integrable for all $s\in [0,T]$},
\end{align} 
then $M^{N_k}_{\cdot}\Rightarrow \mathscr{M}_{\cdot}$ implies that $\mathscr{M}_{\cdot}$ is a continuous $L^2$-martingale with respect to the filtration $\mathcal{F}_t^{\mathscr{M}}$ and that \begin{align*}
\left(M^{N_k}_{\cdot}, \left\lan M^{N_k}\right\ran_{\cdot }\right)\Rightarrow \left(\mathscr{M}_{\cdot},\left\lan \mathscr{M}\right\ran_\cdot\right),
\end{align*}
where $\{\mathcal F_t^{\mathscr M}\}_{t\geq0}$ denotes the
right-continuous filtration generated by $\mathscr M$ as given in \eqref{eq:filtration}.
\end{enumerate}
\end{lemma}

In the following three subsections, we verify conditions in Lemma \ref{lem:CrelM} for the sequence $\left\{M_{\cdot}^{N,\phi}(\psi)\right\}_{N\geq 1}$. 
This establishes the $C$-tightness of the processes $\left\{M_{\cdot}^{N,\phi}(\psi)\right\}_{N\geq 1}$ and $\left\{\mathsf{Z}_\cdot^{N,\phi}(\psi)\right\}_{N\geq 1}$. 
In particular, any limit point is continuous, and hence  the limiting process $\mathscr{Z}^{\vartheta,\phi}_\cdot$ is continuous.

\subsubsection{Proof of \ref{item:Crel1-1} in Lemma \ref{lem:CrelM} for $M^{N,\phi}(\psi)$}

To establish \ref{item:Crel1-1} in Lemma \ref{lem:CrelM} for $M^{N,\phi}(\psi)$, we follow a standard approach based on moment estimates of the quadratic variations.

\begin{lemma}\label{lem:4thQV}
For each $ T> 0$, $p>1$, $\phi\in C_c^+(\R^2)$, and $\psi\in C_b(\R^2)$, 
there exists $C>0$ such that\begin{align}
\IE\left[\left(\left\lan M^{N,\phi}(\psi)\right\ran_{Nt}-\left\lan M^{N,\phi}(\psi)\right\ran_{Ns}\right)^2\right]\leq C|t-s|^{\frac{3}{2}-\frac{1}{2p}}, \quad \text{uniformly in } 0\leq s\leq t\leq T.\label{eq:qv2ndmoment}
\end{align}
The same result holds for the case where $\phi\in C_b^+(\R^2)$ and $\psi\in C_c^2(\R^2)$.
\end{lemma}
Taking $p>2$ and applying the Garsia-Rodemich-Rumsey inequality \cite[Theorem A.1]{FV10}, we obtain \ref{item:Crel1-1}.
\begin{remark}
Lemma \ref{lem:4thQV} provides an upper bound $|t-s|^{\frac{3}{2}-o(1)}$, and it is not expected to  be  optimal.
 Actually, it follows  from \cite[Lemma 4.3]{Tsa24} that 
 the left-hand side of \eqref{eq:qv2ndmoment} admits an upper bound of order $|t-s|^{2-o(1)}$ for $\phi,\psi\in C_c(\R^2)$.  
\end{remark}

We postpone the proof of Lemma \ref{lem:4thQV} to Section \ref{sec:moments}.

Here, we demonstrate an easy estimate of the difference between the values of $\langle M^{N,\phi}(\psi)\rangle_{Nt} $. 
For $0<s<t<T$, \begin{align}
&\frac{1}{4\pi}\lim_{N\to\infty}\mathbb{E}[\langle M^{N,\phi}(\psi)\rangle_{Nt}-\langle M^{N,\phi}(\psi)\rangle_{Ns}]\notag\\
&\leq  \int_{\R^2}\int_{0<u<v<t}\Phi_u(x)^2 G_\vartheta(v-u)\dd u\dd v\dd x- \int_{\R^2}\int_{0<u<v<s}\Phi_u(x)^2 G_\vartheta(v-u)\dd u\dd v\dd x\notag\\
&\leq  \int_{\R^2}\int_{0<u<s<v<t}\Phi_u(x)^2 \widehat{G}_{\vartheta,T}(v-u)\dd u\dd v\dd x+ \int_{\R^2}\int_{s<u<v<t}\Phi_u(x)^2 \widehat{G}_{\vartheta,T}(v-u)\dd u\dd v\dd x\notag\\
&\leq C_{\phi,\vartheta,T,1}\frac{t-s}{\log \frac{e^2 T}{s}}+C_{\phi,\vartheta,T,2}\frac{t-s}{\log \frac{e^2 T}{t-s}},\notag
\end{align}
where we have used $\widehat{G}_{\vartheta,T}(v-u)\leq \widehat{G}_{\vartheta,T}(s-u)$ for $0<u<s<v<T$ in the first term, 
and   $\widehat{G}_{\vartheta,T}(t)=\frac{\dd }{\dd t}\frac{c_{\vartheta,T}}{\log \left(\frac{2T}{t}\right)}$ for $t\in (0,T]$ 
and 
\begin{align*}
\int_s^t \frac{1}{\log \left(\frac{e^2 T}{t-u}\right)}\dd u\sim \frac{t-s}{e^2 T} \frac{1}{\log \left(\frac{e^2 T}{t-s}\right)} \quad \text{as }t-s\searrow 0
\end{align*}
in the second term.

\subsubsection{Proof of \ref{item:Crel1-2}  in Lemma \ref{lem:CrelM} for $M^{N,\phi}(\psi)$}

To prove \ref{item:Crel1-2}  for $M^{N,\phi}_\cdot(\psi)$, we first remark that 
 \begin{align}
\sup_{0\leq k\leq {NT}-1}\left|\delta M_{N,k}^{\phi}(\psi)\right|^3\leq \sum_{k=0}^{{NT}-1}  \left|\delta M_{N,k}^{\phi}(\psi)\right|^3\label{eq:supleqsum}
\end{align}

To prove \eqref{eq:jumcontrol}, we use the following Burkholder type inequality \cite[(PSF) in p.152]{Per02} and \cite[Theorem 21.1]{Bur73}.
\begin{lemma}\label{lem:PSF} 
Let $f:[0,\infty)\to [0,\infty)$ be a continuous increasing function with $f(0)=0$ such that there exists $c_0\in (0,\infty)$ such that  $f(2\lambda)\leq c_0f(\lambda)$ for any $\lambda\geq 0$.

Let $\{M_n\}_{n\geq 0}$ be an $\F_n$-martingale. We set $M_n^*=\sup_{k\leq n}|M_k|$, \begin{align*}
&\langle M\rangle_n:=\sum_{k=1}^n \mathbb{E}\left[(M_k-M_{k-1})^2|\F_{k-1}\right]+E[M_0^2]\\
&d_n^*:=\max_{1\leq k\leq n}|M_k-M_{k-1}|.
\end{align*}
Then, we have \begin{align*}
\mathbb{E}\left[f(M_n^*)\right]\leq c\left(\mathbb{E}\left[f\left(\langle M\rangle_n^{\frac{1}{2}}\right)\right]+\mathbb{E}\left[f(d_n^*)\right]\right).
\end{align*}
\end{lemma}

\begin{proof}[Proof of \eqref{eq:jumcontrol} for $M_{n}^{N,\phi}(\psi)$]
Conditioned on $\F_{n-1}$, $\delta M_{N,n}^{\phi}(\psi)$ is a sum of independent mean-zero random variables $\sum_{y\in \Z^2}\widetilde{M}^{(N,n)}_y(\phi,\psi)$, where \begin{align*}
\widetilde{M}^{(N,n)}_y(\phi,\psi):=\frac{1}{N}\overline{Z}_{N,\frac{n}{N}}^{\phi }(y)\psi_N(y)\xi^{\b_N}_{n,y}. 
\end{align*}
Let $\{\Lambda_k\}_{k\geq 0}$ be an increasing sequence of finite
subsets of $\mathbb Z^2$ such that
 $|\Lambda_k|=k$  and $\bigcup_{k\geq 0}\Lambda_k=\Z^2$. We define the filtration $\mathcal{F}_{n,\Lambda_k}=\mathcal{F}_{n-1}\vee \sigma[\omega_{n,x}, x\in \Lambda_k]$. Then, \begin{align*}
\widetilde{M}_{0}^{(N,n)}(\phi,\psi):=0,\quad \widetilde{M}_{k}^{(N,n)}(\phi,\psi):=\sum_{y\in \Lambda_k}\widetilde{M}_y^{(N,n)}(\phi,\psi)\quad (k\geq 1)
\end{align*}
is an $\mathcal{F}_{n,\Lambda_k}$-martingale, and $\lim_{k\to\infty}\widetilde{M}_{k}^{(N,n)}(\phi,\psi)=\delta M_n^{N,\phi}(\psi)$ 
since, for fixed $N$ and $n$, only finitely many terms are nonzero.

Thus, we can apply Lemma \ref{lem:PSF} with $f(\lambda)=\lambda^3$ ($c(f)=2^3$) 
to $\{\widetilde{M}^{(N,n)}_k(\phi,\psi)\}_{k\geq 0}$. 
Moreover, it follows that \begin{align*}
&\left\langle \widetilde{M}^{(N,n)}(\phi,\psi)\right\rangle_k
=\sum_{y\in \Lambda_k}\frac{\sigma_N^2}{N^2}\overline{Z}_{N,\frac{n}{N}}^{\phi}(y)^2\psi_N(y)^2, 
& \left|\widetilde{M}_y^{(N,n)}(\phi,\psi)\right|=\frac{\sigma_N}{N}\overline{Z}_{N,\frac{n}{N}}^\phi(y)|\psi_N(y)|
\end{align*}
by \eqref{eq:Bernoulli}.

Consequently,
\begin{align*}
\lim_{k\to\infty}
\left\langle
\widetilde M^{(N,n)}(\phi,\psi)
\right\rangle_k
&=
\frac{\sigma_N^2}{N^2}
\sum_{y\in\mathbb Z^2}
\overline Z_{N,n/N}^{\phi}(y)^2\psi_N(y)^2\\
&=
\left\langle M^{N,\phi}(\psi)\right\rangle_n
-
\left\langle M^{N,\phi}(\psi)\right\rangle_{n-1}.
\end{align*}

Applying Lemma \ref{lem:PSF} to
$\{\widetilde M_k^{(N,n)}(\phi,\psi)\}_{k\geq0}$ and then letting $k\to\infty$, 
we obtain \begin{align}
&\mathbb{E}\left[\left|\delta M_{N,n}^{\phi}(\psi)\right|^3\right]\label{eq:3rdbddmoment}\\
&\leq  c\mathbb{E}\left[\left(\langle M^{N,\phi}(\psi)\rangle_{n}-\langle M^{N,\phi}(\psi)\rangle_{n-1}\right)^\frac{3}{2}\right]
+c\mathbb{E}\left[\sum_{y\in \Z^2}\left(\frac{\sigma_N^2}{N^2}{\overline{Z}}_{N;\frac{n}{N}}^{\phi}(y)^2\psi\left(\frac{y}{\sqrt{N}}\right)^2\right)^\frac{3}{2}\right]\notag\\
&\leq c\mathbb{E}\left[\left(\langle M^{N,\phi}(\psi)\rangle_{n}-\langle M^{N,\phi}(\psi)\rangle_{n-1}\right)^2\right]^\frac{3}{4}
+c\|\psi\|_\infty^3\mathbb{E}\left[\sum_{y\in \Z^2}\frac{\sigma_N^3}{N^3}{\overline{Z}}_{N;\frac{n}{N}}^{\phi}(y)^3\right].\notag
\end{align}
Here, the first inequality follows from Lemma \ref{lem:PSF} and
\[
(d_k^*)^3
\leq
\sum_{y\in\Lambda_k}
\left|
\widetilde M_y^{(N,n)}(\phi,\psi)
\right|^3,
\]
while the second follows from H\"older's inequality.

By Lemma \ref{lem:4thQV}, for every $p>1$, there exists
$C=C_{p,T}>0$ such that, uniformly in
$1\leq n\leq  NT$, 
\begin{align*}
\IE\left[
\left(
\left\langle M^{N,\phi}(\psi)\right\rangle_n
-
\left\langle M^{N,\phi}(\psi)\right\rangle_{n-1}
\right)^2
\right]^{3/4}
\leq
C N^{-\frac98+\frac{3}{4p}}.
\end{align*}

In addition, it follows from \cite[Lemma 6.1 (6.4)]{CSZ19b} that for any $\e>0$ and $T\geq 0$, there exists $C>0$ such that \begin{align}
\sum_{1\leq n\leq NT}\sum_{y\in \Z^2}\mathbb{E}\left[{\overline{Z}}_{N;\frac{n}{N}}^{\phi}(y)^3\right]\leq CN^{\frac{5}{2}+\e}.\label{eq:CSZNTmoment}
\end{align}
Taking $p>6$, $\e\in(0,1/2)$, and combining
\eqref{eq:3rdbddmoment} with \eqref{eq:CSZNTmoment}, we obtain
\begin{align*}
\sum_{n=1}^{NT}
\IE\left[\left|\delta M_{N,n}^{\phi}(\psi)\right|^3\right]
&\leq C N^{-\frac{1}{8}+\frac{3}{4p}}
+C\sigma_N^3N^{-\frac{1}{2}+\e}\longrightarrow0.
\end{align*}
Together with \eqref{eq:supleqsum} and Markov's inequality, this proves \eqref{eq:jumcontrol}.

For the case with $\phi\in C_b^+(\R^2)$ and $\psi\in C_c^2(\R^2)$, we replace the second term on the right-hand side of
\eqref{eq:3rdbddmoment}  by \begin{align*}
c\|\phi\|_\infty^3\IE\left[\sum_{y\in \Z^2}\frac{\sigma_N^3}{N^3}{\overline{Z}}_{N;\frac{n}{N}}^{\phi}(y)^3\left|\psi\left(\frac{y}{\sqrt{N}}\right)\right|^3\right],
\end{align*}
whose sum over $1\leq n\leq\lfloor NT\rfloor$ is bounded by $\frac{C}{N^{1-\e}}$ due to \cite[Lemma 6.1 (6.2)]{CSZ19b}.
\end{proof}

\begin{remark}
For \eqref{eq:CSZNTmoment}, we see that \begin{align*}
\mathbb{E}\left[{\overline{Z}}_{N;\frac{n}{N}}^{\phi}(y)^3\right]
&=\mathbb{E}\left[\left({\overline{Z}}_{N;\frac{n}{N}}^{\phi}(y)-\IE\left[{\overline{Z}}_{N;\frac{n}{N}}^{\phi}(y)\right]+\IE\left[{\overline{Z}}_{N;\frac{n}{N}}^{\phi}(y)\right]\right)^3\right]\\
&=\IE\left[\left({\overline{Z}}_{N;\frac{n}{N}}^{\phi}(y)-\IE\left[{\overline{Z}}_{N;\frac{n}{N}}^{\phi}(y)\right]\right)^3\right]\\
&\hspace{2em}+3\IE\left[\left({\overline{Z}}_{N;\frac{n}{N}}^{\phi}(y)-\IE\left[{\overline{Z}}_{N;\frac{n}{N}}^{\phi}(y)\right]\right)^2\right]\IE\left[{\overline{Z}}_{N;\frac{n}{N}}^{\phi}(y)\right]\\
&+\IE\left[{\overline{Z}}_{N;\frac{n}{N}}^{\phi}(y)\right]^3.
\end{align*}
Our setting of $\{\omega_{n,x}\}$, in particular \eqref{eq:omegamoment}, allows us to  use \cite[Lemma 6.1 (6.4)]{CSZ19b} for an upper bound of the first term. 
More precisely, the corresponding contribution is divided into two parts,
 \textit{non-triple intersections} and \textit{triple intersections}, and the latter one vanishes under \eqref{eq:omegamoment}. The second term can be easily estimated by \eqref{eq:varUN} and Lemma \ref{thm:UN}.

\end{remark}

Thus, \ref{item:CrelM1} in Lemma \ref{lem:CrelM} follows when we verify Lemma \ref{lem:4thQV}. 

\subsubsection{Proof of \ref{item:CrelM2} in Lemma \ref{lem:CrelM} for $M^{N,\phi}_\cdot(\psi)$}

\begin{proof}[Proof of \ref{item:CrelM2} in Lemma \ref{lem:CrelM}] 
 Applying Lemma \ref{lem:PSF} with $f(\lambda)=\lambda^3$,
 we obtain 
 \begin{align*}
\IE\left[\sup_{0\leq s\leq T}\left|M_{N;Ns}^{\phi}(\psi)\right|^3\right]
\leq C\IE\left[\left\langle M_{N;\cdot}^{\phi}(\psi)\right\rangle_{NT}^\frac{3}{2}\right]
+C\IE\left[\sup_{1\leq k\leq NT-1}\left| \delta M_{N;k}^{\phi}(\psi)\right|^3\right].
\end{align*}
The second term is uniformly bounded by the estimate obtained in the proof of \ref{item:Crel1-2}.
Moreover,  the first term on the right-hand side is bounded by 
$\dis C \|\psi\|_\infty^3\IE\left[\left\langle M^\phi_{N;\cdot}(1)\right\rangle_{NT}^\frac{3}{2}\right]$
and hence it is bounded by
\begin{align*}
C \|\psi\|_\infty^3\IE\left[\left(\sum_{1\leq k\leq NT}\left| \delta M_{N;k}^{\phi}(1)\right|^2\right)^\frac{3}{2}\right]
\end{align*}
by \cite[Lemma 16.1]{Bur73}. 
Also, combining Theorem 15.1 in \cite{Bur73} and Doob's $L^p$-inequality, this is dominated by \begin{align*}
C\IE\left[\left|M_{N;NT}^{\phi}(1)\right|^3\right]\leq C\IE\left[\mathsf{Z}_{N;T}^{|\phi|}(1)^3\right]+C\mathsf{Z}_{N;0}^{|\phi|}(1)^3.
\end{align*}
Here we used the identity $M_{N,Ns}^{\phi}(1)=\mathsf{Z}_{N;s}^{\phi}(1)-\mathsf{Z}_{N;0}^{\phi}(1)$. 
The uniform boundedness of the right-hand side in $N$ follows from 
\cite[Theorem 1.4]{CSZ19b}. Consequently, for every $s\in[0,T]$, the family
\begin{align*}
\left\{\left|M_{N;Ns}^{\phi}(\psi)\right|^2
+\left\langle M_{N;\cdot}^{\phi}(\psi)\right\rangle_{Ns}:N\geq 1
\right\}
\end{align*}
is uniformly integrable.

For the case where $\phi\in C_b^{+}(\R^2)$ and $\psi\in C_c^2(\R^2)$, we can use the same argument by exchanging the roles of $\phi$ and $\psi$. 
\end{proof}

Thus, all the conditions in Lemma \ref{lem:CrelM} have been verified,
provided that Lemma \ref{lem:4thQV} holds. 
The proof of Theorem
\ref{thm:ContVer} will therefore be complete once Lemma
\ref{lem:4thQV} is proved in Section \ref{sec:moments}.

\subsection{Martingale $ \mathscr{M}_t^{\vartheta,\phi}(\psi)$}\label{subsec:martingale}

The estimates established in the preceding subsections verify all
the conditions of Lemma \ref{lem:CrelM}; hence
Theorem \ref{thm:ContVer} follows.

By the $C$-tightness established in the previous subsection and the
finite-dimensional convergence in Theorem \ref{thm:SHF}, we have
\begin{align*}
\mathsf Z^{N,\phi}
\Rightarrow
\mathscr Z^{\vartheta,\phi}
\quad\text{in }D(M_F(\R^2)).
\end{align*}

Let $\{N_k\}_{k\geq1}$ be a subsequence along which convergence in
distribution holds. By the Skorokhod representation theorem, after
passing to a new probability space, we may assume that
$\mathsf Z^{N_k,\phi}
\longrightarrow
\mathscr Z^{\vartheta,\phi}
\quad\text{a.s. in }D(M_F(\R^2))$.

We next identify the limit of each term on the right-hand side of
\eqref{eq:semimart}.
Indeed, $\mathsf{Z}^{N_k,\phi}_0(\psi)\to \int \phi(x)\psi(x)\dd x$, and Taylor's theorem implies that for each $x\in \Z^2,y\in \Sigma$, and $k\geq 1$, 
there exists $c=c_{N_k,x,y}\in (0,1)$ such that\begin{align*}
\psi_{N_k}(x+y)-\psi_{N_k}(x)=\nabla\psi\left(\frac{x}{\sqrt{N_k}}\right)\cdot \frac{y}{\sqrt{N_k}}+\frac{1}{2N_k}(y_1,y_2)\mathrm{Hess}(\psi)\left(\frac{x+cy}{\sqrt{N_k}}\right)\left(\begin{matrix}y_1\\ y_2\end{matrix}\right),
\end{align*}
where $\mathrm{Hess}(\psi)$ is the Hessian matrix of $\psi$.

Therefore, we have \begin{align*}
\Delta_{N_k}\psi_{N_k}(x)=N\sum_{y\in \Sigma}q_1(y)\left(\psi_{N_k}(x+y)-\psi_{N_k}(x)\right)=\frac{1}{2}\Delta \psi\left(\frac{x}{\sqrt{N_k}}\right)+o(1)
\end{align*} 
uniformly over $x\in\mathbb Z^2$ such that
$x/\sqrt N_k\in K$
for each compact set $K\subset \R^2$.

Combining this locally uniform convergence with
Lemma \ref{lem:tightequiv}, we obtain, almost surely,
 \begin{align*}
\int_0^{t}\sZ_{s}^{N_k,\phi}(\Delta_{N_k}\psi_{N_k})\dd s\to \frac{1}{2}\int_0^t \scZ_{s}^{\vartheta,\phi}(\Delta \psi)\dd s.
\end{align*}

Thus, we obtain that $M^{N_k,\phi}_{N_k\cdot}(\psi)\Rightarrow \mathscr{M}^{\vartheta,\phi}_\cdot(\psi)$ in $D([0,\infty),\R)$ and hence, Lemma \ref{lem:CrelM} implies that\begin{align}
\left(M^{N_k,\phi}_{N_k\cdot}(\psi),\left\langle M^{N_k,\phi}(\psi)\right\rangle_{N_k\cdot}\right)\Rightarrow \left(\mathscr{M}^{\vartheta,\phi}(\psi),\left\langle \mathscr{M}^{\vartheta,\phi}(\psi)\right\rangle_\cdot\right).\label{eq:martingaleconvs}
\end{align}
At this stage, Lemma \ref{lem:CrelM} only shows that
$\mathscr M^{\vartheta,\phi}(\psi)$ is a martingale with respect to
its own natural filtration. It remains to show that it is a
martingale with respect to
$\{\mathcal F_t^{\mathscr Z^{\vartheta,\phi}}\}_{t\geq0}$.

The proof of Theorem \ref{thm:SHFMP} will be complete once the
following two lemmas are proved.

\begin{lemma}\label{lem:martingaleZ}
For any $\vartheta\in \R$, $\phi\in C_c^+(\R^2)$, and $\psi\in C_b^2(\R^2)$, 
$\mathscr{M}^{\vartheta,\phi}_t(\psi)$ is a  continuous $\mathcal{F}_t^{\mathscr{Z}^{\vartheta,\phi}}$-martingale.  
The same result holds for $\phi\in C_b^+(\R^2)$, and $\psi\in C_c^2(\R^2)$.
\end{lemma}

\begin{lemma}\label{lem:QuadConvL2}
For any $\vartheta\in \R$, $\phi\in C_c^+(\R^2)$, and $\psi\in C_b^2(\R^2)$, \begin{align*}
\dis \sup_{t\in [0,T]}\left|\left\lan  \mathscr{M}^{\vartheta,\phi}(\psi)\right\ran_t-\frac{4\pi}{-\log \e}\int_0^t \int_{\R^2}\mathscr{Z}^{\vartheta,\phi}_s\left(p_\e(\cdot-z)\right)^2\psi(z)^2\dd z\dd s\right|\to 0 
\end{align*}
in probability for any $T\geq 0$. 

The same result holds for $\phi\in C_b^+(\R^2)$, and $\psi\in C_c^2(\R^2)$.
\end{lemma}

We also obtain the following formula for the expected quadratic variation.
\begin{corollary}\label{cor:QuadMoment}
For any $\vartheta\in \R$, $\phi\in C_c^+(\R^2)$, $\psi\in C_b^2(\R^2)$, and $t>0$,\begin{align*}
\dis \IE\left[\left\lan  \mathscr{M}^{\vartheta,\phi}(\psi)\right\ran_t\right]={4\pi}\int_{0<u<v<t}\dd u\dd v \int_{\R^2}\dd x\int_{\R^2}\dd y\Phi_u^2(x)G_\vartheta(v-u,y-x)\psi(y)^2. 
\end{align*} 
The same result holds for $\phi\in C_b^+(\R^2)$, and $\psi\in C_c^2(\R^2)$.
\end{corollary}
\begin{proof}
Since we have proved that  $\left\langle M^{N,\phi}(\psi)\right\rangle_{Nt}$ is uniformly integrable, we have \begin{align*}
\IE\left[\left\lan  \mathscr{M}^{\vartheta,\phi}(\psi)\right\ran_t\right]=\lim_{N\to \infty}\IE\left[\left\langle M^{N,\phi}(\psi)\right\rangle_{Nt}\right].
\end{align*}
The desired formula then follows by the same argument as in the proof of \cite[Theorem 1.2]{CSZ19b}.
\end{proof}

\subsubsection{Proof of Lemma \ref{lem:martingaleZ}}

Since $\mathsf Z_{k/N}^{N,\phi}$ is $\mathcal F_k$-measurable for
each $k$, we have
\begin{align*}
\mathcal F_{n/N}^{\mathsf Z^{N,\phi}}
\subset
\mathcal F_n.
\end{align*}

Let $0\leq s_1\leq \dots\leq s_m\leq s<t$, $m\geq 1$ and $F$ be a bounded continuous function on $M_F(\R^2)^m$. Then, we have \begin{align*}
&\IE\left[\left(M^{N,\phi}_{Nt}(\psi)-M^{N,\phi}_{Ns}(\psi)\right)F\left(\mathsf{Z}^{N,\phi}_{s_1},\dots,\mathsf{Z}^{N,\phi}_{s_m}\right)\right]=0
\end{align*}
Since $F$ is bounded and the martingale increments are uniformly
integrable, we may pass to the limit in the preceding identity, that is  \begin{align*}
&\IE\left[\left(\mathscr{M}^{\vartheta,\phi}_{t}(\psi)-\mathscr{M}^{\vartheta,\phi}_{s}(\psi)\right)F\left(\mathscr{Z}^{\vartheta,\phi}_{s_1},\dots,\mathscr{Z}^{\vartheta,\phi}_{s_m}\right)\right]=0.
\end{align*} 
By a monotone-class argument, the preceding identity extends to every
bounded random variable measurable with respect to the raw natural
filtration generated by
$\{\mathscr Z_r^{\vartheta,\phi}:0\leq r\leq s\}$.
The continuity of
$\mathscr M^{\vartheta,\phi}(\psi)$, together with the right
continuity and completion of the filtration, then yields the
martingale property with respect to
$\{\mathcal F_t^{\mathscr Z^{\vartheta,\phi}}\}_{t\geq0}$.
This completes the proof of Lemma \ref{lem:martingaleZ}.
\begin{flushright}\toybox\end{flushright}

\subsubsection{Proof of Lemma \ref{lem:QuadConvL2}}
For every $t\geq0$ and $\e>0$, we know that \begin{align*}
\int_0^t \int_{\R^2} \scZ^{\vartheta,\phi}_{s}(p_\e(\cdot-z))^2 \psi(z)^2\dd z\dd s=\lim_{k\to \infty}\int_0^ {\frac{\lr{N_kt}}{N_k}}\int_{\R^2} \overline{\sZ}^{N_k,\phi}_{s}(p_\e(\cdot-z))^2 \psi(z)^2\dd z\dd s\quad \text{a.s.}
\end{align*}
For $N\geq1$ and $\e>0$, define the mollified
quadratic-variation functional by \begin{align*}
\left\langle M^{N,\phi}(\psi)\right\rangle_t^{(\e)}:=\frac{4\pi}{-\log \e}\int_0^ {\frac{\lr{Nt}}{N}}\int_{\R^2} \overline{\sZ}^{N,\phi}_{s}(p_\e(\cdot-z))^2 \psi(z)^2\dd z\dd s.
\end{align*}

Then, 
Fatou's lemma with Lemma \ref{lem:CrelM} implies that 
\begin{align*}
&\IE\left[	\left(\frac{{4\pi}}{-\log \e}\int_0^t\int_{\R^2} \scZ^{\vartheta,\phi}_{s}(p_\e(\cdot-z))^2 \psi(z)^2\dd z\dd s- \left\lan \mathscr{M}^{\vartheta,\phi}(\psi)\right\ran_{t}\right)^2		\right]\notag\\
&\leq \varliminf_{k\to\infty}\IE\left[\left(\left\langle M^{N_k,\phi}(\psi)\right\rangle_t^{(\e)}-\left\langle M^{N_k,\phi}(\psi)\right\rangle_{Nt}\right)^2\right].
\end{align*}

We prove the following lemma. 

\begin{lemma}\label{lem:QuadProcEst}
We have \begin{align}
{ \lim_{\e\to 0}\varlimsup_{N\to\infty}\IE\left[\left(\left\langle M^{N,\phi}(\psi)\right\rangle_t^{(\e)}-\left\langle M^{N,\phi}(\psi)\right\rangle_{Nt}\right)^2\right]=0}\label{eq:limit-discrete}
\end{align}
for $t\geq 0$, $\phi\in C_c^+(\R^2)$, and  $\psi\in C_b^2(\R^2)$. The same result holds for $\phi\in C_b^+(\R^2)$, and $\psi\in C_c^2(\R^2)$.
\end{lemma}
The proof of Lemma \ref{lem:QuadProcEst} is postponed to Section \ref{sec:momentsbound} , and we complete the proof of Lemma \ref{lem:QuadConvL2} by using Lemma \ref{lem:QuadProcEst}. 

\begin{proof}[Proof of Lemma \ref{lem:QuadConvL2}]

Let $D\subset [0,T]$ be a countable dense subset. Then, it follows from Lemma \ref{lem:QuadProcEst} that for any sequence $\{\e_n\}_{n\geq 1}$ with $\e_n\to 0$, there exists a subsequence $\{\e_{n_k}\}$ such that 
\begin{align*}
\frac{{4\pi}}{-\log \e_{n_k}}\int_0^t\int_{\R^2} \scZ^{\vartheta,\phi}_{s}(p_{\e_{n_k}}(\cdot-z))^2 \psi(z)^2\dd z\dd s\to \left\lan \mathscr{M}^{\vartheta,\phi}_{\cdot}(\psi)\right\ran_{t}
\end{align*}
for $t\in D$ a.s. 
Since each approximating process is non-decreasing and the limiting
process is continuous, convergence on the dense set $D$ implies
uniform convergence on $[0,T]$. Furthermore, since every sequence $\{\e_n\}$ tending to $0$ has a subsequence
along which the convergence is uniform on $[0,T]$ almost surely,
the full family converges uniformly on $[0,T]$ in probability.
This proves Lemma \ref{lem:QuadConvL2}.

\end{proof}

\begin{remark}\label{rem:moment4th}
We remark that every term in expectation in \eqref{eq:qv2ndmoment} and  \eqref{eq:limit-discrete} can be described via partition functions from \eqref{eq:QVprocess}. More precisely, we need look at  
\begin{align}
\frac{\sigma_N^4}{N^2}\int_{[s,t]^2} \dd u \dd v \sum_{y^1,y^2}\IE\left[\overline{Z}_{N;u}^{\phi}(y^1)^2\overline{Z}_{N;v}^{\phi}(y^2)^2\right]\psi_N(y^1)^2\psi_N(y^2)^2\label{eq:4thmomentQV}
\end{align}
for Lemma \ref{lem:4thQV}, and the linear combination of 
\begin{align}
&\frac{\sigma_N^4}{N^2}\int_{[0,t]^2} \dd u \dd v \sum_{y^1,y^2}\IE\left[\overline{Z}_{N;u}^{\phi}(y^1)^2\overline{Z}_{N;v}^{\phi}(y^2)^2\right]\psi_N(y^1)^2\psi_N(y^2)^2,\label{eq:QVmoments1}\\
&\frac{\sigma_N^2}{N^3}\int_{[0,t]^2} \dd u \dd v \sum_{y^1,y^2,y^3}\IE\left[\overline{Z}_{N;u}^{\phi}(y^1)\overline{Z}_{N;u}^{\phi}(y^2)\overline{Z}_{N;v}^{\phi}(y^3)^2\right]\label{eq:QVmoments2}\\
&\hspace{3em}\times \int_{\R^2}\dd zp_\e\left(z,\frac{y^1}{\sqrt{N}}\right)p_\e\left(z,\frac{y^2}{\sqrt{N}}\right)\psi(z)^2 \psi_N(y^3)^2,\notag\\
&\frac{1}{N^4}\int_{[0,t]^2} \dd u \dd v \sum_{y^1,y^2,y^3,y^4}\IE\left[\overline{Z}_{N;u}^{\phi}(y^1)\overline{Z}_{N;u}^{\phi}(y^2)\overline{Z}_{N;v}^{\phi}(y^3)\overline{Z}_{N;v}^{\phi}(y^4)\right]\label{eq:QVmoments3}\\
&\hspace{1.5em}\times \int_{(\R^2)^2}\hspace{-.5em}\dd z_1\dd z_2p_\e\left(z_1,\frac{y^1}{\sqrt{N}}\right)p_\e\left(z_1,\frac{y^2}{\sqrt{N}}\right)p_\e\left(z_2,\frac{y^3}{\sqrt{N}}\right)p_\e\left(z_2,\frac{y^4}{\sqrt{N}}\right)\psi(z_1)^2 \psi(z_2)^2.\notag
\end{align}
for Lemma \ref{lem:QuadProcEst}.

Thus, we will focus entirely on computing moments of partition
functions in the proofs of Lemmas \ref{lem:4thQV} and
\ref{lem:QuadProcEst}.
\end{remark}

\section{Proof of Theorem \ref{thm:martinalemeasure}}\label{sec:Quadconv}

In this section, we prove Theorem
\ref{thm:martinalemeasure} by extending the martingale map from
smooth test functions to bounded Borel functions and identifying the
quadratic variation of the extension.

\begin{proof}[Proof of  Theorem \ref{thm:martinalemeasure}]
First, we prove the existence of an extension of $\mathscr{M}^{\vartheta,\phi}_t(\psi)$ for $\psi\in \mathcal{B}_b(\R^2)$.

To extend the definition of the martingale to bounded Borel functions, we approximate $\psi$ by smooth functions. 
Let $\psi$ be a bounded Borel function.
Then, Lusin's theorem and Tietze extension theorem imply that there exists a sequence $\{\hat{\psi}_n\}$ in $C_b(\R^2)$ such that \begin{align*}
\hat{\psi}_n(x)\to \psi(x)\quad \text{a.e.~$x$ and  }\|\hat{\psi}_n\|_\infty=\|\psi\|_\infty \text{ for all $n\geq 1$.}
\end{align*}
(See \cite[Remark 1.3.30]{Tao11}.)
Also, we know any bounded continuous function $f$ can be approximated by $p_\epsilon*f$ uniformly on  compact sets.
 
For each $\psi\in \mathcal{B}_b(\R^2)$, we define $\psi_n(x)=p_{\epsilon_n}*\hat{\psi}_n(x)\in C_b^\infty(\R^2)$, where $\epsilon_n$ is chosen such that \begin{align*}
p_{\epsilon_n}*\hat{\psi}_n(x)\to \psi(x)\quad \text{for Lebesgue-a.e.~$x$ and  }\|p_{\epsilon_n}*\hat{\psi}_n\|_\infty\leq \|\psi\|_\infty \text{ for all $n\geq 1$.}
\end{align*}

Then, it follows from  the Burkholder-Davis-Gundy inequality and Corollary \ref{cor:QuadMoment} that for $n,m\geq 1$, and $t\geq 0$ \begin{align}
&\IE\left[\sup_{0\leq s\leq t}\left(\mathscr{M}_s^{\vartheta,\phi}(\psi_n)-\mathscr{M}_s^{\vartheta,\phi}(\psi_m)\right)^2\right]\label{eq:CauchyPsi}\\
&\qquad \leq C\IE\left[\left\langle\mathscr{M}^{\vartheta,\phi}(\psi_n-\psi_m)\right\rangle_t\right]\notag\\
&\qquad = 4\pi C\int_{0<u<v<t}\dd u\dd v\int_{\R^2\times \R^2}\dd x\dd y\Phi_u(x)^2 G_{\vartheta}(v-u,y-x)\left(\psi_n(y)-\psi_m(y)\right)^2\notag
\end{align} 
for some constant $C>0$. 
Thus,
$\{\mathscr M^{\vartheta,\phi}(\psi_n)\}_{n\geq1}$
is a Cauchy sequence in $\mathcal M_c^2$ on every compact time interval.
We denote its limit by
$\mathscr M^{\vartheta,\phi}(\psi)$.

By polarization and the Kunita--Watanabe inequality, it follows that for $n,m\geq 1$\begin{align*}
\left|\left\langle\mathscr{M}^{\vartheta,\phi}(\psi_n)\right\rangle_s -\left\langle\mathscr{M}^{\vartheta,\phi}(\psi_m)\right\rangle_s\right|\leq \left\langle\mathscr{M}^{\vartheta,\phi}(\psi_n+\psi_m)\right\rangle_t^\frac{1}{2}\left\langle\mathscr{M}^{\vartheta,\phi}(\psi_n-\psi_m)\right\rangle_t^\frac{1}{2}.
\end{align*}
for $0\leq s\leq t$, and hence \begin{align}
&\IE\left[\sup_{0\leq s\leq t}\left|\left\langle\mathscr{M}^{\vartheta,\phi}(\psi_n)\right\rangle_s -\left\langle\mathscr{M}^{\vartheta,\phi}(\psi_m)\right\rangle_s\right|\right]\label{eq:CauchyQuad}\\
&\leq E\left[\left\langle\mathscr{M}^{\vartheta,\phi}(\psi_n+\psi_m)\right\rangle_t\right]^\frac{1}{2}E\left[\left\langle\mathscr{M}^{\vartheta,\phi}(\psi_n-\psi_m)\right\rangle_t\right]^\frac{1}{2}.\notag
\end{align}
Since we have already shown that the right-hand side converges to zero, $\dis \left\{\left\langle \mathscr{M}^{\vartheta,\phi}(\psi_n)\right\rangle_\cdot\right\}_{n\geq 1}$ is an $L^1$-Cauchy sequence and $C$-tight. 
The limit, denoted by $\dis \left\langle \mathscr{M}^{\vartheta,\phi}(\psi)\right\rangle_\cdot$, defines the quadratic variation of $\mathscr{M}^{\vartheta,\phi}_\cdot(\psi)$ since $\dis \mathcal{M}_t^{\vartheta,\phi}(\psi)^2-\left\langle \mathcal{M}^{\vartheta,\phi}(\psi)\right\rangle_t$ is a martingale. 

The same estimate applied to two different approximating sequences
shows that the limit is independent of the chosen approximation.
Linearity follows by approximating linear combinations, and
uniqueness follows from the density argument above.
The same estimate also shows
that the resulting extension is linear and unique in
$\mathcal M_c^2$.

Finally, taking $n\to\infty$ in the quadratic-variation formula for
$\psi_n$, using \eqref{eq:CauchyQuad} and dominated convergence,
gives the asserted formula for
$\langle\mathscr M^{\vartheta,\phi}(\psi)\rangle$.

For the extension from $C_c^2(\R^2)$ to $\mathcal{B}_c(\R^2)$, we adapt the same argument by using an approximation $\psi_n=\psi*f_n$ for a probability density $f\in C_c^\infty(\R^2)$ and $f_n(x)=nf(\sqrt{n}x)$ instead of $\psi_n$.
\end{proof}

\section{Proof of Theorem \ref{thm:peaks}}\label{sec:peaks}

\begin{proof}[Proof of Theorem \ref{thm:peaks}]
\ref{item:LARGE1} ~ Suppose that there exists $c>0$ such that \begin{align*}
\varlimsup_{\lambda\to \infty}\varlimsup_{\e\to 0}\iint_{\mathcal{T}_{s,t}(A,f,\lambda,\e)}\mathscr{Z}^{\vartheta,\phi}_u(f_\e(\cdot-x))\dd x\dd u>c
\end{align*}
with positive probability.

On this event, there exists a sequence $\lambda_j\to\infty$ such that \begin{align*}
\left\langle \mathscr{M}^{\vartheta,\phi}(A)\right\rangle_t-\left\langle \mathscr{M}^{\vartheta,\phi}(A)\right\rangle_s
&=-\lim_{\e\to 0}\frac{4\pi}{\log \e}\int_s^t \int_{A}\left(\mathscr{Z}^{\vartheta,\phi}_u(f_{\e}(\cdot-x))\right)^2\dd x\dd u\\
&\geq \varlimsup_{\e\to 0}{4\pi\lambda_j}\iint_{\mathcal{T}_{s,t}(A,f,\lambda_j,\e)}\mathscr{Z}^{\vartheta,\phi}_u(f_{\e}(\cdot-x))\dd x\dd u\geq 4\pi \lambda_j c
\end{align*}
for every $j$, which is impossible because the left-hand side is finite almost surely.

\ref{item:LARGE2}~
Choose $x_0\in A$ and $r>0$ such that
$
\overline B(x_0,r)\subset A,
$
where $B_x(r)=\{y\in \R^2:|y-x|<r\}$ is the open ball centered at $x$ with radius $r>0$.
Let $\psi\in C_b^2(\R^2)$ be positive functions such that \begin{align*}
&0\leq \psi(x)\leq 1 \quad \text{for  $x\in \R^2$}
\quad \text{and}\quad 
\psi(x)=\begin{cases}
0\quad &\text{for $x\in A^c$}\\
1\quad &\text{for $x\in \overline{B}_{x_0}(r)$}.
\end{cases}
\end{align*}

We focus on the martingale $\mathscr{M}_t^{\vartheta,\phi}(\psi)$. 
Then, we have \begin{align*}
\IP\left(\left\langle \mathscr{M}^{\vartheta,\phi}(\psi)\right\rangle_t-\left\langle \mathscr{M}^{\vartheta,\phi}(\psi)\right\rangle_s>0\right)>0.
\end{align*}
Indeed, we have \begin{align*}
\IP\left(\left\langle \mathscr{M}^{\vartheta,\phi}(\psi)\right\rangle_t-\left\langle \mathscr{M}^{\vartheta,\phi}(\psi)\right\rangle_s>0\right)
\geq \frac{1}{4}\frac{\IE\left[\left\langle \mathscr{M}^{\vartheta,\phi}(\psi)\right\rangle_t-\left\langle \mathscr{M}^{\vartheta,\phi}(\psi)\right\rangle_s\right]^2}
{\IE\left[\left(\left\langle \mathscr{M}^{\vartheta,\phi}(\psi)\right\rangle_t-\left\langle \mathscr{M}^{\vartheta,\phi}(\psi)\right\rangle_s\right)^2\right]}
\end{align*}
from the Paley-Zygmund inequality. 
Moreover,  Lemma \ref{lem:QuadProcEst}, \eqref{eq:martingaleconvs}, 
and   Fatou's lemma imply that the denominator on the right-hand side is finite. 
In addition, it follows from Corollary \ref{cor:QuadMoment} that the numerator is strictly positive.

Moreover, Remark \ref{rem:compactsupp} below implies that there exists a deterministic sequence
$\e_n\to 0$
such that \begin{align}
\left\langle \mathscr{M}^{\vartheta,\phi}(\psi)\right\rangle_t-\left\langle \mathscr{M}^{\vartheta,\phi}(\psi)\right\rangle_s
&=-\lim_{n\to\infty }\frac{4\pi}{\log \e_n}\int_s^t \int_{A}\left(\mathscr{Z}^{\vartheta,\phi}_u(f_{\e_n}(\cdot-x))\right)^2\psi(x)^2\dd x\dd u,\label{eq:Quadasconv}
\end{align}
almost surely.

 We set
\begin{align*}
&I_{s,t,\lambda}^{(1)}(\e):=-\frac{4\pi}{\log \e}\iint_{\mathcal{T}_{s,t}(A,f,\lambda,\e)}\left(\mathscr{Z}^{\vartheta,\phi}_u(f_\e(\cdot-x))\right)^2\psi(x)^2\dd x\dd u\\
&I_{s,t,\lambda}^{(2)}(\e):=-\frac{4\pi}{\log \e}\iint_{[s,t]\times A\backslash \mathcal{T}_{s,t}(A,f,\lambda,\e)}\left(\mathscr{Z}^{\vartheta,\phi}_u(f_\e(\cdot-x))\right)^2\psi(x)^2\dd x\dd u.
\end{align*} 
We remark that if $\dis \varliminf_{n \to \infty}I_{s,t,\lambda}^{(1)}(\e_n)>0$ 
for some $\lambda>0$, then $\varlimsup_{n\to\infty}\mathcal{T}_{s,t}(A,f,\lambda,\e_n)$ is nonempty.

Thus, it is enough to prove that 
\begin{align*}
\mathcal{E}:=\left\{		\left\langle \mathscr{M}^{\vartheta,\phi}(\psi)\right\rangle_t-\left\langle \mathscr{M}^{\vartheta,\phi}(\psi)\right\rangle_s>0		\right\}
\overset{\text{a.s.}}{\subset}\bigcup_{m\geq 1}\left\{\varliminf_{n\to\infty}I_{s,t,\lambda_m}^{(1)}(\e_n)>0\right\}
\end{align*}
for $\lambda_m=\frac{1}{m}$.

By passing to a further subsequence, if necessary, we may also assume
that
\begin{align*}
\int_{s}^t \mathscr{Z}^{\vartheta,\phi}_u(\psi)\dd u=\lim_{\e_n\to 0}\int_{s}^t \int_{\R^2}\mathscr{Z}^{\vartheta,\phi}_u(f_{\e_n}(\cdot-x))\psi(x)\dd x\dd u,\quad \text{a.s.}
\end{align*}
 We set \begin{align*}
&J_{s,t,\lambda}^{(1)}(\e):=\iint_{\mathcal{T}_{s,t}(A,f,\lambda,\e)}\mathscr{Z}^{\vartheta,\phi}_u(f_\e(\cdot-x))\psi(x)^2\dd x\dd u\\
&J_{s,t,\lambda}^{(2)}(\e):=\iint_{[s,t]\times A\backslash \mathcal{T}_{s,t}(A,f,\lambda,\e)}\mathscr{Z}^{\vartheta,\phi}_u(f_\e(\cdot-x))\psi(x)^2\dd x\dd u.
\end{align*}

Since we have $-\frac{\mathscr{Z}_u^{\vartheta,\phi}(f_\e(\cdot-x))}{\log \e}\leq \lambda$ for $(u,x)\in [s,t]\times A\backslash \mathcal{T}_{s,t}(A,f,\lambda,\e)$, \begin{align*}
I_{s,t,\lambda}^{(2)}(\e_n)&=-\frac{4\pi}{\log \e_n}\iint_{[s,t]\times A\backslash \mathcal{T}_{s,t}(A,f,\lambda,\e_n)}\left(\mathscr{Z}^{\vartheta,\phi}_u(f_{\e_n}(\cdot-x))\right)^2\psi(x)^2\dd x\dd u \\
&\leq  4\pi\lambda \iint_{[s,t]\times A\backslash \mathcal{T}_{s,t}(A,f,\lambda,\e_n)}\mathscr{Z}^{\vartheta,\phi}_u(f_{\e_n}(\cdot-x))\psi(x)\dd x\dd u 
=4\pi \lambda J_{s,t,\lambda}^{(2)}(\e_n).
\end{align*}
On $\mathcal E$, the quantity
\begin{align*}
\sup_n J_{s,t,\lambda}^{(2)}(\e_n)
\end{align*}
is finite. Hence, for each $\omega\in\mathcal E$, there exists
$m=m(\omega)$ such that
\begin{align*}
4\pi\lambda_m
\varlimsup_{n\to\infty}
J_{s,t,\lambda_m}^{(2)}(\e_n)
\leq
\frac{1}{2}
\left(
\langle \mathscr{M}^{\vartheta,\phi}(\psi)\rangle_t-\langle \mathscr{M}^{\vartheta,\phi}(\psi)\rangle_s
\right).
\end{align*}
Therefore, 
we have
\begin{align*}
\varlimsup_{n\to \infty}I_{s,t,\lambda_m}^{(2)}(\e_n)\leq \frac{1}{2}\left( \left\langle \mathscr{M}^{\vartheta,\phi}(\psi)\right\rangle_t-\left\langle \mathscr{M}^{\vartheta,\phi}(\psi)\right\rangle_s\right),
\end{align*} 
on $\mathcal{E}$, and hence, it follows that \begin{align*}
\varliminf_{n\to \infty}I_{s,t,\lambda_m}^{(1)}(\e_n)\geq  \frac{1}{2}\left( \left\langle \mathscr{M}^{\vartheta,\phi}(\psi)\right\rangle_t-\left\langle \mathscr{M}^{\vartheta,\phi}(\psi)\right\rangle_s\right)>0\quad \text{on }\mathcal{E}.
\end{align*}

Hence, on $\mathcal E$, there exists $m\geq1$ such that
\[
\limsup_{n\to\infty}
\mathcal T_{s,t}(A,f,\lambda_m,\e_n)
\ne\emptyset.
\]
Since \(\IP(\mathcal E)>0\), assertion
\ref{item:LARGE2} follows.
\end{proof}

\begin{appendix}

\section{Moments of partition functions}\label{sec:ChaosExpansion}\label{sec:MomentsChaos}

As mentioned in Remark \ref{rem:moment4th}, the fourth moments of partition functions are crucial for proving Lemma \ref{lem:4thQV} and Lemma \ref{lem:QuadProcEst}. 
In this section, we introduce the polynomial chaos expansion of the partition functions
 and derive the moment estimates needed in the proofs of these lemmas.

From now on, we will omit the parameter $N$ in the notation if it is clear from the context.

\subsection{Chaos expansion and moments}\label{label:ChaosEx}\label{subsec:4}
Let $\mathbb{T}$ be a countable set and  $\{\omega_{a}\}_{a\in \mathbb{T}}$ be independent Bernoulli-distributed random variables with $\IP(\omega_a=1)=\IP(\omega_a=-1)=\frac{1}{2}$.

For a finite subset $F\subset \mathbb{T}$, we define \begin{align*}
\omega_F=\prod_{a\in F}\omega_a
\end{align*}
and the polynomial chaos $P(\omega)$ is defined as a linear combination of $\{\omega_F:F\subset \mathbb{T}\text{ is finite}\}$, i.e.~\begin{align*}
P(\omega)=\sum_{F\subset \mathbb{T}:\text{finite}}a_F\omega_F
\end{align*}
for some $\{a_F\}_{F\subset \mathbb{T}:\text{finite}}$ such that only finitely many coefficients $a_F$ are nonzero. We use the convention $\omega_\emptyset =1$.

Then,
\begin{align*}
\mathbb E[\omega_{F_1}\cdots\omega_{F_k}]
=
\begin{cases}
1,
& \begin{smallmatrix}\text{if every $a\in F_1\cup\cdots\cup F_k$ belongs to }\\ \text{an even
number of the sets $F_1,\dots,F_k$,}\end{smallmatrix}\\
0,
& \text{otherwise.}
\end{cases}
\end{align*}

Now, we give the polynomial chaos expansion of partition functions. We have for $t>0$, $y\in \Z^2$, $\phi\in C_c(\R^2)$, 
\begin{align*}
&\overline{Z}_{N;t}^{\phi}(y)=q_{Nt}(\phi_N,y)\\
&+\sum_{k=1}^\infty \sum_{\bsm 1\leq n_1<\dots<n_k<Nt\\ x_1,\dots,x_k\in \Z^2\esm}\xi_{n_1,x_1}^{\beta_N} q_{n_1}(\phi_N,x_1)\left(\prod_{i=2}^{k}\xi_{n_i,x_i}^{\beta_N}q_{n_{i-1},n_i}(x_{i-1},x_i)\right)q_{n_k,Nt}(x_{k},y)\\
&=\sum_{k=0}^\infty\sum_{\begin{smallmatrix}\bm{A}\subset \mathbb{N}\times \mathbb{Z}^2\\ |\bm {A}|=k\end{smallmatrix}}a_{\bm{A}}^{(t)}(\phi,y) \xi_{\bm{A}}^{\beta_N},
\end{align*}
where for $k=0$, the inner sum consists only of $\bm A=\emptyset$  and for $k\geq 1$, the summation of $\bm{A}$ is taken over the set of  ${\bm{A}}=(A_1,\dots,A_{|\bm{A}|})$ with $A_i=(n_i,x_i)\in \mathbb{N}\times \mathbb{Z}^2$ ($1\leq n_1<\dots<n_{|\bm{A}|}<Nt$) and for such ${\bm{A}}$, we define \begin{align*}
&q_{n}(\phi_N,y)=\sum_{x\in \Z^2} \phi_N(x)q_{n}(x,y)\quad \text{for $n\geq 0$, $y\in \Z^2$}\\
&\xi_{\bm{A}}^{\beta_N}=\prod_{(n,x)\in {\bm{A}}}\xi_{n,x}^{\beta_N}
\end{align*}
and 
\begin{align*}
a_{\bm{A}}^{(t)}(\phi,y)=\begin{cases}
q_{Nt}(\phi_N,y)\quad &\text{if } {\bm{A}}=\emptyset\\ 
q_{n_1}(\phi_N,x_1)\left(\prod_{i=2}^{|{\bm{A}}|}q(A_{i-1},A_i)\right)q(A_{|{\bm{A}}|},A_{|{\bm{A}}|+1})\quad &\text{if }{\bm{A}}=(A_1,\dots,A_{|{\bm{A}}|})\\
0\quad &\text{otherwise}.
\end{cases}
\end{align*}
Here, we adopt the convention that $\prod_{i=2}^{1}(\cdots)=1$ for convenience,  and we set \begin{align*}
&A_{|\bm{A}|+1}=(Nt,y)\\
&q(A_{i-1},A_i)=q_{n_i-n_{i-1}}(x_{i-1},x_i).
\end{align*}
Hereafter, we consider only  a finite set
$\bm A\subset\mathbb{N}\times \mathbb{Z}^2$ containing at most one point
at each time, and we write
\begin{align*}
\bm A=\{A_1,\dots,A_{|\bm A|}\},
\qquad
A_i=(n_i,x_i),
\qquad
n_1<\cdots<n_{|\bm A|},
\end{align*}
for its time-ordered representation.

For our convenience, we introduce a set of finite subsets of indices \begin{align*}
&\mathtt{F}(T)=\mathtt{F}_{T}:=\left\{\bm{A}\subset \mathbb{N}\times \mathbb{Z}^2: \begin{array}{l}\bm{A}=(A_1,\dots,A_{|\bm{A}|}),\quad A_i=(n_i,x_i)\in \mathbb{N}\times \Z^2, \\
1\leq n_1<\dots<n_{|\bm{A}|}<NT\end{array}\right\}
\end{align*}
for $T\geq 0$, where we set $|\bm{A}|=0$ for $\bm{A} =\emptyset$. We define $\dis \mathtt{F}=\bigcup_{T\geq 0}\mathtt{F}_{T}$.

Now, we fix $T>0$. The above argument gives that for any $0\leq s\leq t\leq T$\begin{align}
&\IE\left[\overline{Z}_{N;s}^{\phi}(y^a)\overline{Z}_{N;s}^{\phi}(y^b)\overline{Z}_{N;t}^{\phi}(y^c)\overline{Z}_{N;t}^{\phi}(y^d)\right]\label{eq:partmomentlevel1}\\
&=\sum_{\bsm \bm{A},\bm{B}\in \mathtt{F}({s})\esm}\sum_{\bsm \bm{C},\bm{D}\in \mathtt{F}(t)\esm} a_{\bm{A}}^{(s)}(\phi,y^a)a_{\bm{B}}^{(s)}(\phi,y^b)a_{\bm{C}}^{(t)}(\phi,y^c)a_{\bm{D}}^{(t)}(\phi,y^d) \IE\left[\xi_{\bm{A}}^{\beta_N} \xi_{\bm{B}}^{\b_N}\xi_{\bm{C}}^{\b_N}\xi_{\bm{D}}^{\b_N}\right].\notag
\end{align}

\begin{remark}
We remark that this representation of the fourth moment of partition function can be regarded as the summation of the product of transition kernels of random walks starting from $y^a,y^b,y^c,y^d$.

\end{remark}

Now, we  focus on the finite subsets $\bm{A},\bm{B},\bm{C},\bm{D} \in \mathtt{F}$ that contribute to the summation in the right-hand side. 

We say that
$\bm A,\bm B,\bm C,\bm D$
have the even-intersection property if every point in their union
belongs to an even number of these four sets.

Then, we have from \eqref{eq:omegamoment} \begin{align*}
{\small \IE\left[ \xi_{\bm{A}}^{\b_N} \xi_{\bm{B}}^{\b_N} \xi_{\bm{C}}^{\b_N} \xi_{\bm{D}}^{\b_N}\right]=\begin{cases}
\sigma_N^{|\bm{A}|+|\bm{B}|+|\bm{C}|+|\bm{D}|}\,&\text{if $\bm A,\bm B,\bm C,\bm D$ have the
even-intersection property,}\\
0\quad &\text{otherwise}\\
\end{cases}}
\end{align*}

Thus, $\bm{A},\bm{B},\bm{C},\bm{D}$ that have even intersections contribute to the summation of $\bm{A},\bm{B},\bm{C},\bm{D}$ in \eqref{eq:partmomentlevel1}.

\begin{remark}
  This implies that the associated random walks meet at each time $n$ such that $(n,x)$ belongs to $\bm{A}\cup \bm{B}\cup \bm{C}\cup \bm{D}$, and the number of random walks that meet at time $n$ is even. 
\end{remark}

The above argument yields that \begin{align}
 &\eqref{eq:4thmomentQV}\label{eq:partmoments3}
=\frac{\sigma_N^4}{N^4}\sum_{\bsm Ns\leq Nu\leq Nt\\ Ns\leq Nv\leq Nt\esm}\sum_{y^1,y^2}\sum_{\bsm \bm{A},\bm{B}\in \mathtt{F}(u)\\\bm{C},\bm{D}\in \mathtt{F}(v)\\ \text{even intersections}\esm}\\
&\hspace{8em}\sigma_N^{|\bm{A}|+\dots+|\bm{D}|}
a^{(u)}_{\bm{A}}(\phi,y^1)a^{(u)}_{\bm{B}}(\phi,y^1)a^{(v)}_{\bm{C}}(\phi,y^2)a^{(v)}_{\bm{D}}(\phi,y^2).\notag
\end{align}
We can write \eqref{eq:QVmoments1}-\eqref{eq:QVmoments3} in similar ways, but we omit giving them here.

Hereafter, we assume that
$\bm A,\bm B,\bm C,\bm D$
have the even-intersection property.

We write the elements of
$\bm A\cup\bm B\cup\bm C\cup\bm D$
in chronological order as
\begin{align*}
\bm A\cup\bm B\cup\bm C\cup\bm D
=
\{(n_i,x_i):1\leq i\leq k\},
\qquad
n_1\leq\cdots\leq n_k,
\end{align*}
where for the case $n_i=n_{i+1}$, we may choose $x_i\ne x_{i+1}$ such that $(n_i,x_i)\in \bm{A}$ 
since  if $n_i=n_{i+1}$ for some $i$, then $(n_i,x_i)$ belong to two of $\bm{A},\bm{B},\bm{C},\bm{D}$ and  $(n_{i+1},x_{i+1})$ belongs to the other two.  

We denote the distinct intersection times by
\begin{align*}
\left.
(\bm A\cup\bm B\cup\bm C\cup\bm D)
\right|_{\mathbb N}
:=
\{n_1,\dots,n_k\}
=
\{m_1<\cdots<m_l\}.
\end{align*}

For each intersection time
$m\in\{m_1,\dots,m_l\}$, we have two cases:
\begin{enumerate}[label=(\arabic*)]
  \item\label{item:pairs1} there exist two distinct sets
$\bm E,\bm F\in \{\bm{A},\bm{B},\bm{C},\bm{D}\}$ such that $\bm{E}\cap \bm{F}\cap \{(m,x):x\in \Z^2\}\ne \emptyset$,
   while $\bm{E}'\cap \{(m,x):x\in \Z^2\}=\emptyset$ for every $\bm{E}'\in \{\bm{A},\bm{B},\bm{C},\bm{D}\}\backslash \{\bm{E},\bm{F}\}$,
  \item\label{item:pairs2} there exist two pairs 
$(\bm E,\bm F)$ and $(\bm G,\bm H)$ such that
$\{\bm E,\bm F,\bm G,\bm H\}
=
\{\bm A,\bm B,\bm C,\bm D\},$
and
$\bm E\cap\bm F\cap\{(m,x):x\in\Z^2\}\ne\emptyset,$
$\bm G\cap\bm H\cap\{(m,x):x\in\Z^2\}\ne\emptyset.$
\end{enumerate}

We call the case \ref{item:pairs1} (resp.~\ref{item:pairs2}) one pair (resp.~two pairs) for $m$.

We denote by $[EF]$ the one pair for $m$ if the case \ref{item:pairs1} holds for $\bm{E},\bm{F}$, and by $[EF][GH]$ the two pairs for $m$ if the case \ref{item:pairs2} holds for $\bm{E},\bm{F},\bm{G},\bm{H}$.

\begin{remark}\label{rem:quadruple}
We note that in case \ref{item:pairs2}, we may have $\bm{A}\cap \bm{B}\cap \bm{C}\cap \bm{D}\cap \{(m,x):x\in \Z^2\}\not =\emptyset$.

  However, the contribution from $\bm{A},\bm{B},\bm{C},\bm{D}$ with $\bm{A}\cap\bm{B}\cap \bm{C}\cap \bm{D}\ne \emptyset$ is negligible. 
  Indeed, it is bounded above by the sum of the quadruple-intersection terms in
the expansion of $\|\psi\|_\infty^4\IE\left[\sZ_{N;t}^{\phi}(1)^4\right]$. 
The proof of \cite[Theorem 6.1]{CSZ23} shows that this contribution
is negligible; see the argument following \cite[Proposition 6.6]{CSZ23}.

\end{remark}

Then, we can write \eqref{eq:partmoments3} as 
\begin{align*}
 &\frac{\sigma_N^4}{N^4}\sum_{\bsm Ns\leq Nu\leq Nt\\ Ns\leq Nv\leq Nt\esm}\sum_{y^1,y^2}
 \hspace{-4em}\sum_{\bsm \bm{A},\bm{B}\in \mathtt{F}(u)\\\bm{C},\bm{D}\in \mathtt{F}(v)\\ \text{even intersections}\\ \text{only one pair exists at each intersection time}\esm}
 \hspace{-3em}\sigma_N^{|\bm{A}|+\dots+|\bm{D}|}
  a^{(u)}_{\bm{A}}(\phi,y^1)a^{(u)}_{\bm{B}}(\phi,y^1)a^{(v)}_{\bm{C}}(\phi,y^2)a^{(v)}_{\bm{D}}(\phi,y^2)\\
 &+\frac{\sigma_N^4}{N^4}\sum_{\bsm Ns\leq Nu\leq Nt\\ Ns\leq Nv\leq Nt\esm}\sum_{y^1,y^2}
 \hspace{-4em}\sum_{\bsm \bm{A},\bm{B}\in \mathtt{F}(u)\\\bm{C},\bm{D}\in \mathtt{F}(v)\\ \text{even intersections}\\ \text{two pairs occur at some intersection time}\esm}
 \hspace{-3em}\sigma_N^{|\bm{A}|+\dots+|\bm{D}|}
  a^{(u)}_{\bm{A}}(\phi,y^1)a^{(u)}_{\bm{B}}(\phi,y^1)a^{(v)}_{\bm{C}}(\phi,y^2)a^{(v)}_{\bm{D}}(\phi,y^2)\\
&=:\mathcal{I}_1^{(N)}+\mathcal{I}_2^{(N)}.
\end{align*}

In the next subsection, we will look at the $\mathcal{I}_1^{(N)}$ in detail and give an upper bound for it. 

The term $\mathcal I_2^{(N)}$ will be estimated in
Section \ref{subsubsec:I2bound}, where it is shown to be negligible.

\subsection{Upper bound for $\mathcal{I}_1^{(N)}$}\label{subsec:upperboundI1}

For $\mathcal{I}_1^{(N)}$, we can write \begin{align*}
  \bm{A}\cup \bm{B}\cup \bm{C}\cup \bm{D}|_\N=\{(n_1,\dots,n_k):1\leq n_1<\dots<n_k\}
\end{align*}
for some $k\in \mathbb{N}$ and $n_1,\dots,n_k\in \N$. 
We now define the collection
$\mathcal P(\bm A,\bm B,\bm C,\bm D)$
of intersection pairs as follows. For each $n_i$, we have one pair for $n_i$ since we are considering $\mathcal{I}_1^{(N)}$. 
Then, we obtain two of $\bm{A},\bm{B},\bm{C},\bm{D}$ that have an intersection at time $n_i$. 
We call this pair of sets a pair for $n_i$. 

Thus, we obtain an associated pairing sequence, for example,
$([AB],[BC],\dots,[CD],[AB]).$

Next, we consider the summation over the space variables in $a^{(u)}_{\bm{A}}(\phi,y^1)a^{(u)}_{\bm{B}}(\phi,y^1)a^{(v)}_{\bm{C}}(\phi,y^2)a^{(v)}_{\bm{D}}(\phi,y^2)$ associated with the consecutive pairs 
in $\mathcal{P}(\bm{A},\bm{B},\bm{C},\bm{D})$, e.g.~$(\dots,[AB],[AB],\dots)$. 
Recalling the definition of $U_N(n,x)$ in \eqref{eq:UN}, we find that the summation over the space variables associated with the consecutive pairs  is equal to $U_N(n_j-n_{i},x_j-x_{i})$ for some $i<j$.

Reorganizing the sum over
$\bm A,\bm B,\bm C,\bm D$
according to their associated pairing sequences and using the notation introduced in \cite[(2.6)]{LZ26}, we obtain \begin{align}
  \mathcal{I}_1^{(N)}&=\frac{1}{N^4}
  \sum_{k=2}^\infty \sum_{\bsm [i_1,j_1],\dots,[i_k,j_k]\in \mathcal{P}(\bm{A},\bm{B},\bm{C},\bm{D})\\ [i_r,j_r]\ne [i_{r+1},j_{r+1}] \text{ for $r=1,\dots,k-1$}\\ [i_{k-1},j_{k-1}],[i_k,j_k]\in \{[AB],[CD]\}\esm}
  \sum_{\bsm 1\leq m_1<n_1<m_2<\dots<m_k<n_k\leq Nt\\ Ns <n_{k-1}<n_k\leq Nt\\ x_1,\dots,x_k,y_1,\dots,y_k\in \Z^2\esm}\label{eq:partmoment34}\\
 & 
 \prod_{r=1}^k U_{m_{r},n_r}(x_{r},y_r)\psi_N(y_{k-1})^2\psi_N(y_k)^2\notag\\
 &\times \left(\prod_{1\leq r\leq k-1} q_{m_{r+1}-n_{\mathtt{p}(i_{r+1})}}(y_{\mathtt{p}(i_{r+1})},x_{r+1})q_{m_{r+1}-n_{\mathtt{p}(j_{r+1})}}(y_{\mathtt{p}(j_{r+1})},x_{r+1})\right)\notag,
\end{align}
and for a pair $[i_r,j_r]$, we set \begin{align*}
  \mathtt{p}(i_r):=\max\{\ell\in\{1,\dots,r-1\}:i_r\in\{i_\ell,j_\ell\}\}
\end{align*}
if the set on the right-hand side is nonempty, and set
$\mathtt p(i_r):=0$ otherwise.
We define $\mathtt p(j_r)$ analogously.
In other words, $\mathtt{p}(i_r)$ (resp.~$\mathtt{p}(j_r)$) is the index of the last pair before $[i_r,j_r]$ that contains $i_r$ (resp.~$j_r$) if such a pair exists, and otherwise, we set $\mathtt{p}(i_r)=0$ (resp.~$\mathtt{p}(j_r)=0$).
Furthermore, we set $m_0=n_0=0$ and use the convention
\begin{align*}
q_{m_r-n_0}(y_0,x_r)
:=
q_{m_r}(\phi_N,x_r)
\end{align*}
whenever the predecessor index is $0$.

\begin{remark}
By construction, $\mathtt{p}(i_k)\leq k-1$ and $\mathtt{p}(j_k)\leq k-1$ ($k\geq 1$) and one of them is at most $k-2$ for $k\geq 2$. 
Moreover, we note that  if $\mathtt{p}(i_k)\leq k-2$ and $\mathtt{p}(j_k)\leq k-2$, then one of them is $k-2$ since $[i_{k-1},j_{k-1}]\ne [i_k,j_k]$.
 Indeed, in this case the pairs
$[i_{k-1},j_{k-1}]$ and $[i_k,j_k]$ are complementary, while
$[i_{k-2},j_{k-2}]\ne[i_{k-1},j_{k-1}]$; hence
$[i_{k-2},j_{k-2}]$ contains either $i_k$ or $j_k$.

We do not distinguish $[i_k,j_k]$ from $[j_k,i_k]$.
For later convenience, however, we choose the ordering so that $\mathtt{p}(i_k)\leq \mathtt{p}(j_k)$.
\end{remark}

For simplicity of notation, we write 
\begin{align*}
\mathcal{P}_{f,AB,CD}:=\left\{
\mathtt w=([i_1,j_1],\dots,[i_k,j_k]):
  \begin{array}{l}
k\geq 2 \\    
{}[i_r,j_r]\in  \mathcal{P}(\bm{A},\bm{B},\bm{C},\bm{D})\text{ for $r=1,\dots,k$}\\
{} [i_r,j_r]\ne [i_{r+1},j_{r+1}] \text{ for $r=1,\dots,k-1$}\\ 
{}\{[i_{k-1},j_{k-1}],[i_k,j_k]\}=\{[AB],[CD]\}
  \end{array}
\right\}
\end{align*}
the set of pairing sequences whose last two pairings are $[AB]$ and $[CD]$, in either order. 
For $\mathtt w=([i_1,j_1],\dots,[i_k,j_k])
\in\mathcal P_{f,AB,CD}$, we write $|\mathtt{w}|=k$.

\vspace{1em}

Now, we give an upper bound of \eqref{eq:partmoment34} by taking summation in spatial variables $x_i,y_i$.

For the case $\psi\in C_b(\R^2)$, \eqref{eq:partmoment34} is bounded  by the expression obtained by replacing $\psi_N(y_{|\mathtt{w}|-1})^2\psi_N(y_{|\mathtt{w}|})^2$ by $\|\psi\|_\infty^4$.

First, we take summation in the order of $y_{k}$$\rightarrow$$x_{k}$$\rightarrow$$y_{k-1}$$\rightarrow\dots$ as follows: 
We remark that $x_i$ ($1\leq i\leq k$) appears once in a factor $U$ and twice in factors $q$, and the same holds for $y_i$ ($1\leq i\leq k-2$).

The summand of \eqref{eq:partmoment34} has the form \[
F_k(\bm{m},\bm{n},x_1,\dots,x_{k},y_1,\dots,y_{k-1})U_{m_{k},n_{k}}(x_{k},y_{k}),
\]
where we set $\bm{m}=(m_1,\dots,m_k)$, $\bm{n}=(n_1,\dots,n_k)$, 
and hence the summation of \eqref{eq:partmoment34} in $y_{k}$   is dominated by 
\[
 F_k(\bm{m},\bm{n},x_1,\dots,x_{k},y_1,\dots,y_{k-1})U_{m_{k},n_{k}}.
 \] 
We note that   $x_{k}$ appears as $q_{m_k-n_{\mathtt{p}(i_k)}}(y_{\mathtt{p}(i_k)},x_k)q_{m_k-n_{\mathtt{p}(j_k)}}(y_{\mathtt{p}(j_k)},x_k)$ in $F(\bm{m},\bm{n},x_1,\dots,x_{k},y_1,\dots,y_{k-1})$.

We know that for $\ell \geq 1$,\begin{align}
\sum_{x_{\ell}}q_{m_\ell-n_{\mathtt{p}(i_\ell)}}(y_{\mathtt{p}(i_\ell)},x_\ell)q_{m_\ell-n_{\mathtt{p}(j_\ell)}}(y_{\mathtt{p}(j_\ell)},x_\ell)&=q_{2m_{\ell}-n_{\mathtt{p}(i_\ell)}-n_{\mathtt{p}(j_\ell)}}(y_{\mathtt{p}(i_\ell)},y_{\mathtt{p}(j_\ell)})\label{eq:QQconv1}\\
&\leq \frac{C_{q,1}}{2m_{\ell}-n_{\mathtt{p}(i_\ell)}-n_{\mathtt{p}(j_\ell)}}\notag
\end{align}
if $\mathtt{p}(j_\ell)\geq \mathtt{p}(i_\ell)\geq 1$, where $C_{q,1}$ is independent of
$m_\ell,n_{\mathtt p(i_\ell)},n_{\mathtt p(j_\ell)}$
and the spatial variables, 
\begin{align}
\sum_{x_{\ell}}q_{m_\ell-n_{\mathtt{p}(i_\ell)}}(y_{\mathtt{p}(i_\ell)},x_\ell)q_{m_\ell-n_{\mathtt{p}(j_\ell)}}(y_{\mathtt{p}(j_\ell)},x_\ell)&=q_{2m_{\ell}-n_{\mathtt{p}(j_\ell)}}(\phi_N,y_{\mathtt{p}(j_\ell)})\leq C_{q,2}\label{eq:QQconv2}
\end{align}
 if $\mathtt{p}(i_\ell)=0<\mathtt{p}(j_\ell)$, where $C_{q,2}$ is a constant  depending only on $\int \phi (x)\dd x$, and \begin{align}
\sum_{x_{\ell}}q_{m_\ell-n_{\mathtt{p}(i_\ell)}}(y_{\mathtt{p}(i_\ell)},x_\ell)q_{m_\ell-n_{\mathtt{p}(j_\ell)}}(y_{\mathtt{p}(j_\ell)},x_\ell)&=\sum_{y\in \Z^2}q_{2m_\ell}(\phi_N,y)\phi_N(y)\leq C_{q,3}N\label{eq:QQconv3}
\end{align}
 if $\mathtt{p}(i_\ell)=\mathtt{p}(j_\ell)=0$, where $C_{q,3}$ is a constant depending only on $\int \phi (x)\dd x$. 
 Thus, these are uniformly bounded in $x_1,\dots,x_k$.
 We denote by \begin{align*}
 \widetilde{q}(\ell)=\begin{cases}
 \frac{C_{q,1}}{2m_\ell-n_{\mathtt{p}(i_\ell)}-n_{\mathtt{p}(j_\ell)}}\quad &\text{if }  \mathtt{p}(j_\ell)\geq \mathtt{p}(i_\ell)\geq 1\\
 C_{q,2}&\text{if }0=\mathtt{p}(i_\ell)<\mathtt{p}(j_\ell)\\
 C_{q,3}N&\text{if }\mathtt{p}(i_\ell)=\mathtt{p}(j_\ell)=0.
 \end{cases}
 \end{align*}

Therefore, it follows that \begin{align*}
  &\sum_{x_k}F_k(\bm{m},\bm{n},x_1,\dots,x_{k},y_1,\dots,y_{k-1})U_{m_{k},n_{k}}\\
  &\leq U_{m_{k},n_{k}}\widetilde{q}(k) F_{k-1}(\bm{m},\bm{n},x_1,\dots,x_{k-1},y_1,\dots,y_{k-2})U_{m_{k-1},n_{k-1}}(x_{k-1},y_{k-1}).
\end{align*}

 By induction, we can obtain an upper bound of \eqref{eq:partmoment34}.

To give it, we divide $\widetilde{\mathcal{P}}_{f,AB,CD}$  into two disjoint sets \begin{align*}
&\widetilde{\mathcal{P}}_{\alpha}=\{\mathtt{w}=([i_1,j_1],[i_2,j_2],\dots,[i_k,j_k])\in {\mathcal{P}}_{f,AB,CD}:\text{$\{i_1,j_1,i_2,j_2\}\ne \{A,B,C,D\}$}\},\\
&\widetilde{\mathcal{P}}_{\beta}=\{\mathtt{w}=([i_1,j_1],[i_2,j_2],\dots,[i_k,j_k])\in {\mathcal{P}}_{f,AB,CD}:\text{$\{i_1,j_1,i_2,j_2\}=\{A,B,C,D\}$}\}.
\end{align*} 
For $\mathtt{w}\in \widetilde{\mathcal{P}}_{\alpha}$, the spatial summation requires one application of
\eqref{eq:QQconv3} and two applications of
\eqref{eq:QQconv2}.
For
$\mathtt w\in\widetilde{\mathcal P}_\beta$,
it requires two applications of \eqref{eq:QQconv3} and no
application of \eqref{eq:QQconv2}.

We define by $l(\mathtt{w})=\inf \{1\leq \ell\leq |\mathtt{w}|:A,B,C,D\in \{i_1,j_1,\dots,i_\ell,j_\ell\}\}$ 
the first index by which every label
$A,B,C,D$ has appeared in at least one pairing.

\begin{figure}[t]
    \begin{tabular}{cc}
      \begin{minipage}[t]{0.45\hsize}
        \centering
        \includegraphics[keepaspectratio, scale=0.22]{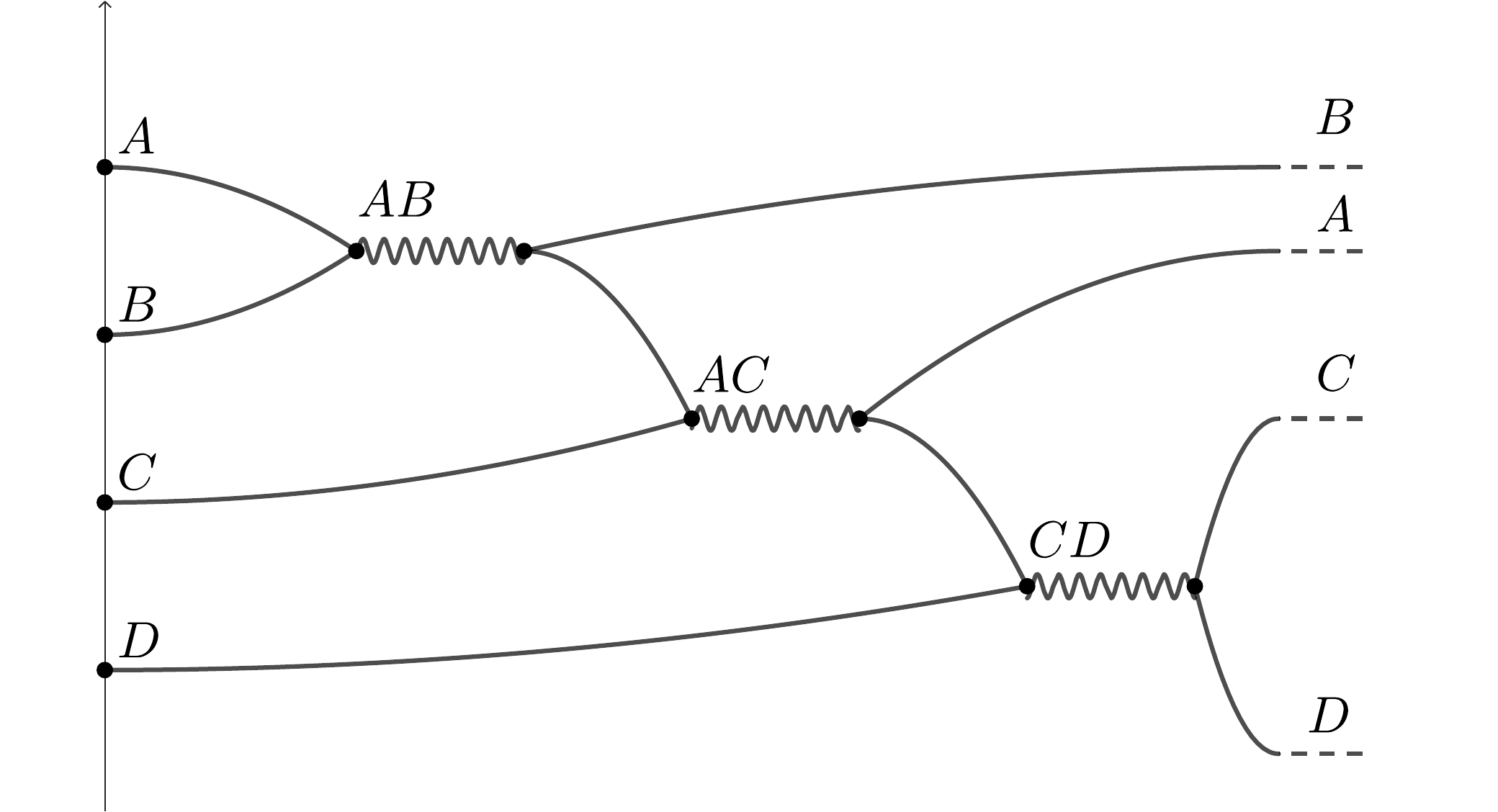}
        \caption{Pairing structure corresponding to the Type-$\alpha$ contribution
in \eqref{eq:partmoment35}.}
        \label{fig1}
      \end{minipage} &
      \begin{minipage}[t]{0.45\hsize}
        \centering
        \includegraphics[keepaspectratio, scale=0.22]{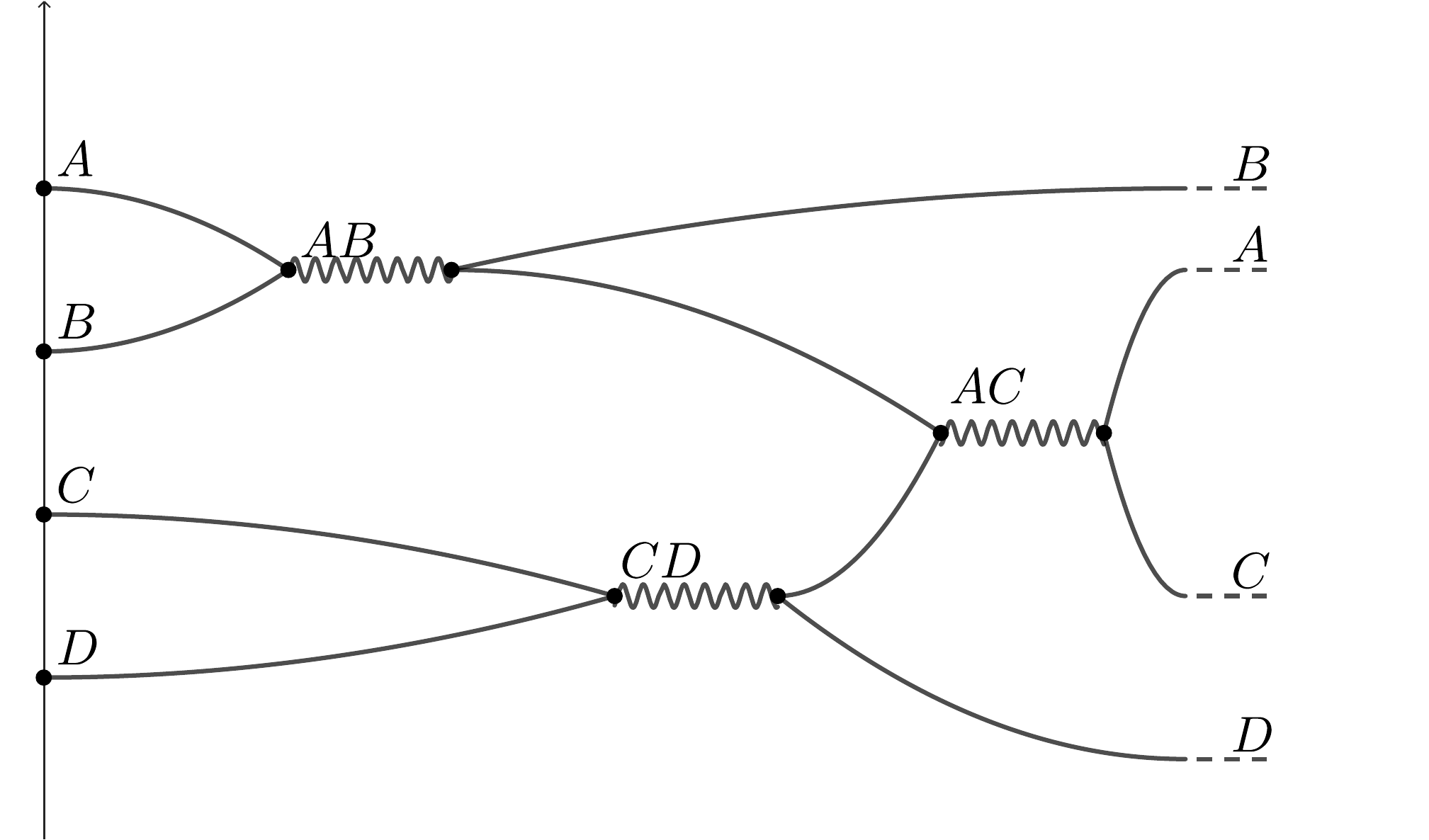}
        \caption{Pairing structure corresponding to the Type-$\beta$ contribution
in \eqref{eq:partmoment36}.}
        \label{fig2}
      \end{minipage}
    \end{tabular}
  \end{figure}

Thus, we can find that \eqref{eq:partmoment34} is dominated  by \begin{align}
&\frac{2C_{q,2}C_{q,3}}{N^3}\sum_{\mathtt{w}\in \widetilde{\mathcal{P}}_{\alpha}}
\sum_{\bsm 1\leq m_1<n_1<m_2<\dots<m_{|\mathtt{w}|}<n_{|\mathtt{w}|}\leq Nt\\ Ns <n_{|\mathtt{w}|-1}<n_{|\mathtt{w}|}\leq Nt\esm}
\prod_{\bsm j=3,\dots,|\mathtt{w}|\\  j\ne l(\mathtt{w})\esm} \widetilde{q}(j) \prod_{j=1}^{|\mathtt{w}|}U_{m_j,n_j}\notag\\
&+\frac{2C_{q,3}^2}{N^2}\sum_{\mathtt{w}\in \widetilde{\mathcal{P}}_{\beta}}
\sum_{\bsm 1\leq m_1<n_1<m_2<\dots<m_{|\mathtt{w}|}<n_{|\mathtt{w}|}\leq Nt\\ Ns <n_{|\mathtt{w}|-1}<n_{|\mathtt{w}|}\leq Nt\esm}
\prod_{\bsm j=3,\dots,|\mathtt{w}|\esm} \widetilde{q}(j) \prod_{j=1}^{|\mathtt{w}|}U_{m_j,n_j}\notag\\
&\leq \frac{2C_{q,2}C_{q,3}e^{2T\lambda}}{N^3}\sum_{\mathtt{w}\in \widetilde{\mathcal{P}}_{\alpha}}
\sum_{\bsm 1\leq m_1<n_1<m_2<\dots<m_{|\mathtt{w}|}<n_{|\mathtt{w}|}\leq Nt\\ Ns <n_{|\mathtt{w}|-1}<n_{|\mathtt{w}|}\leq Nt\esm}
\prod_{\bsm j=3,\dots,|\mathtt{w}|\\  j\ne l(\mathtt{w})\esm} \widetilde{q}(j)\prod_{j=1}^{|\mathtt{w}|}e^{-\lambda\frac{n_j-m_j}{N}}U_{m_j,n_j}\label{eq:partmoment35}\tag{Type-$\alpha$}\\
&+\frac{2C_{q,3}^2e^{2T\lambda}}{N^2}\sum_{\mathtt{w}\in \widetilde{\mathcal{P}}_{\beta}}
\sum_{\bsm 1\leq m_1<n_1<m_2<\dots<m_{|\mathtt{w}|}<n_{|\mathtt{w}|}\leq Nt\\ Ns <n_{|\mathtt{w}|-1}<n_{|\mathtt{w}|}\leq Nt\esm}
\prod_{\bsm j=3,\dots,|\mathtt{w}|\esm} \widetilde{q}(j)\prod_{j=1}^{|\mathtt{w}|}e^{-\lambda\frac{n_j-m_j}{N}}U_{m_j,n_j}.\label{eq:partmoment36}\tag{Type-$\beta$}
\end{align}
Here   $\lambda>0$ is a constant (chosen later) which will play the same role as the one introduced in \cite{CSZ19b}.

For \eqref{eq:partmoment35}, 
we have four cases
\begin{enumerate}[label=($\alpha$-\arabic*)]
\item\label{A-1} $m_1<m_2<m_{l(\mathtt{w})}\leq Ns$
\item\label{A-2} $m_1<m_2\leq Ns<m_{l(\mathtt{w})}$
\item\label{A-3} $m_1\leq Ns<m_2<m_{l(\mathtt{w})}$
\item\label{A-4} $Ns<m_1<m_2<m_{l(\mathtt{w})}$
\end{enumerate}
and  for \eqref{eq:partmoment36}, we have three cases
\begin{enumerate}[label=($\beta$-\arabic*)]
\item\label{B-1} $m_1<m_2\leq Ns$
\item\label{B-2} $m_1\leq Ns<m_2$
\item\label{B-3} $Ns<m_1<m_2$.
\end{enumerate}

Also, we split the summation according to the position of $Ns$ relative to the intervals $[m_i,n_i]$.
\begin{definition}
For each sequence $(\bm{m},\bm{n})$, 
there exists a unique index
$i=i(\bm m,\bm n)\in\{1,\dots,|\mathtt w|-1\}$ such that one of the following holds: \begin{enumerate}[label=$\mathrm{(s}$-$\mathrm{\arabic*)}$]
\item\label{item:s-1} $m_{i(\bm{m},\bm{n})} \leq Ns< n_{i(\bm{m},\bm{n})}$	
\item\label{item:s-2} $ n_{i(\bm{m},\bm{n})-1}\leq Ns<m_{i(\bm{m},\bm{n})}$.
\end{enumerate} 
\end{definition}

This gives fourteen cases in total.
Since they are treated by the same changes of variables and
convolution estimates, we present the details only for \ref{item:s-1}$\times$\ref{A-1}, and \ref{item:s-2}$\times$\ref{B-1}.

\section{Proof of Lemma \ref{lem:4thQV} }\label{sec:moments}
In this section, we complete the proof of Lemma \ref{lem:4thQV}.

\subsection{Upper bound for  $\mathcal{I}_1^{(N)}$}
As mentioned at the end of the previous section, we give  upper bounds of \ref{item:s-1}$\times$\ref{A-1} 
and \ref{item:s-2}$\times$\ref{B-1}.
The remaining cases are treated in the same way: only the ranges of
the variables before and after $Ns$ change, while the convolution
estimates and the resulting powers of $t-s$ remain unchanged.

To derive moment bounds, we use estimates from \cite{CSZ19b},
where third moments are bounded in terms of multivariate integrals.
The method is essentially the same as the one in \cite{CSZ19b} but our integrands are complicated.

\begin{remark}
  For $\phi\in C_b(\R^2)$ and $\psi\in C_c^2(\R^2)$, we take the summation in the order of $x_1\to y_1\to x_2\to \dots$ and then a similar argument holds. 
\end{remark}

\subsubsection{\ref{item:s-1}$\times $\ref{A-1}-case}

Fix $i(\bm m,\bm n)=i$. We make the following change of time
variables:
\begin{align*}\begin{cases}
u_k=n_k-m_k &\text{ for $1\leq k\leq |\mathtt{w}|$, $k\ne i$ }\\
u_i=Ns-m_i\\
\widetilde{u}_i=n_i-Ns\\
v_k=m_{k+1}-n_{k}		&\text{ for $1\leq k\leq |\mathtt{w}|-1$ }.
\end{cases}
\end{align*}

Then, we change variables from $(\bm m,\bm n)$ to
$(\bm u,\bm v,\widetilde u_i)$ and enlarge the summation domain as
follows:
 \begin{align*}
D(i,\mathtt{w}):=\left\{\begin{matrix}u_1,\dots,u_i,v_1,\dots,v_{i-1}\in \{0,\dots,NT\}\\
\widetilde{u}_i,u_{i+1},\dots,u_{|\mathtt{w}|},v_{i},\dots,v_{|\mathtt{w}|-1}\in \{0,\dots,(t-s)N\}\end{matrix}\right\}.
\end{align*}

Since $\mathtt p(i_r)\leq r-2$ and $\mathtt p(j_r)\leq r-1$, it follows that \begin{align}
&2m_r-n_{\mathtt p(i_r)}-n_{\mathtt p(j_r)}\geq 2m_r-n_{r-1}-n_{r-2}\geq 2v_{r-1}+v_{r-2}\label{eq:mnuv}
\intertext{and}
\notag&U_{n_r-m_r}=\begin{cases}
U_{u_r}\quad &\textrm{for $1\leq r\leq |\mathtt{w}|$, $r\ne i$}\\
U_{u_i+\widetilde{u}_i}\quad &\textrm{for $r=i$}.
\end{cases}
\end{align}

Thus, we have
\begin{align*}
&\textrm{\ref{item:s-1}$\times$\ref{A-1}}\\
&\leq \frac{C_{q,2}C_{q,3}e^{2T\lambda}}{N^3}\sum_{k=4}^\infty \sum_{l=3}^{k-1}
\sum_{\bsm \mathtt{w}\in \widetilde{\mathcal{P}}_{\alpha}\\|\mathtt{w}|=k\\l(\mathtt{w})=l \esm}
\sum_{i=l(\mathtt{w})}^{k-1}
\sum_{\bm{u},\widetilde{u},\bm{v}\in D(i,\mathtt{w})}\\
&\hspace{5em}\prod_{\bsm j=3\\ j\ne l(\mathtt{w})\esm}^{k} \left(\frac{C_{q,1}}{2v_{j-1}+v_{j-2}}\wedge 1\right)
\left(\prod_{\bsm j=1,\dots,k\\ j\ne  i\esm}e^{-\lambda\frac{u_j}{N}}U_{u_j}\right)e^{-\lambda\frac{\widetilde{u}_{i}+u_{i}}{N}}U_{\widetilde{u}_{i}+u_{i}},
\end{align*}
where $|\mathtt{w}|\geq 4$ since $|\mathtt{w}|>l(\mathtt{w})>2$ for \ref{A-1}.

We use the following result.

\begin{lemma}\label{lem:Uintgral}\cite[Lemma 5.3 and (5.37)]{CSZ19b}
For each $\lambda\geq 1$ and $T\geq 1$, there exists a constant $\mathtt{c}_\vartheta<\infty$ such that for any $N\geq 1$\begin{align*}
\sum_{u=1}^{NT} e^{-\lambda \frac{u}{N}}U_u\leq \frac{\mathtt{c}_\vartheta}{2+\log \lambda}.
\end{align*}
\end{lemma}

\begin{corollary}\label{cor:Uint}
For each $T\geq 1$ and $t\geq 0$, and $\lambda\geq1$, there exists $\mathtt{c}_{\vartheta,T,t}<\infty$ such that for any $N\geq 1$\begin{align*}
\frac{1}{N}\sum_{v=0}^{\lfloor Nt\rfloor}\sum_{u=1}^{\lfloor NT\rfloor} e^{-\lambda \frac{u+v}{N}}U_{u+v}\leq \frac{\mathtt{c}_{\vartheta,T,t} \lfloor t+\frac{1}{N}\rfloor}{2+\log \lambda}.
\end{align*}
\end{corollary}

Thus, \begin{align}
&\textrm{\ref{item:s-1}$\times$\ref{A-1}}\label{eq:inequ1-1}\\
&\leq \frac{C_{q,2}C_{q,3}e^{2T\lambda}(t-s)}{N^2}\notag\\
&\qquad \sum_{k=4}^\infty  
\left(\frac{c_{\vartheta}}{2+\log \lambda}\right)^{k}
\sum_{l=3}^{k-1}\sum_{\bsm \mathtt{w}\in \widetilde{\mathcal{P}}_{\alpha}\\|\mathtt{w}|=k \\l(\mathtt{w})=l\esm}
\sum_{i=l(\mathtt{w})}^{k-1}\sum_{\bm{v}\in D_{\bm{v}}(i,\mathtt{w})}
\prod_{\bsm j=3 \\ j\ne l(\mathtt{w})\esm}^{k} \left(\frac{C_{q,1}}{2v_{j-1}+v_{j-2}}\wedge 1\right),\notag
\end{align}
where we set \begin{align}
D_{\bm{v}}(i,\mathtt{w}):=\left\{\begin{matrix}v_1,\dots,v_{i-1}\in \{0,\dots,NT\}\\
v_{i},\dots,v_{|\mathtt{w}|-1 }\in \{0,\dots,(t-s)N\}\end{matrix}\right\}.\label{eq:dvis}
\end{align}

Now, we focus on \begin{align*}
\sum_{\bm{v}\in D_{\bm{v}}(i,\mathtt{w})}\prod_{\bsm j=3\\ j\ne l(\mathtt{w})\esm}^{k}\left(\frac{C_{q,1}}{2v_{j-1}+v_{j-2}}\wedge 1\right)
\end{align*}
which is the summation of the product of $k-3$ terms with respect to $k-1$ variables. 
By the AM-GM inequality $a+b\geq 2\sqrt{ab}$ for $a\geq 0 $, $b\geq 0$, it is dominated by  
\begin{align}
&\frac{1}{(2N)^{k-3}}\sum_{\bm{v}\in D_{\bm{v}}(i,\mathtt{w})}\prod_{\bsm j=1\\ j\ne l(\mathtt{w})-2\esm}^{k-2}\left(\frac{C_{q,1}}{\sqrt{\frac{v_{j+1}}{N}}\sqrt{\frac{v_{j+1}}{N}+\frac{v_{j}}{N}}}\wedge N\right)\label{eq:integralstree1}\\
&\leq C_{q,1}^{k-3}N^2 \int_{[0,s]^{i-1}\times [0,t-s]^{k-i}}d\bm{v}\prod_{\bsm j=1\\ j\ne l(\mathtt{w})-2\esm}^{k-2}\frac{1}{\sqrt{v_{j+1}}\sqrt{v_{j+1}+v_{j}}}.\notag
\end{align}

Now, we integrate it in order from $v_1$ to $v_{k-1}$.

Since $v_1$ appears as $\frac{1}{\sqrt{v_1+v_2}}$ in integrand of \eqref{eq:integralstree1},  \begin{align*}
&\int_{[0,s]^{i-1}\times [0,t-s]^{k-i}}d\bm{v}\prod_{\bsm j=1\\ j\ne l(\mathtt{w})-2\esm}^{|\mathtt{w}|-2}\frac{1}{\sqrt{v_{j+1}}\sqrt{v_{j+1}+v_{j}}}\\
&\leq 2\sqrt{2T}\int_{[0,T]^{i-2}\times [0,t-s]^{k-i}}d\bm{v}\left(\prod_{j=2}^{l(\mathtt{w})-3}
\frac{1}{\sqrt{v_j}\sqrt{v_j+v_{j+1}}}\right)\frac{1}{\sqrt{v_{l(\mathtt{w})-2}}}\\
&\hspace{5em}\frac{1}{\sqrt{v_{l(\mathtt{w})-1}+v_{l(\mathtt{w})}}}\left(\prod_{j'=l(\mathtt{w})}^{k-2}\frac{1}{\sqrt{v_{j'}}\sqrt{v_{j'}+v_{j'+1}}}\right)\frac{1}{\sqrt{v_{k-1}}}.
\end{align*}

To estimate the integral in the right-hand side, we use the results in \cite{CSZ19b}.

We define \begin{align*}
&\phi_t^{(0)}(u):=1,\\
&\phi_t^{(k)}(u):=\int_0^t \frac{1}{\sqrt{s(s+u)}}\phi_t^{(k-1)}(s)ds,\quad \text{for $k\geq 1$, $u\geq 0$, and $t> 0$}.
\end{align*}
Then, by induction on $k$,
\begin{align*}
\phi_t^{(k)}(u)
=
\int_0^1
\frac{
\phi_t^{(k-1)}(tr)
}{
\sqrt{r(r+u/t)}
}\dd r
=
\int_0^1
\frac{
\phi_1^{(k-1)}(r)
}{
\sqrt{r(r+u/t)}
}\dd r
=
\phi_1^{(k)}\left(\frac{u}{t}\right).
\end{align*}

\begin{lemma}[{\cite[Lemma 5.4]{CSZ19b}}]
For all $k\in \mathbb{N}$, \begin{align*}
\phi^{(k)}_1(v)\leq 32^k \sum_{i=0}^k\frac{1}{i!}\left(\frac{1}{2}\log \frac{e^2}{v}\right)^i\leq 32^k\frac{e}{\sqrt{v}},
\end{align*}
for all $v\in (0,1)$.
\end{lemma}

In particular, \begin{align*}
\phi_T^{(k)}(v)\leq 32^kp^k \sum_{i=0}^k \frac{1}{i!}\left(\frac{1}{2p}\log \frac{Te^2}{v}\right)^i\leq (32p)^kT^{\frac{1}{2p}}e^\frac{1}{p}v^{-\frac{1}{2p}}
\end{align*}
for $T>0$ and $p\geq 1$.

Therefore, \begin{align}
&\int_{[0,T]^{l(\mathtt{w})-3}}\prod_{j=2}^{l(\mathtt{w})-3}\frac{1}{\sqrt{v_{j}}\sqrt{v_j+v_{j+1}}}\frac{1}{\sqrt{v_{l(\mathtt{w})-2}}}dv_2\dots dv_{l(\mathtt{w})-2}\notag\\
&\leq (32p)^{l(\mathtt{w})-4}T^\frac{1}{2p}e^\frac{1}{p}\int_0^T \frac{1}{v^{\frac{1}{2}+\frac{1}{2p}}}dv\leq \frac{1}{\frac{1}{2}-\frac{1}{2p}}(32p)^{l(\mathtt{w})-3}T^{\frac{1}{2}}e^\frac{1}{p}.\label{eq:ineqL}
\end{align} 
for $p>1$.

When $l(\mathtt w)=3$, the corresponding product is empty and the
integral is bounded directly by $C_pT^{1/2}$. Thus, after enlarging
the constant, the same form of estimate may be used for all
$l(\mathtt w)\geq3$.

Also, we have {\small \begin{align}
&\int_{[0,T]^{i-l(\mathtt{w})+1}\times [0,t-s]^{k-i}}\frac{1}{\sqrt{v_{l(\mathtt{w})-1}+v_{l(\mathtt{w})}}}\left(\prod_{j'=l(\mathtt{w})}^{k-2}\frac{1}{\sqrt{v_{j'}}\sqrt{v_{j'}+v_{j'+1}}}\right)\frac{1}{\sqrt{v_{k-1}}}dv_{l(\mathtt{w})-1}\dots dv_{k-1}\label{eq:ineqL2}\\
&\leq  2\sqrt{T}(32p)^{k-l(\mathtt{w})-1}T^\frac{1}{2p}e^\frac{1}{p}\int_0^{t-s}\frac{1}{v^{\frac{1}{2}+\frac{1}{2p}}}dv\notag\\
&\leq \frac{2}{\frac{1}{2}-\frac{1}{2p}}T^{\frac{1}{2}+\frac{1}{2p}}(32p)^{k-l(\mathtt{w})-1}e^\frac{1}{p}(t-s)^{\frac{1}{2}-\frac{1}{2p}}.\notag
\end{align}}

Combining \eqref{eq:inequ1-1}, \eqref{eq:integralstree1}, \eqref{eq:ineqL}, and \eqref{eq:ineqL2},
we obtain that \begin{align*}
\textrm{\ref{item:s-1}$\times$\ref{A-1}}&\leq C_{p}C_{q,2}C_{q,3}e^{2\lambda T}(t-s)^{\frac{3}{2}-\frac{1}{2p}}\\
&\qquad \sum_{k=4}^\infty  \left(\frac{c_{\vartheta}}{2+\log \lambda}\right)^{k}\sum_{l=3}^{k-1}
\sum_{\bsm \mathtt{w}\in \widetilde{\mathcal{P}}_{\alpha}\\|\mathtt{w}|=k \\l(\mathtt{w})=l\esm}\sum_{i=l(\mathtt{w})}^{k-1}(32p)^{k-5}T^{\frac{1}{2p}+\frac{3}{2}}\\
&\leq C_{p}C_{q,2}C_{q,3}e^{2\lambda T}(t-s)^{\frac{3}{2}-\frac{1}{2p}}\sum_{k=4}^\infty  \left(\frac{c_{\vartheta}}{2+\log \lambda}\right)^{k}k\,6^k(32p)^{k-5}T^{\frac{1}{2p}+\frac{3}{2}}\\
&\leq C_{p,T,\lambda}(t-s)^{\frac{3}{2}-\frac{1}{2p}}
\end{align*}
for $\lambda$ large enough, where $C_p$ is a constant depending only on $p>1$,  $C_{p,T,\lambda}$ is a constant depending on $p,T,\lambda$, and where we used the crude bounds
$\#\{\mathtt w:|\mathtt w|=k\}\leq6^k$
and
$\sum_i1\leq k$.

\subsubsection{\ref{item:s-2}$\times$\ref{B-1}-case}

We fix $i(\bm{m},\bm{n})=i$. 
We make the following change of time
variables:

\begin{align*}\begin{cases}
u_k=n_k-m_k,		&\text{ for $1\leq k\leq |\mathtt{w}|$ }\\
v_k=m_{k+1}-n_{k}&\text{for $1\leq k\leq |\mathtt{w}|-1$, $k\ne i-1$}\\
v_{i-1}=Ns-n_{i-1}\\
\widetilde{v}_{i-1}=m_{i}-Ns.
\end{cases}
\end{align*}

Then, we replace the variables of the summation from $(\bm{m},\bm{n})$ to $(\bm{u},\bm{v},\widetilde{v}_{i-1})$ and enlarge the range of them as follows:
\begin{align*}
D_{2}(i,\mathtt{w}):=\left\{\begin{matrix}u_1,\dots,u_{|\mathtt{w}|},v_1,\dots,v_{i-1}\in \{0,\dots,NT\}\\
\widetilde{v}_{i-1},v_{i},\dots,v_{|\mathtt{w}|-1}\in \{0,\dots,(t-s)N\}\end{matrix}\right\}.
\end{align*} 

Using \eqref{eq:mnuv}, it follows that 
\begin{align*}
&\textrm{\ref{item:s-2}$\times$\ref{B-1}}\\
&\leq \frac{C_{q,3}^2e^{2T\lambda}}{N^2}
\sum_{k=4}^\infty \sum_{\bsm \mathtt{w}\in \widetilde{\mathcal{P}}_{\beta}\\ |\mathtt{w}|=k\esm}
\sum_{i=3}^{k-1}
\sum_{\bsm (\bm{u},\bm{v},\widetilde{v}_{i-1})\in D_2(i,\mathtt{w})\esm}
\prod_{\bsm j=3,\dots,k \esm} \left(\frac{C_{q,1}}{2v_{j-1}+v_{j-2}}\wedge 1\right)
\prod_{j=1}^{k}e^{-\lambda\frac{u_j}{N}}U_{u_j}
,\notag
\end{align*}
where we remark that $i(\bm m,\bm n)\leq|\mathtt w|-1$
because $n_{|\mathtt w|-1}>Ns$. Also, we remark that the summand does not contain $\widetilde{v}_{i-1}$.

Lemma \ref{lem:Uintgral} yields that \begin{align*}
&\textrm{\ref{item:s-2}$\times$\ref{B-1}}\\
&\leq \frac{C_{q,3}^2e^{2T\lambda}(t-s)}{N}
\sum_{k=4}^\infty \left(\frac{c_{\vartheta}}{2+\log \lambda}\right)^k
\sum_{\bsm \mathtt{w}\in \widetilde{\mathcal{P}}_{\beta}\\ |\mathtt{w}|=k\esm}
\sum_{i=3}^{k-1}\sum_{\bsm \bm{v}\in D_{\bm{v}}(i,\mathtt{w})\esm}
\prod_{\bsm j=3,\dots,k \esm} \left(\frac{C_{q,1}}{2v_{j-1}+v_{j-2}}\wedge 1\right)
\end{align*}
where $D_{\bm{v}}(i,\mathtt{w})$ is defined in \eqref{eq:dvis}.

Then, a similar argument to the analysis of \eqref{eq:integralstree1} yields that \begin{align*}
&\sum_{\bsm \bm{v}\in D_{\bm{v}}(i,\mathtt{w})\esm}
\prod_{\bsm j=3,\dots,k \esm}
 \left(\frac{C_{q,1}}{2v_{j-1}+v_{j-2}}\wedge 1\right)\\
 &\qquad \leq C_{q,1}^{k-2}N\int_{[0,s]^{i-1}\times [0,t-s]^{k-i}}d\bm{v} \prod_{j=1}^{k-2}\frac{1}{\sqrt{v_{j+1}}\sqrt{v_{j+1}+v_j}}.
\end{align*}

The remaining analysis is identical to that following
\eqref{eq:integralstree1}, and we omit the details.

Consequently, we  obtain that \begin{align*}
&\textrm{\ref{item:s-2}$\times$\ref{B-1}}\leq C_{p,T,\lambda}(t-s)^{\frac{3}{2}-\frac{1}{2p}}
\end{align*}
for $\lambda$ large enough, where $C_{p,T,\lambda}$ is a constant depending on $p,T,\lambda$.

\subsubsection{Upper bound of $\mathcal{I}_2^{(N)}$}\label{subsubsec:I2bound}

To estimate $\mathcal{I}_2^{(N)}$, we use the bound of $\mathcal{I}_1^{(N)}$. 

We note that modifying the above arguments yields that \begin{align*}
  \mathcal{I}_2^{(N)} &\leq \frac{2C_{q,2}C_{q,3}e^{2T\lambda}}{N^3}\sum_{\mathtt{w}\in \widetilde{\mathcal{P}}_{\alpha}}
\sum_{ (\bm{m},\bm{n})\in D_e(T)}
\prod_{\bsm j=3,\dots,|\mathtt{w}|\\  j\ne l(\mathtt{w})\esm} \widetilde{q}(j)\prod_{j=1}^{|\mathtt{w}|}e^{-\lambda\frac{n_j-m_j}{N}}U_{m_j,n_j}\\
&+\frac{2C_{q,3}^2e^{2T\lambda}}{N^2}\sum_{\mathtt{w}\in \widetilde{\mathcal{P}}_{\beta}}
\sum_{(\bm{m},\bm{n})\in D_e(T)}
\prod_{\bsm j=3,\dots,|\mathtt{w}|\esm} \widetilde{q}(j)\prod_{j=1}^{|\mathtt{w}|}e^{-\lambda\frac{n_j-m_j}{N}}U_{m_j,n_j},
\end{align*}
where \begin{align*}
 D_e(T) =\left\{(\bm{m},\bm{n}):\bsm \bm{m}=(m_1,\dots,m_k),\bm{n}=(n_1,\dots,n_k) \text{ for $k\geq 2$},\\  1\leq m_1<n_1\leq m_2<n_2\leq\cdots\leq m_k<n_k\leq Nt,\\ \text{$n_i=m_{i+1}$ for some $1\leq i\leq k-1$}\esm\right\}.
\end{align*}

For each possible index $i$ satisfying $n_i=m_{i+1}$, we eliminate
the variable $m_{i+1}$ from the summation. The resulting sum is
bounded by the corresponding estimate for $\mathcal I_1^{(N)}$,
with one fewer independent time variable and therefore one additional
factor $N^{-1}$. The remaining analysis is identical to that following
\eqref{eq:integralstree1}, and we omit the details. 

Consequently, we obtain \begin{align*}
 \mathcal{I}_2^{(N)} \leq C_{p,T,\lambda}(t-s)^{\frac{3}{2}-\frac{1}{2p}}.
\end{align*}

Combining these
bounds with the estimate for $\mathcal I_2^{(N)}$, we obtain, for
every $p>1$,
\[
\IE\left[
\left(
\langle M^{N,\phi}(\psi)\rangle_{Nt}
-
\langle M^{N,\phi}(\psi)\rangle_{Ns}
\right)^2
\right]
\leq
C_{p,T}|t-s|^{\frac{3}{2}-\frac{1}{2p}},
\]
uniformly in $N$ and $0\leq s\leq t\leq T$.
This proves Lemma \ref{lem:4thQV}.

\section{Proof of Lemma \ref{lem:QuadProcEst}}\label{sec:momentsbound}

As noted in Remark \ref{rem:moment4th}, 
it is enough to focus on \eqref{eq:QVmoments1}-\eqref{eq:QVmoments3}. 

The third-moment limit for
$\overline Z^\phi_{N;s}(\psi)$ was established in \cite{CSZ19b},
and explicit higher-moment formulas in the continuous setting were
obtained in \cite{GQT21}.

First, we introduce some notation.

We define the set of pairing sequences on $\{A,B,C,D\}$ by\begin{align*}
\mathcal{P}_{f}:=\left\{
  \mathtt w=([i_1,j_1],\dots,[i_k,j_k]):
  \begin{array}{l}
k\geq 1,\\ 
{} [i_r,j_r]\in \mathcal{P}(\bm{A},\bm{B},\bm{C},\bm{D}), \text{ for $r=1,\dots,k$}\\
{} [i_r,j_r]\ne [i_{r+1},j_{r+1}] \text{ for $r=1,\dots,k-1$}
  \end{array}
\right\}
\end{align*}
and for each $\mathtt{w}\in  {\mathcal{P}}_f$ and $E\in \{A,B,C,D\}$, we define \begin{align*}
k^E:=\sup\{\ell:E\in [i_\ell,j_\ell]\in \mathtt{w} \}
\end{align*}
if the set is nonempty and we set $k^E=0$ if it is  empty. 
Thus, $k^E$ is the last index at which $E$ appears in
$\mathtt w$.

We also define   
\begin{align*}
\Theta_\mathtt{w}(\bm{u},\bm{v},\bm{x},\bm{y})
=\prod_{r=1}^{|\mathtt{w}|}p_{u_r-v_{\mathtt{p}(i_r)}}\left(y_{\mathtt{p}(i_r)},x_r\right)
p_{u_r-v_{\mathtt{p}(j_r)}}\left(y_{\mathtt{p}(j_r)},x_r\right)
\end{align*}
for  $\phi\in C_c(\R^2)$, $0<u_1<v_1<\dots<u_{|\mathtt w|}<v_{|\mathtt w|}$, $x_1,y_1,\dots,x_{|\mathtt w|},y_{|\mathtt w|}\in \R^2$, and $\mathtt{w}\in {\mathcal{P}}_f $, 
where we set $p_{u_r-v_0}(y_0,x_r)
:=
\Phi_{u_r}(x_r)$ for $v_0:=0$.

The limits of \eqref{eq:QVmoments1}-\eqref{eq:QVmoments3} are given as follows. 
\begin{lemma}\label{lem:4thmomentlimitTT}
Let $\phi\in C_c(\R^2)$ and $\psi\in C_b^2(\R^2)$. For each $t\geq 0$  and for each $\varepsilon>0$,  \begin{align}
&\lim_{N\to \infty}\IE\left[\left(\int_0^{\frac{\lfloor Nt\rfloor}{N}}\int_{\R^2}\overline{\mathsf{Z}}^{\phi}_{N;s}(p_\e(\cdot-z))^2\psi(z)^2\dd z\dd s\right)^2\right]\label{eq:momentTT1}\\
&=\int_{[0,t]^2}\prod_{i=1}^2\left(\int\Phi_{s_i+\e}(y)^2\psi(y)^2\dd y\right)\dd s_1\dd s_2
\notag\\
&+\sum_{k\geq 1}\sum_{\bsm \mathtt{w}\in \mathcal{P}_{f}\\|\mathtt{w}|=k\esm } 
(4\pi)^k\iint_{\bsm  0<u_1<v_1<\dots<u_{k-1}<v_{k-1}<u_k<v_k<t\\ v_{k^A}\vee v_{k^B}<\sigma<t,v_{k^C}\vee v_{k^D}<\tau<t \esm}\dd\bm{u}\dd\bm{v}\dd \sigma\dd \tau
\int _{\left(\R^{2}\right)^{2k+2}}\dd \bm{x}\dd \bm{y}\dd z\dd w\notag\\
&\hspace{3em}\prod_{i=1}^kG_\vartheta(v_i-u_i,y_i-x_i)\Theta_\mathtt{w}(\bm{u},\bm{v},\bm{x},\bm{y})\psi(z)^2\psi(w)^2\notag\\
&\hspace{4em}p_{\sigma-v_{k^A}+\e}(y_{k^A},z)
p_{\sigma-v_{k^B}+\e}(y_{k^B},z)
p_{\tau-v_{k^C}+\e}(y_{k^C},w)
p_{\tau-v_{k^D}+\e}(y_{k^D},w).\notag
\end{align}
\end{lemma}

\begin{lemma}\label{lem:4thmomentlimitT1}
Let $\phi\in C_c(\R^2)$ and $\psi\in C_b^2(\R^2)$. For each $t\geq 0$  and for each $\varepsilon>0$, \begin{align}
&\lim_{N\to \infty}\IE\left[\left(\int_0^{\frac{\lfloor Nt\rfloor}{N}}\int_{\R^2}\overline{\mathsf{Z}}^{\phi}_{N;s}(p_\e(\cdot-z))^2\psi(z)^2\dd z\dd s\right)\left\langle M^{N,\phi}(\psi) \right\rangle_t\right]\label{eq:momentT1-1}\\
&=\sum_{k\geq 1}\sum_{\bsm \mathtt{w}\in \mathcal{P}_{f}\\|\mathtt{w}|=k,s_k=[C,D]\esm } (4\pi)^k\iint_{\bsm 0<u_1<v_1<\dots<u_{k-1}<v_{k-1}<u_k<v_k<t\\ v_{k-1}<\sigma<t \esm}\dd\bm{u}\dd\bm{v}\dd \sigma\int _{\left(\R^{2}\right)^{2k+1}}\dd \bm{x}\dd \bm{y}\dd z
\notag\\
&\hspace{1.5em}\prod_{i=1}^k G_\vartheta(v_i-u_i,y_i-x_i)\Theta_\mathtt{w}(\bm{u},\bm{v},\bm{x},\bm{y})\psi(y_k)^2
p_{\sigma-v_{k^A}+\e}(y_{k^A},z)p_{\sigma-v_{k^A}+\e}(y_{k^B},z)\psi(z)^2\notag\\
&+\sum_{k\geq 1}\sum_{\bsm \mathtt{w}\in \mathcal{P}_{f}\\|\mathtt{w}|=k,\\ s_{k-1}=[C,D],s_k={A,B}\esm } (4\pi)^k\iint_{\bsm 0<u_1<v_1<\dots<u_{k-1}<v_{k-1}<u_k<v_k<t\\ v_{k}<\sigma<t \esm}\dd \bm{u}\dd \bm{v}\dd \bm{\sigma}\int _{\left(\R^{2}\right)^{2k+1}}\dd \bm{x}\dd \bm{y}\dd z\notag\\
&\hspace{1.5em}\prod_{i=1}^k G_\vartheta(v_{i}-u_i,y_i-x_i)
\Theta_\mathtt{w}(\bm{u},\bm{v},\bm{x},\bm{y})p_{\sigma-v_{k}+\e}\left(y_{k},z\right)^2\psi(z)^2\psi(y_{k-1})^2.\notag
\end{align}
\end{lemma}

\begin{lemma}\label{lem:4thmomentlimit11}
Let $\phi\in C_c(\R^2)$ and $\psi\in C_b^2(\R^2)$. For each $t\geq 0$, \begin{align}
&\lim_{N\to \infty}\IE\left[\langle M^{N,\phi}(\psi)\rangle_t^2\right]\label{eq:moments111}\\
&=\sum_{k\geq 2}\sum_{\bsm \mathtt{w}\in \mathcal{P}_{f}\\|\mathtt{w}|=k,\\ s_{k-1}=[A,B],s_k=[C,D] \esm } 
(4\pi)^k\iint_{\bsm 0<u_1<v_1<\dots<u_{k-1}<v_{k-1}<u_k<v_k<t\esm}\dd \bm{u}\dd \bm{v}
\int _{\left(\R^{2}\right)^{2k}}\dd \bm{x}\dd \bm{y}\notag\\
&\hspace{3em}\prod_{i=1}^k G_\vartheta(v_{i}-u_i,y_i-x_i)
\Theta_{\mathtt{w}}(\bm{u},\bm{v},\bm{x},\bm{y})\psi(y_{k-1})^2\psi(y_{k})^2.\notag
\end{align}
\end{lemma}

We give an outline of the proof of Lemma \ref{lem:4thmomentlimitT1} and 
omit the proofs of Lemma \ref{lem:4thmomentlimitTT} and \ref{lem:4thmomentlimit11} since the argument are essentially the same.

\begin{remark}\label{rem:4thmoments}
The summations in \eqref{eq:momentTT1}-\eqref{eq:moments111} converge absolutely.
 Therefore, we have  the following alternative representation of \eqref{eq:highermomentsrepre} for $h=4$: if $\varphi_i=\varphi$ and $\psi_i=\psi$ for $i=1,\dots,4$
\begin{align}
&\eqref{eq:highermomentsrepre}\label{eq:4thmomentsrepre}\\
&=\int_{(\R^2)^4}\prod_{E\in\{A,B,C,D\}}\phi(x_E)p_t(x_E,y_E)\psi(y_E)\dd \bm{x}\dd \bm{y}\notag\\
&+\sum_{k\geq 1}\sum_{\bsm \mathtt{w}\in\widetilde{\mathcal{P}}_f\\ |\mathtt{w}|=k\esm}(4\pi)^k\idotsint\limits_{0<u_1<v_1<\dots<u_k<v_k<t}\dd \bm{u}\dd \bm{v}\int_{(\bm{R}^2)^{2k}}\dd \bm{x}\dd \bm{y}\notag\\
&\hspace{4em}\prod_{i=1}^k G_\vartheta(v_i-u_i,y_i-x_i)\Theta_\mathtt{w} (\bm{u},\bm{v},\bm{x},\bm{y})\prod_{E\in \{A,B,C,D\}}\Psi_{t-v_{k^E}}(y_{k^E}).\notag
\end{align}
When $\psi\equiv 1$ and $\varphi\geq 0$, \eqref{eq:momentTT1}-\eqref{eq:moments111} and \eqref{eq:4thmomentsrepre}  have the same representation. 
\end{remark}

\subsection{Proof of Lemma \ref{lem:4thmomentlimitT1}}\label{subsec:momentlimit}
We first give the chaos expansion of the expectation of 
\begin{align*}
E\left[\left(\int_0^{\frac{\lfloor Nt\rfloor}{N}}\int_{\R^2}\overline{\mathsf{Z}}^{\phi}_{N;s}(p_\e(\cdot-z))^2\psi(z)^2\dd z\right)\left\langle M^{N,\phi}(\psi) \right\rangle_t\right]
\end{align*}
in a similar way to the argument in Section \ref{sec:ChaosExpansion}. 
We use the label $A,B$ derived from the ``random walks" in $\int_0^{\frac{\lfloor Nt\rfloor}{N}}\int_{\R^2}\overline{\mathsf{Z}}^{\phi}_{N;s}(p_\e(\cdot-z))^2\psi(z)^2\dd z$ and $C,D$ derived from the ``random walks" in $\left\langle M^{N,\phi}(\psi) \right\rangle_t$.

We remark that after the last intersection between $A$ and $B$, each of them may meet $C$ or $D$ but after the last intersection between $C$ and $D$, neither $C$ nor $D$ will meet  other particles. Thus, $\mathtt{w}\in \mathcal{P}_f$ contributing the chaos expansion should  satisfies one of
\begin{enumerate}[label=(\arabic*)]
\item $\mathtt{w}_{|\mathtt{w}|}=[C,D]$
\item $\mathtt{w}_{|\mathtt{w}|-1}=[C,D]$ and $\mathtt{w}_{|\mathtt{w}|}=[A,B]$.
\end{enumerate}
As we have seen in Remark \ref{rem:quadruple} and Section \ref{subsubsec:I2bound}, even intersections other than ${\mathcal{P}}_f$ in the chaos expansion of moments are negligible. 
Therefore, we have the following representation of the moment. We omit its proof since it is essentially the same as the discussion in \eqref{eq:partmoment34}.
\begin{lemma}\label{lem:chaosexpansionT1}
Let $\phi\in C_c(\R^2)$ and $\psi\in C_b^2(\R^2)$. For each $t\geq 0$, we have 
\begin{align}
&\IE\left[\left(\int_0^{\frac{\lfloor Nt\rfloor}{N}}\int_{\R^2}\overline{\mathsf{Z}}^{\phi}_{N;s}(p_\e(\cdot-z))^2\psi(z)^2\dd z\dd s\right)\left\langle M^{N,\phi}(\psi) \right\rangle_t\right]\label{eq:chaosexpansionT1-1}\\
&=\frac{1}{N^5}\sum_{k\geq 1}\sum_{\bsm \mathtt{w}\in {\mathcal{P}}_f\\ |\mathtt{w}|=k,[i_k,j_k]=[CD]\esm}
\sum_{\bsm x_1,y_1\dots,x_k,y_k\esm}
\sum_{\bsm 1\leq m_1<n_1<\dots<m_{k}<n_k\leq Nt\\ n_{k-1}\leq u\leq Nt\esm}\notag\\
&\hspace{4em}\prod_{i=1}^{k}U(n_i-m_i,y_i-x_i)\widetilde{\Theta}^{(N)}(\bm{m},\bm{n},\bm{x},\bm{y})\psi_N(y_k)^2\notag\\
&\hspace{5em}\times \int_{\R^2}\dd zQ_{A}(\bm{y},u,\e,\psi,z)Q_{B}(\bm{y},u,\e,\psi,z)\notag\\
&+\frac{1}{N^5}\sum_{k\geq 1}\sum_{\bsm \mathtt{w}\in \mathcal{P}_f\\ |\mathtt{w}|=k,[i_{k-1},j_{k-1}]=[CD],[i_k,j_k]=[AB]\esm}\sum_{x_1,y_1,\dots,x_k,y_k}\sum_{\bsm 1\leq m_1<n_1<\dots<m_{k}<n_k\leq Nt\\ n_k< u\leq Nt\esm}\notag\\
&\hspace{4em}\prod_{i=1}^{k}U(n_i-m_i,y_i-x_i)\widetilde{\Theta}^{(N)}(\bm{m},\bm{n},\bm{x},\bm{y})\psi_N(y_k)^2\notag\\
&\hspace{5em}\int_{\R^2}\dd z Q_{A}(\bm{y},u,\e,\psi,z)Q_{B}(\bm{y},u,\e,\psi,z)\notag\\
&\hspace{5em} +o(1),\notag
\end{align}
where we set \begin{align*}
&\widetilde{\Theta}^{(N)}(\bm{m},\bm{n},\bm{x},\bm{y})=\prod_{r=1}^{|\mathtt{w}|}q_{m_r-n_{\mathtt{p}(i_r)}}(y_{\mathtt{p}(i_r)},x_r)q_{m_r-n_{\mathtt{p}(j_r)}}(y_{\mathtt{p}(j_r)},x_r)\\
&Q_E(\bm{y},u,\e,\psi,z)=\sum_{w\in \Z^2}q_{u-n_{k^E}}(y_{k^E},w)p_{\e}\left(\frac{w}{\sqrt{N}}-z\right)\psi(z)\quad \text{for $E=A,B$},
\end{align*}
and when $k^E=0$, we use the convention
\begin{align*}
n_0:=0,
\qquad
q_{u-n_0}(y_0,w):=
q_u(\phi_N,w).
\end{align*}
\end{lemma}

We use the labels $A$ and $B$ for the two replicas coming from the
mollified measure, while $C$ and $D$ correspond to the two replicas
coming from the quadratic variation. Accordingly, the terminal
kernel $Q_E$ is needed only for $E=A,B$.

Thus, it is enough to see the limit of the right-hand side of \eqref{eq:chaosexpansionT1-1} for the proof of Lemma \ref{lem:4thmomentlimitT1}. Here, we give an idea of the proof of this convergence since it is essentially the same as the proof of \cite[(5.3)]{CSZ19b}.

Indeed, we may regard it as a Riemann-sum approximation for some function due to the following approximation:\begin{align*}
&q_n(\phi_N,x)\sim \int_{\R^2}\phi(y)p_{\frac{n}{N}}\left(\frac{x}{\sqrt{N}}-y\right)\dd y \\
&\sum_{y\in \Z^2}q_n(x,y)p_\e\left(\frac{y}{\sqrt{N}}-z\right)\psi(z)\sim p_{\frac{n}{N}+\e}\left(z-\frac{x}{\sqrt{N}}\right)\psi(z)\\
&q_{n}(x,y)\sim \frac{1}{N}p_{\frac{n}{N}}\left(\frac{y-x}{\sqrt{N}}\right)\\
&U_N\left(n,z\right)\sim \frac{4\pi}{N^2}G_\vartheta\left(\frac{n}{N},\frac{z}{\sqrt{N}}\right).
\end{align*}

The approximation of
$\widetilde{\Theta}^{(N)}
(\bm m,\bm n,\bm x,\bm y)$
produces the factor
\begin{align*}
N^{-\sum_{E\in\{A,B,C,D\}}
(r^E(\mathtt w)-1)}
=
N^{-(2k-4)},
\end{align*}
where for each $\mathtt w=([i_1,j_1],\dots,[i_k,j_k])
\in\mathcal P_f$ and $E\in\{A,B,C,D\}$, we define the number of appearances of $E$ in $\mathtt{w}$ by
\begin{align*}
r^E(\mathtt w)
&:=
\#\left\{
1\leq r\leq|\mathtt w|:
E\in\{i_r,j_r\}
\right\} 
\end{align*}
and we used 
\begin{align*}
\sum_{E\in\{A,B,C,D\}}
r^E(\mathtt w)
=
2k,
\end{align*}
since each of the $k$ pairings contains exactly two labels.
On the other hand, the approximation of
\begin{align*}
\prod_{i=1}^{k}
U(n_i-m_i,y_i-x_i)
\end{align*}
produces the factor $N^{-2k}$.
Together with the prefactor $N^{-5}$ and the Riemann-sum
normalization, this yields the required overall scaling.

To prove the convergence, we first look at the second term of \eqref{eq:chaosexpansionT1-1} with the near diagonal sets are cut off, i.e. \begin{align*}
\left\{\begin{array}{l}\dis n_i-m_i\geq \delta N ,\quad (1\leq i\leq |\mathtt{w}|),\\ \dis m_{i+1}-n_i\geq \delta N, \quad(0\leq i \leq |\mathtt{w}|-1),\\u-n_{|\mathtt{w}|}>\delta N\end{array}\right\}
\end{align*} for some fixed $\delta>0$. It  approximates \eqref{eq:chaosexpansionT1-1} uniformly in $N$ since we know that the boundedness of the fourth moment of $\overline{\mathsf{Z}}^{\phi}_{N;s}(1)$ and the second moment of $\langle M^{N,\phi}(\psi) \rangle_t$. 

Also, we can find from Remark \ref{rem:4thmoments} that the second term of \eqref{eq:momentT1-1} is approximated by the one restricted by \begin{align*}
\left\{\begin{array}{l}v_i-u_i>\e,\quad (1\leq i\leq  |\mathtt{w}|),\\ u_{i+1}-v_i>\e,\quad ( 0\leq i\leq |\mathtt{w}|-1),\\ u-v_{|\mathtt{w}|}>\e\end{array}\right\}.
\end{align*}

We can use the same arguments as in  \cite[p.415]{CSZ19b}, the summations  or integrals with the restricted spatial variables $\bm{x}$ and $\bm{y}$ to the set \begin{align*}
&\left\{\begin{array}{l}|x_1|\leq M\sqrt{N}, |y_i-x_i|\leq M\sqrt{N}, \quad (1\leq i\leq |\mathtt{w}|), \\ |x_{i+1}-y_i|\leq M\sqrt{N},\quad  (0\leq i\leq |\mathtt{w}|-1)\end{array}\right\}\quad &\text{for \eqref{eq:momentT1-1}}\\
&\left\{\begin{array}{l}|x_1|\leq M, |y_i-x_i|\leq M, \quad (1\leq i\leq |\mathtt{w}|),\\  |x_{i+1}-y_i|\leq M, \quad (0\leq i\leq |\mathtt{w}|-1)\end{array}\right\}\quad &\text{for \eqref{eq:chaosexpansionT1-1}}
\end{align*} 
for large $M>0$ are negligible as $N\to \infty$ and then $M\to \infty$. \toybox

\begin{proof}[Proof of Lemma \ref{lem:QuadProcEst}]
Lemma \ref{lem:4thmomentlimitTT}-\ref{lem:4thmomentlimit11} yield that \begin{align*}
&\lim_{N\to\infty}\IE\left[\left(\left\langle M^{N,\phi}(\psi)\right\rangle_{t}^{(\e)}-\left\langle M^{N,\phi}(\psi)\right\rangle_{Nt}\right)^2\right]\\
&=\frac{(4\pi)^2}{(-\log \e)^2} \eqref{eq:momentTT1}-2\frac{4\pi}{-\log \e}\eqref{eq:momentT1-1}+\eqref{eq:moments111}.
\end{align*}

Thus, it is enough to see \begin{align*}
\frac{(4\pi^2)^2}{(-\log \e)^2} \eqref{eq:momentTT1}-2\frac{4\pi^2}{-\log \e}\eqref{eq:momentT1-1}+\eqref{eq:moments111}\to 0
\end{align*}
as $\e\to 0$.

For \eqref{eq:momentTT1}, we first remark that  \begin{align*}
&\frac{1}{-\log \e}\left|\int_{v_{k^A}\vee v_{k^B}<\sigma<t}\dd \sigma \int_{\R^2}\dd z p_{\sigma-v_{k^A}+\e}\left(y_{k^A},z\right)p_{\sigma-v_{k^B}+\e}\left(y_{k^B},z\right)\psi(z)^2\right|\\
&\leq \frac{1}{-\log \e}\|\psi\|_\infty \int_{v_{k^A}\vee v_{k^B}<\sigma<t}\dd \sigma p_{2\sigma-v_{k^A}-v_{k^B}+2\e}(0)\leq \frac{\|\psi\|^2_\infty}{2\pi}C_t
\end{align*}
for $0<v_{k^A}\vee v_{k^B}<t$ and some $C_t>0$. Moreover, we find that \begin{align*}
&\frac{1}{-\log \e}\int_{v_{k^A}\vee v_{k^B}<\sigma<t}\dd \sigma \int_{\R^2}\dd z 
p_{\sigma-v_{k^A}+\e}\left(y_{k^A},z\right)p_{\sigma-v_{k^B}+\e}\left(y^B_{k^B},z\right)\psi(z)^2\\
&\to \begin{cases}
\frac{1}{4\pi}\psi(v_{k^A})^2 \quad &\text{if }v_{k^A}=v_{k^B} \text{ and }y_{k^A}=y_{k^B} \Leftrightarrow \text{$k^A=k^B$}\\
0&\text{otherwise}.
\end{cases}
\end{align*}
by Lemma \ref{lem:psiepsilon}. 
It does hold for $C,D$.

 In particular, they converge to $\frac{1}{4\pi}\psi(z_A)^2$ and $\frac{1}{4\pi}\psi(z_C)^2$ if and only if $k^A=k^B$ and $k^C=k^D$ ($\Leftrightarrow$ $(k^A,k^B,k^C,k^D)=(k-1,k-1,k,k)$ or $(k,k,k-1,k-1)$).
 
The dominated convergence theorem yields \begin{align*}
&\frac{1}{(-\log \e)^2}\eqref{eq:momentTT1}\\
&\to \sum_{k\geq 2}\sum_{\bsm \mathtt{w}\in {\mathcal{P}}_f\\ |\mathtt{w}|=k\\ [i_{k-1},j_{k-1}]=[A,B],[i_{k},j_{k}]=[C,D]\esm}(4\pi)^{k-2}\idotsint\limits_{0<u_1<v_1<\dots<u_k<v_k<t}\dd \bm{u}\dd \bm{v}\int_{(\bm{R}^2)^{2k}}\dd \bm{x}\dd \bm{y}\notag\\
&\hspace{4em}\prod_{i=1}^k G_\vartheta(v_i-u_i,y_i-x_i)\Theta_{\mathtt{w}} (\bm{u},\bm{v},\bm{x},\bm{y})\psi(y_{k-1})^2\psi(y_k)^2\\
&+\sum_{k\geq 2}\sum_{\bsm \mathtt{w}\in {\mathcal{P}}_f\\ |\mathtt{w}|=k\\ [i_{k-1},j_{k-1}]=[C,D],[i_{k-1},j_{k-1}]=[A,B]\esm}(4\pi)^{k-2}\idotsint\limits_{0<u_1<v_1<\dots<u_k<v_k<t}\dd \bm{u}\dd \bm{v}\int_{(\bm{R}^2)^{2k}}\dd \bm{x}\dd \bm{y}\notag\\
&\hspace{4em}\prod_{i=1}^k G_\vartheta(v_i-u_i,y_i-x_i)\Theta_{\mathtt{w}} (\bm{u},\bm{v},\bm{x},\bm{y})\psi(y_{k-1})^2\psi(y_k)^2\\
&=\frac{1}{(4\pi)^2}\eqref{eq:moments111}
\end{align*}

Applying the same argument to \eqref{eq:momentT1-1},  we obtain that \begin{align*}
&\frac{1}{(-\log \e)}\eqref{eq:momentT1-1}\\
&\to \sum_{k\geq 2}\sum_{\bsm \mathtt{w}\in {\mathcal{P}}_f\\ |\mathtt{w}|=k\\ [i_{k-1},j_{k-1}]=[A,B],[i_{k},j_{k}]=[C,D]\esm}(4\pi)^{k-1}\idotsint\limits_{0<u_1<v_1<\dots<u_k<v_k<t}\dd \bm{u}\dd \bm{v}\int_{(\bm{R}^2)^{2k}}\dd \bm{x}\dd \bm{y}\notag\\
&\hspace{4em}\prod_{i=1}^k G_\vartheta(v_i-u_i,y_i-x_i)\Theta_\mathtt{w} (\bm{u},\bm{v},\bm{x},\bm{y})\psi(y_{k-1})^2\psi(y_k)^2\\
&+\sum_{k\geq 2}\sum_{\bsm \mathtt{w}\in {\mathcal{P}}_f\\ |\mathtt{w}|=k\\ [i_{k-1},j_{k-1}]=[C,D],[i_{k},j_{k}]=[A,B]\esm}(4\pi)^{k-1}\idotsint\limits_{0<u_1<v_1<\dots<u_k<v_k<t}\dd \bm{u}\dd \bm{v}\int_{(\bm{R}^2)^{2k}}\dd \bm{x}\dd \bm{y}\notag\\
&\hspace{4em}\prod_{i=1}^k G_\vartheta(v_i-u_i,y_i-x_i)\Theta_\mathtt{w} (\bm{u},\bm{v},\bm{x},\bm{y})\psi(y_{k-1})^2\psi(y_k)^2\\
&=\frac{1}{4\pi}\eqref{eq:moments111}.
\end{align*}

\end{proof}

\begin{proof}[Proof of \eqref{eq:QuadGenlim} in Theorem \ref{thm:martinalemeasure}]
It is enough to show that \begin{align*}
\lim_{\e\to 0}\IE\left[\left(\frac{4\pi}{-\log \e}\int_0^t\int_{\R^2}\mathscr{Z}^{\vartheta,\phi}_s(p_\e(\cdot-z))^2\psi(z)^2\dd z\dd s-\left\langle \mathscr{M}^{\vartheta,\phi}(\psi)\right\rangle_t\right)^2\right]=0,
\end{align*}
and it follows once we modify the proof of Theorem \ref{thm:SHFMP} by using  the following lemma.
\end{proof}

\begin{lemma}\label{lem:psiepsilon}
Let $\psi\in \mathcal{B}_b(\R^2)$ and $T>0$. Then, for each  $0<t<T$, there exists $C_{T,\psi}$ such that \begin{align}
&\sup_{0<\e<\frac{1}{2}}\left|\frac{-1}{\log \e}\int_{\R^2}\dd z\int_{0}^t\dd s\frac{1}{s+\e} p_{\frac{s+2\e}{2}}(z-x)\psi(z)\right|\leq C_{T,\psi}\label{eq:e0conv1}
\intertext{for each $x\in \R^2$ and }
&\lim_{\e\to 0}\frac{-1}{\log \e}\int_{\R^2}\dd z\int_0^t\dd s\frac{1}{s+\e} p_{\frac{s+\e}{2}}(z-x)\psi(z)=\psi(x)\label{eq:e0conv2}
\end{align}
for a.e.~$x$ for $0<t\leq T$.
\end{lemma}
\begin{proof}
It is straightforward to see that\begin{align*}
&\left|\frac{-1}{\log \e}\int_{\R^2}\dd z\int_{0}^t\dd s\frac{1}{s+\e} p_{\frac{s+\e}{2}}(z-x)\psi(z)\right|\\
&\leq \|\psi\|_\infty \frac{-1}{\log \e}\left[\log(t+\e)-\log (\e)\right]\leq C_T\|\psi\|_\infty
\end{align*}
for some $C_T>0$.

Also, \begin{align*}
&\int_{\R^2}\dd z\int_0^t\dd s\frac{1}{s+\e} p_{\frac{s+\e}{2}}(z-x)\psi(z)\\
&=\int_{\R^2}\dd z\psi(z) \frac{1}{\pi |z-x|^2}\left(\exp\left(-\frac{|z-x|^2}{t+\e}\right)-\exp\left(-\frac{|z-x|^2}{\e}\right)\right)\\
&=\int_{0}^\infty \dd r\int_{[0,2\pi]}\dd \theta \psi\left(x+r(\cos \theta,\sin\theta)\right) \frac{1}{\pi r}\left(\exp\left(-\frac{r^2}{t+\e}\right)-\exp\left(-\frac{r^2}{\e}\right)\right).
\end{align*}
Therefore, we obtain from  l'H\^{o}pital's rule that \begin{align*}
&\lim_{\e\to 0}\frac{1}{-\log \e}\int_{\R^2}\dd z\int_0^t\dd s\frac{1}{s+\e} p_{\frac{s+\e}{2}}(z-x)\psi(z)\\
&=-\lim_{\e\to 0}\e \int_{0}^\infty \dd r\int_{[0,2\pi]}\dd \theta \psi\left(x+r(\cos \theta,\sin\theta)\right) \frac{1}{\pi r}\frac{\dd}{\dd \e}\left(\exp\left(-\frac{r^2}{t+\e}\right)-\exp\left(-\frac{r^2}{\e}\right)\right)
\end{align*}
if the limit in the right-hand side exists for each $t>0$.

It follows from the Lebesgue  differentiation theorem that it is equal to \begin{align*}
\lim_{\e\to 0} \int_{0}^\infty \dd r\int_{[0,2\pi]}\dd \theta \psi\left(x+r(\cos \theta,\sin\theta)\right) \frac{r}{\pi \e}\exp\left(-\frac{r^2}{\e}\right)=\lim_{\e\to 0}\int_{\R^2}p_{\e}(z-x)\psi(z)\dd z=\psi(x)
\end{align*}
a.e.~$x$. 

Therefore, \eqref{eq:e0conv2} holds for $\psi\in \mathcal{B}_b(\R^2)$. 
\end{proof}

\begin{remark}\label{rem:compactsupp}
In Theorem \ref{thm:SHFMP}, the quadratic variation is approximated
using the heat-kernel mollifier
$p_\e(\cdot-x)$. The same quadratic variation can be recovered using
a compactly supported mollifier of the form $f_\e(\cdot-x)
=
\frac1\e
f\left(
\frac{\cdot-x}{\sqrt\e}
\right),
$
where $f\in C_c^+(\R^2)$ and $\int f=1$.

Indeed, we can modify the proof by using the following lemma instead of Lemma \ref{lem:psiepsilon}.
\end{remark}
\begin{lemma}\label{lem:approxfun} Let $f\in C_c^+(\mathbb R^2)$ satisfy
\[
\int_{\mathbb R^2}f(y)\dd y=1.
\]
Let $\psi\in C_b^2(\R^2)$,  and $T>0$. Then, for each $0<t<T$, there exists $C_{T,f,\psi}$ such that 
 \begin{align}
&\sup_{0<\e<\frac{1}{2}}\left|\frac{-1}{\log \e}\int_{\R^2}\dd z\int_0^t\dd s\int_{\R^2\times \R^2}\hspace{-1.5em}\dd w\dd w' p_s(w-x)p_s(w'-x)\frac{1}{\e}f\left(\frac{w-z}{\sqrt{\e}}\right)\frac{1}{\e}f\left(\frac{w'-z}{\sqrt{\e}}\right)\psi(z)\right|\label{eq:approxfun1}\\
&\leq C_{T,f,\psi}\notag
\end{align}
 and \begin{align}
&\lim_{\e\to 0}\frac{-1}{\log \e}\int_{\R^2}\dd z\int_0^t\dd s\int_{\R^2\times \R^2}\hspace{-1.5em}\dd w\dd w' p_s(w-x)p_s(w'-x)\frac{1}{\e}f\left(\frac{w-z}{\sqrt{\e}}\right)\frac{1}{\e}f\left(\frac{w'-z}{\sqrt{\e}}\right)\psi(z)\label{eq:approxfun2}
\\
&=\frac{\psi(x)}{4\pi}\notag
\end{align}
for each $x\in \R^2$.
\end{lemma}
\begin{proof}
\eqref{eq:approxfun1} follows from \eqref{eq:e0conv1} and there exists a constant $C>0$ such that $f(x)\leq Cp_1(x)$ for any $x\in \R^2$.

For \eqref{eq:approxfun2}, we see that \begin{align*}
&\int_{\R^2}\dd z\int_0^t\dd s\int_{\R^2\times \R^2}\dd w\dd w' p_s(w-x)p_s(w'-x)\frac{1}{\e}f\left(\frac{w-z}{\sqrt{\e}}\right)\frac{1}{\e}f\left(\frac{w'-z}{\sqrt{\e}}\right)\psi(z)\\
&=\int_{\R^2}\dd z\int_0^t\dd s\int_{\R^2\times \R^2}\dd y\dd y' p_s(z+\sqrt{\e}y-x)p_s(z+\sqrt{\e}y'-x)f\left(y\right)f\left(y'\right)\psi(z)\\
&=\e\int_{\R^2}\dd \widetilde{z}\int_0^t\dd s\int_{\R^2\times \R^2}\dd y\dd y' p_s(\sqrt{\e}(y+\widetilde{z}))p_s(\sqrt{\e}(y'+\widetilde{z}))f\left(y\right)f\left(y'\right)\psi(x+\sqrt{\e}\widetilde{z})\\
&=\int_{\R^2}\dd \widetilde{z}\int_0^\frac{t}{\e}\dd u\int_{\R^2\times \R^2}\dd y\dd y' p_u(y+\widetilde{z})p_u(y'+\widetilde{z})f\left(y\right)f\left(y'\right)\psi(x+\sqrt{\e}\widetilde{z})\\
&=\int_{\R^2}\dd \widetilde{z}\int_0^\frac{t}{\e}\dd u\int_{\R^2\times \R^2}\dd y\dd y' p_u(y+\widetilde{z})p_u(y'+\widetilde{z})f\left(y\right)f\left(y'\right)\left(\psi(x)+O(\sqrt{\e}|\widetilde{z}|)\right),
\end{align*}
where the error term is bounded in absolute value by $\|\nabla\psi\|_\infty\sqrt\e\,|\widetilde z|$.
Moreover, \begin{align*}
&\int_{\R^2}\dd \widetilde{z}\int_0^\frac{t}{\e}\dd u\int_{\R^2\times \R^2}\dd y\dd y' p_u(y+\widetilde{z})p_u(y'+\widetilde{z})f\left(y\right)f\left(y'\right)\psi(x)\\
&=\psi(x)\int_0^\frac{t}{\e}\dd u\int_{\R^2\times \R^2}\dd y\dd y' p_{2u}(y-y')f\left(y\right)f\left(y'\right)\\
&=\frac{\psi(x)}{4\pi}\int_{\R^2\times \R^2}\dd y\dd y' f\left(y\right)f\left(y'\right)\int_{\frac{\e|y-y'|^2}{4t}}^\infty \frac{e^{-u}}{u}\dd u.
\end{align*}
Hence, after division by $-\log\e$, it converges to $\psi(x)/(4\pi)$ by l'H\^opital's rule.

Moreover, the contribution of the error term is bounded by
\begin{align*}
&\int_{(\R^2)^3}
\dd\widetilde z\dd y\dd y'\,
|\widetilde z|\,
p_u(y+\widetilde z)
p_u(y'+\widetilde z)
f(y)f(y')\leq
C_f\left(
\mathbf 1_{\{0<u\leq1\}}
+
u^{-1/2}\mathbf 1_{\{u>1\}}
\right).
\end{align*}
Indeed, for $0<u\leq1$, the bound follows from
$\|p_u*f\|_\infty\leq\|f\|_\infty$ and the compact support of $f$.
For $u\geq1$, the Cauchy--Schwarz inequality and
$\int_{\R^2}f(y)\dd y=1$ give
\begin{align*}
\left(
\int_{\R^2}
p_u(y+\widetilde z)f(y)\dd y
\right)^2
&\leq
\int_{\R^2}
p_u(y+\widetilde z)^2f(y)\dd y.
\end{align*}
Consequently, after the change of variables
$w=y+\widetilde z$,
\begin{align*}
&\int_{(\R^2)^3}
|\widetilde z|\,
p_u(y+\widetilde z)
p_u(y'+\widetilde z)
f(y)f(y')
\dd\widetilde z\dd y\dd y'\\
&\leq
\int_{\R^2}|w|p_u(w)^2\dd w+
\left(
\int_{\R^2}|y|f(y)\dd y
\right)
\int_{\R^2}p_u(w)^2\dd w\\
&\leq
C_f\left(u^{-1/2}+u^{-1}\right)
\leq
C_fu^{-1/2}.
\end{align*}

This proves \eqref{eq:approxfun2}.

\end{proof}

\end{appendix}

\textbf{Acknowledgements} This work was supported by JSPS KAKENHI Grant Numbers JP22K03351 and
JP23K22399. The author is grateful to Francesco Caravenna for helpful discussions
concerning the spaces $\mathscr C^{s,b}$ and to Nikos Zygouras for
useful comments on the manuscript.
The author used Microsoft Copilot as an assistive tool for proofreading
the manuscript and checking its English presentation, notation, and
internal consistency. The author critically reviewed the suggestions,
independently verified the mathematical content, and assumes full
responsibility for the final version.


\begin{thebibliography}{CSZ19b}

\bibitem[Alm05]{Alm05}
Alexandre Almeida.
\newblock Wavelet bases in generalized {B}esov spaces.
\newblock {\em J. Math. Anal. Appl.}, 304(1):198--211, 2005.

\bibitem[Bur73]{Bur73}
D.~L. Burkholder.
\newblock Distribution function inequalities for martingales.
\newblock {\em Ann. Probability}, 1:19--42, 1973.

\bibitem[Che26]{Che26}
Yu-Ting Chen.
\newblock Martingale problem of the two-dimensional stochastic heat equation at criticality.
\newblock arxiv:2504.21791, 2026.

\bibitem[CSZ]{CSZ24}
Francesco Caravenna, Rongfeng Sun, and Nikos Zygouras.
\newblock In {\em Statistical Mechanics and Stochastic PDEs}, volume 2401 of {\em Lecture Notes in Mathematics}. Springer.
\newblock arxiv:2412.10311.

\bibitem[CSZ19a]{CSZ19a}
Francesco Caravenna, Rongfeng Sun, and Nikos Zygouras.
\newblock The {D}ickman subordinator, renewal theorems, and disordered systems.
\newblock {\em Electron. J. Probab.}, 24:Paper No. 101, 40, 2019.

\bibitem[CSZ19b]{CSZ19b}
Francesco Caravenna, Rongfeng Sun, and Nikos Zygouras.
\newblock On the moments of the {$(2+1)$}-dimensional directed polymer and stochastic heat equation in the critical window.
\newblock {\em Comm. Math. Phys.}, 372(2):385--440, 2019.

\bibitem[CSZ23]{CSZ23}
Francesco Caravenna, Rongfeng Sun, and Nikos Zygouras.
\newblock The critical 2d {S}tochastic {H}eat {F}low.
\newblock {\em Invent. Math.}, 233(1):325--460, 2023.

\bibitem[CSZ25]{CSZ25}
Francesco Caravenna, Rongfeng Sun, and Nikos Zygouras.
\newblock Singularity and regularity of the critical 2d stochastic heat flow.
\newblock arXiv:2504.06128, 2025.

\bibitem[Daw93]{Daw93}
Donald~A. Dawson.
\newblock Measure-valued {M}arkov processes.
\newblock In {\em \'Ecole d'\'Et\'e{} de {P}robabilit\'es de {S}aint-{F}lour {XXI}---1991}, volume 1541 of {\em Lecture Notes in Math.}, pages 1--260. Springer, Berlin, 1993.

\bibitem[DF94]{DF94}
Donald~A. Dawson and Klaus Fleischmann.
\newblock A super-{B}rownian motion with a single point catalyst.
\newblock {\em Stochastic Process. Appl.}, 49(1):3--40, 1994.

\bibitem[FL06]{FL04}
Walter Farkas and Hans-Gerd Leopold.
\newblock Characterisations of function spaces of generalised smoothness.
\newblock {\em Ann. Mat. Pura Appl. (4)}, 185(1):1--62, 2006.

\bibitem[FM17]{FM17}
Marco Furlan and Jean-Christophe Mourrat.
\newblock A tightness criterion for random fields, with application to the {I}sing model.
\newblock {\em Electron. J. Probab.}, 22:Paper No. 97, 29, 2017.

\bibitem[FV10]{FV10}
Peter~K. Friz and Nicolas~B. Victoir.
\newblock {\em Multidimensional stochastic processes as rough paths}, volume 120 of {\em Cambridge Studies in Advanced Mathematics}.
\newblock Cambridge University Press, Cambridge, 2010.
\newblock Theory and applications.

\bibitem[GN25]{GN25}
Shirshendu Ganguly and Kyeongsik Nam.
\newblock Sharp moment and upper tail asymptotics for the critical $2d$ stochastic heat flow, 2025.
\newblock arxiv:2507.22029.

\bibitem[GN26]{GN26}
Shirshendu Ganguly and Kyeongsik Nam.
\newblock Logarithmic intermittency of the critical 2d shf, 2026.
\newblock arxiv:2608.12270.

\bibitem[GQT21]{GQT21}
Yu~Gu, Jeremy Quastel, and Li-Cheng Tsai.
\newblock Moments of the 2{D} {SHE} at criticality.
\newblock {\em Probab. Math. Phys.}, 2(1):179--219, 2021.

\bibitem[Hai14]{Hai14}
M.~Hairer.
\newblock A theory of regularity structures.
\newblock {\em Invent. Math.}, 198(2):269--504, 2014.

\bibitem[JS87]{JS87}
Jean Jacod and Albert~N Shiryaev.
\newblock {\em Limit theorems for stochastic processes}.
\newblock Springer, 1987.

\bibitem[Kah85]{Kah85}
Jean-Pierre Kahane.
\newblock Sur le chaos multiplicatif.
\newblock {\em Ann. Sci. Math. Qu\'ebec}, 9(2):105--150, 1985.

\bibitem[Kal17]{Kal17}
Olav Kallenberg.
\newblock {\em Random measures, theory and applications}, volume~77 of {\em Probability Theory and Stochastic Modelling}.
\newblock Springer, Cham, 2017.

\bibitem[LL10]{LL10}
Gregory~F. Lawler and Vlada Limic.
\newblock {\em Random walk: a modern introduction}, volume 123 of {\em Cambridge Studies in Advanced Mathematics}.
\newblock Cambridge University Press, Cambridge, 2010.

\bibitem[LX26]{LX26}
Ziyang Liu and Zuodi Xie.
\newblock Moments comparison inequalities for critical 2d stochastic heat flow, polymers and gmc, 2026.
\newblock arxiv:2609.12820.

\bibitem[LZ26]{LZ26}
Ziyang Liu and Nikos Zygouras.
\newblock On the {M}oments of the {M}ass of {S}hrinking {B}alls {U}nder the {C}ritical 2{$d$} {S}tochastic {H}eat {F}low.
\newblock {\em Comm. Math. Phys.}, 407(6):Paper No. 132, 2026.

\bibitem[Mit83]{Mit83}
Itaru Mitoma.
\newblock Tightness of probabilities on {$C([0,1];{\cal S}\sp{\prime} )$}\ and {$D([0,1];{\cal S}\sp{\prime} )$}.
\newblock {\em Ann. Probab.}, 11(4):989--999, 1983.

\bibitem[Mou01]{Mou01}
Susana Moura.
\newblock Function spaces of generalised smoothness.
\newblock {\em Dissertationes Math. (Rozprawy Mat.)}, 398:88, 2001.

\bibitem[Per02]{Per02}
Edwin Perkins.
\newblock Dawson-{W}atanabe superprocesses and measure-valued diffusions.
\newblock In {\em Lectures on probability theory and statistics ({S}aint-{F}lour, 1999)}, volume 1781 of {\em Lecture Notes in Math.}, pages 125--324. Springer, Berlin, 2002.

\bibitem[RV14]{RV14}
R\'emi Rhodes and Vincent Vargas.
\newblock Gaussian multiplicative chaos and applications: a review.
\newblock {\em Probab. Surv.}, 11:315--392, 2014.

\bibitem[Spi76]{Spi76}
Frank Spitzer.
\newblock {\em Principles of random walk}, volume Vol. 34 of {\em Graduate Texts in Mathematics}.
\newblock Springer-Verlag, New York-Heidelberg, second edition, 1976.

\bibitem[Tao11]{Tao11}
Terence Tao.
\newblock {\em An introduction to measure theory}, volume 126 of {\em Graduate Studies in Mathematics}.
\newblock American Mathematical Society, Providence, RI, 2011.

\bibitem[Tsa24]{Tsa24}
Li-Cheng Tsai.
\newblock Stochastic heat flow by moments, 2024.
\newblock arxiv:2410.14657.

\bibitem[Wal86]{Wal86}
John~B. Walsh.
\newblock An introduction to stochastic partial differential equations.
\newblock In {\em \'Ecole d'\'et\'e{} de probabilit\'es de {S}aint-{F}lour, {XIV}---1984}, volume 1180 of {\em Lecture Notes in Math.}, pages 265--439. Springer, Berlin, 1986.

\end{thebibliography}
\end{document}